  \documentclass[final,1p,times]{elsarticle}

\usepackage{amsmath}
\usepackage{amssymb} 
\usepackage{algorithm}
\usepackage{algorithmic}
\usepackage{multirow}
\usepackage{pstool}
\usepackage{multicol}
\allowdisplaybreaks
\allowdisplaybreaks

\extrafloats{100}

\newcommand{\Var}{\operatorname{Var}}

\newcommand{\bE}{{\mathbb{E}}}

\newcommand{\bR}{{\mathbb{R}}}
\newcommand{\bN}{{\mathbb{N}}}

\newcommand{\cC}{\mathcal{C}}

\newcommand{\cF}{\mathcal{F}}
\newcommand{\cH}{\mathcal{H}}
\newcommand{\cJ}{\mathcal{J}}

\newcommand{\cN}{\mathcal{N}}

\newcommand{\beq}{\begin{equation}}
\newcommand{\eeq}{\end{equation}}
\newcommand{\ba}{\begin{array}}
	\newcommand{\ea}{\end{array}}

\usepackage{framed,multirow}

\usepackage{latexsym}

\usepackage{url}
\usepackage{xcolor}
\definecolor{newcolor}{rgb}{.8,.349,.1}

\journal{Journal of Computational Physics}

\begin{document}

\begin{frontmatter}



\title{Optimal design of acoustic metamaterial cloaks under uncertainty}

 \author[label1]{Peng Chen}
 \author[label2]{Michael R. Haberman}
 \author[label1,label2,label3]{Omar Ghattas}
 \address[label1]{Oden Institute for Computational Engineering \& Sciences, The University of Texas at Austin, Austin, TX 78712 ({peng@oden.utexas.edu})}
 \address[label2]{Walker Department of Mechanical Engineering, The University of Texas at Austin, Austin, TX 78712 ({haberman@utexas.edu})}
 \address[label3]{Department of Geological Sciences, The University of Texas at Austin, Austin, TX 78712 ({omar@oden.utexas.edu})}
\fntext[label1]{
This research was partially funded by the Department of Energy, Office of Science, Office of Advanced Scientific Computing Research, Mathematical Multifaceted Integrated Capability Centers (MMICCS) program under award DE-SC0019303; the Simons Foundation under award 560651; the Air Force Office of Scientific Research, Computational Mathematics program under award FA9550- 17-1-0190; and the National Science Foundation, Division of Advanced Cyberinfrastructure under award ACI-1550593.
}

\begin{abstract}
In this work, we consider the problem of optimal design of an acoustic
cloak under uncertainty and develop scalable approximation and
optimization methods to solve this problem. The design variable is
taken as an infinite-dimensional spatially-varying field that
represents the material property, while an additive
infinite-dimensional random field represents the variability of the
material property or the manufacturing error. Discretization of this
optimal design problem results in high-dimensional design variables
and uncertain parameters. To solve this problem, we develop a
computational approach based on a Taylor approximation and an
approximate Newton method for optimization, which is based on a
Hessian derived at the mean of the random field. We show our approach
is scalable with respect to the dimension of both the design variables
and uncertain parameters, in the sense that the necessary number of
acoustic wave propagations is essentially independent of these
dimensions, for numerical experiments with up to one million design
variables and half a million uncertain parameters. We demonstrate
that, using our computational approach, an optimal design of the
acoustic cloak that is robust to material uncertainty is achieved in a
tractable manner. The optimal design under uncertainty problem is
posed and solved for the classical circular obstacle surrounded by a
ring-shaped cloaking region, subjected to both a single-direction
single-frequency incident wave and multiple-direction
multiple-frequency incident waves. Finally, we apply the method to a
deterministic large-scale optimal cloaking problem with complex
geometry, to demonstrate that the approximate Newton method's Hessian
computation is viable for large, complex problems.

\end{abstract}




\end{frontmatter}



\section{Introduction}
\label{sec:introduction}

Research on acoustic and elastic metamaterials is a product of a unique combination of technological advances that have been made over the last three decades to achieve extraordinary redirection, absorption, or amplification of acoustic or elastic wave disturbances by designing the sub-wavelength structure of the medium through which mechanical disturbances propagate \cite{HabermanGuild16, CummerChristensenAlu16, MaSheng16}. Of specific interest here is the field of acoustic and elastic metamaterials that make use of the convergence of novel concepts in physics with advances in technology and computational methods, primarily the field of additive manufacturing (AM) \cite{Mueller08,Matthews15,WegenerLinden10} and access to robust computational tools \cite{AlexanderianPetraStadlerEtAl17,AlexanderianGloorGhattas15,ChenVillaGhattas19,BashirWillcoxGhattasEtAl08}. The simultaneous rapid maturation of AM and computational methods allows researchers to rapidly simulate, build, and test elaborate structures for acoustic wave manipulation that follow from rigorous mathematical predictions such as transformation acoustics (TA) \cite{Cummer13, ChenChan10}. TA is a mathematical approach that uses coordinate transformations to map the physical space to a different space of interest using a one-to-one map between the two domains \cite{Cummer13, CummerPopaSchurigEtAl08, Norris08, ChenChan10}. The mathematical map is then used to determine the material properties in the region of the transformed fields that produce the same effect in the physical space. For example, mapping the acoustic field in a finite volume surrounding a small scatterer to that of a shell surrounding a larger object allows one to determine the material properties within the shell that produce a cloak capable of rerouting acoustic waves around the large object. This approach provides a forward model for the determination of the material properties required to generate an acoustic cloak using TA.

 However, the true research challenge is to define material microstructure that generates effective material properties that meet the prescription provided by TA for the frequencies of interest using existing materials and manufacturing methods. Coordinate transformation and its application to the manipulation of electromagnetic waves preceded the application of coordinate transformations to acoustic waves. Indeed the concept of transformation acoustics was initially facilitated by utilizing the direct analogy that exists in two dimensions between Maxwell's equations and the equations describing acoustic wave propagation \cite{Cummer13, ChenChan10}. Slight differences in the coordinate transformation were found for 3D geometries, arising from the fundamental differences of 3D wave propagation for transverse (electromagnetic) and longitudinal (acoustic) waves \cite{CummerPopaSchurigEtAl08}. Although highly anisotropic, such effective fluids could be theoretically realized using alternating layered structures with ordinary (isotropic) fluid-like properties \cite{TorrentSanchez-Dehesa08} or frequency-dependent waveguide designs \cite{ZhangXiaFang11}. However, coordinate transformation of elastic materials, in which compressional and transverse shear waves co-exist, require a far more exotic metamaterial with a fundamentally new type of microstructure: pentamode (PM) materials. PM materials are defined as materials whose stiffness tensors have only one non-null eigenvalue out of a possible six \cite{Norris08}. In other words, these materials have five deformation modes that can be imposed without storing energy in the material \cite{Norris08} and they can therefore be thought of as quasi-fluids. While the transformation acoustics provides an exact analytical solution for the material properties required to create a perfect cloak, it suffers from several serious drawbacks. The primary problem with this approach is that while it prescribes the material properties required to achieve cloaking, it cannot define the sub-wavelength structure that will generate the required properties. In this sense, coordinate transformation methods are simply analytical methods to solve forward problems and shed no light on how the behavior of interest can be generated. The vast majority of acoustic cloaking research has therefore relied on physical insight and researcher creativity to find material structures that generate the required material properties to achieve cloaking of an object.

A powerful technique to achieve cloaking can be accomplished using a plasmonic cloak, also known as a scattering cancellation (SC) cloak, which is a non-resonant means of eliminating the field scattered from an object, thereby hiding it from detection. This was originally applied and demonstrated for electromagnetic waves \cite{AluEngheta05, RainwaterKerkhoffMelinEtAl12} using plasmonic materials to achieve the necessary cloaking layer properties. The SC approach to cloaking was subsequently shown to be an effective means of cloaking acoustic waves \cite{GuildHabermanAlu11, GuildAluHaberman11}. Unlike cloaks developed using a coordinate transformation approach \cite{Leonhardt06, Norris08}, only the scattered field in the surrounding medium is eliminated, and therefore this solution does not limit the incident wave from interacting with the object. As a result, there is no restriction on the use of isotropic materials to create a plasmonic cloak, and it may be used to suppress the scattering from sensors  \cite{AluEngheta09, AluEngheta10, GuildAluHaberman14}. In previous work, the composition of an SC cloak for cylindrical or spherical objects was found by minimizing the total scattering cross-section of the object and cloak by varying the number, radius, and material properties of layers surrounding the object to be cloaked \cite{AluEngheta05, GuildHabermanAlu11, GuildAluHaberman11,GuildHabermanAlu12}. This approach was later extended to the design of cloaks for non-spherical objects and collections of objects \cite{GuildHicksHabermanEtAl15}. Further, the SC method is well-suited for numerical approaches to determining material property distributions required to achieve cloaking. It has been applied to design three-dimensional cloaks with unidirectional performance \cite{SanchisGarcia-ChocanoLlopis-PontiverosEtAl13} and two-dimensional cloaks that exploit B\'{e}zier scatterers in the cloaking region to minimize the scattered field \cite{LuSanchisWenEtAl18}. A similar computational approach employs a gradient-based optimization algorithm to minimize the total scattering cross section (TCSC) of a collection of rigid or elastic cylinders surrounded cylindrical scatterers that collectively act as a unidirectional cloak \cite{AmirkulovaNorris20}. Similar work by Andkj\ae r and Sigmund employed topology optimization to design a cloak that used a small number of scatterers in the cloaking region to conceal a circular region in two-dimensional space from detection via airborne sound \cite{AndkjerSigmund13}. Each of these contributions employ numerical optimization to determine the geometry and properties in the cloaking region. However, these works and many others have only paid cursory attention to the influence that variation in material properties or geometry may have on cloaking performance. Further, to the authors' knowledge, there has been no effort to study how cloaking design may change when variability is accounted for in the design. Given that fabrication of cloaks must consider real-world variation in as-built material properties or achievable levels of manufacturing precision, addressing this problem is central to the creation of reliable acoustics cloaks.

In the optimal design of acoustic cloaks, uncertainties may arise from
various sources, including material property variability and flaws or
deviations introduced by the manufacturing process. It is therefore
important to take uncertainties into account in order to design a
robust cloak that can cancel the scattered wave as much as possible
under different realizations of the uncertainty. For this purpose we
consider the problem of optimal design of an acoustic cloak under
uncertainty. While our methodology can be applied more generally, the
case considered here is that of time-harmonic acoustic wave
propagation and scattering from an impenetrable obstacle. The wave
motion in the background medium and the cloak is described by the
Helmholtz equation with varying wavenumber, i.e., a spatially-varying
sound speed in the cloaking region. We model the sound speed in the
cloak as a perturbation of the sound speed in the host homogeneous
medium by an exponential factor, which is taken as an
infinite-dimensional spatially-varying design variable field.  The
uncertain parameter is modeled as a Gaussian random field that is
additive to the design variable supported in the same cloak
region. The objective for the optimal design is to minimize the
scattered wave outside the obstacle and cloak region, for which we
take a suitable norm of this quantity as the design objective. Since
the design objective depends on the uncertain parameter through the
Helmholtz equation, it is also an uncertain or random function. To
account for this uncertainty in the optimal design, we consider both
the mean and the variance of the design objective and minimize an
objective functional including a weighted combination of the two. To
promote the sparsity of the design material, we add a weighted
$L^1$-norm of the design variable as a penalty to the objective
functional.

The optimal design under uncertainty problem presented above leads to
a random partial differential equation (PDE)-constrained optimization
problem, which after appropriate discretization results in
high-dimensional uncertain parameters and optimization
variables. Solution of this class of problems faces enormous
challenges, and has received increasing attention in recent years
\cite{BorziSchulzSchillingsEtAl10, SchillingsSchmidtSchulz11,
  HouLeeManouzi11,GunzburgerLeeLee11,RosseelWells12,
  KouriHeinkenschloosVanBloemenWaanders12,TieslerKirbyXiuEtAl12,ChenQuarteroniRozza13,LassilaManzoniQuarteroniEtAl13,
  ChenQuarteroni14, KouriHeinkenschlossRidzalEtAl14, KunothSchwab13,
  NgWillcox14, ChenQuarteroniRozza16, KunothSchwab16, KouriSurowiec16,
  BennerOnwuntaStoll16, AlexanderianPetraStadlerEtAl17,
  AliUllmannHinze17, ChenVillaGhattas19, ChenGhattas20a}.
One prominent challenge is the evaluation of the high-dimensional
integral involved in the mean and variance of the design objective. A
straightforward approach is to use Monte Carlo integration, which
amounts to the sample average approximation (SAA) method, which has a
convergence rate ($O(M^{-1/2})$ with $M$ samples) that does not depend
on the parameter dimension. Nevertheless, its convergence is often too
slow, so a large number of samples is required to achieve a certain
required accuracy. Since one PDE has to be solved for each sample,
this leads to an optimization constrained by a large number of PDEs,
and thus this method is usually computationally prohibitive. As an
alternative, rapidly-convergent methods such as stochastic Galerkin
and collocation have been applied \cite{HouLeeManouzi11,
  GunzburgerLeeLee11, TieslerKirbyXiuEtAl12, RosseelWells12,
  KouriHeinkenschloosVanBloemenWaanders12,
  LassilaManzoniQuarteroniEtAl13, ChenQuarteroniRozza13,
  ChenQuarteroni14, KunothSchwab13, KunothSchwab16}. However, they
often face the curse of dimensionality, i.e., the computational
complexity grows exponentially with respect to the uncertain parameter
dimension, which prevents their use for problems with high-dimensional
uncertain parameters.

Another challenge is that discretization of the
design variable field leads to a high-dimensional optimization
problem. A simple steepest descent based method will require far too many
optimization iterations to converge, while a Newton method may
converge rapidly but require the computation of the Hessian of the
objective functional acting in given directions, which is often too
complex for sophisticated approximations of the objective functional
as we employ here. In this work, we propose a computational approach
based on a Taylor approximation for the evaluation of the
high-dimensional integral in the objective functional and an
approximate Newton method for the high-dimensional optimization
problem. We employ the Taylor approximation based optimization
strategy proposed in \cite{AlexanderianPetraStadlerEtAl17,
  ChenVillaGhattas19, ChenGhattas20a}, by which we approximate the
design objective by its (quadratic) Taylor expansion with respect to
the uncertain parameter and compute the trace of the preconditioned
Hessian resulting from this approximation by a randomized singular
value decomposition (SVD) algorithm. The computational complexity
measured in terms of the number of PDE solves depends only on
the---often small and dimension independent---number of dominant
eigenvalues of the preconditioned Hessian, and not on the nominal
large uncertain parameter dimension. Thus this approximation is
scalable with respect to the parameter dimension. To solve the
high-dimensional optimization problem, we propose an approximate
Newton method in which the Hessian of the objective functional based on
the quadratic Taylor approximation, which is too complex to compute,
is approximated by that of the deterministic objective functional,
i.e., one that is evaluated at the mean of the random variable.
Provided the uncertainty is not too large, e.g., the noise-to-signal
ratio or the ratio between the magnitude of the uncertain parameter
and that of the design variable is less than $20\%$ in our
application, the deterministic Hessian provides a good approximation
of the true Hessian, thus leading to an optimization method that
is (effectively) scalable with respect to the optimization variable
dimension. 

We apply the proposed computational approach to the optimal design
under uncertainty of the acoustic cloak in several different
settings. A classical circular obstacle surrounded by a ring-shaped
cloak is used to demonstrate the efficacy of optimal design under
uncertainty and our scalable computational approach.  First, we
consider a deterministic approximation of the objective functional,
which results in a deterministic optimal design problem. In this
setting, the scattered field is efficiently eliminated by the
optimization. Second, we compare this design with the optimal design
under uncertainty and show that the latter achieves a significant
reduction in variability of the scattered field relative to the
deterministic optimal design. Third, to demonstrate the scalability of
the Taylor approximation and the approximate Hessian-based Newton
optimization methods, we solve the optimal design problem for a
sequence of refined finite element discretizations with dimension up
to half a million for the uncertain parameters and one million for the
design variables. Scalability with respect to the dimensions of the
uncertain parameters and the design variables is demonstrated by
dimension-independence of (1) the convergence of the optimizer, (2)
the spectral decay of the eigenvalues of the preconditioned Hessian of
the design objective with respect to the uncertain parameters, and (3)
the accuracy of the Taylor approximation. Fourth, we extend the
optimal design problem with single direction, single frequency
incident wave to one with multiple directions and multiple
frequencies, and demonstrate the efficacy of the acoustic
cloak. Finally, we consider a more complex geometry representative of
a stealth aircraft, for which we also obtain an effective acoustic
cloak.

The rest of the paper is organized as follows: In Section
\ref{sec:design}, we formulate the optimal design of the acoustic
cloak under uncertainty problem, including the governing Helmholtz PDE
constraint, the uncertain parameters and design variables, and the
formulation of the mean-variance objective functional and sparsifying
penalty term. Section \ref{sec:approximation} presents the
approximation methods of the mean-variance functional, including the
deterministic approximation, the sample average approximation, and the
Taylor approximation with randomized SVD computation of the resulting
trace. The optimization method is presented in Section
\ref{sec:optimization}, in which the computation of the gradient and
(approximate) Hessian of the objective functional with respect to the
design variables, as well as the approximate Newton method itself, are
derived. Several numerical experiments for the optimal design of an
acoustic cloak are presented in Section \ref{sec:numerics}, which is
followed by conclusions in Section \ref{sec:conclusion}.

\section{Problem formulation}
\label{sec:design}
In this section, we formulate the problem of optimal design of an acoustic cloak under uncertainty. The forward problem consists of time-harmonic acoustic wave scattering in an inhomogeneous medium described by the Helmholtz equation, in a region truncated by perfectly matched layer. We describe the representation of the design variables and uncertain parameters, the mean-variance objective functional, and the formulation of the optimal design under uncertainty problem.

\subsection{Acoustic wave scattering}
The time-harmonic acoustic wave scattering of an incident wave in a host medium from an impenetrable obstacle surrounded by an inhomogeneous metameterial medium is governed by the following Helmholtz equation \cite{ColtonKress12}:
\begin{subequations}\label{eq:Helmholtz} 
	\beq\label{eq:HelmholtzA}
	\Delta u + k^2 u  =  (k_0^2 - k^2)u^{\text{inc}} \quad \text{ in } \bR^d\setminus {D}_o,
	\eeq
	\beq\label{eq:HelmholtzB}
	\nabla u \cdot n = - \nabla  u^{\text{inc}} \cdot n \quad \text{ on }  \partial D_o, 
	\eeq
	\beq\label{eq:HelmholtzC}
	\lim_{r \to \infty } r^{(d-1)/2} \left(\frac{\partial u}{\partial r } - i k u \right)  = 0,
	\eeq
\end{subequations}
where $\bR^d$ is the physical space of dimension $d =2, 3$, $D_o\subset \bR^d$ is the region of the obstacle with boundary $\partial D_o$. $u^{\text{inc}} $ is the incident wave given by $u^{\text{inc}} = e^{i k_0 x \cdot b}$ in direction $b\in \bR^d$ with the complex unit $i = \sqrt{-1}$; $u$ is the scattered wave; the total wave is given by $u^t = u + u^{\text{inc}}$. In addition, $k_0$ is the wavenumber in the background medium given by the positive constant $k_0 =\omega/c_0$ with frequency $\omega$ and constant speed of sound $c_0$ in the host medium, while $k(x) = \omega/c(x)$ is a spatially-varying wavenumber in the inhomogeneous medium. $c(x)$ denotes the speed of sound at $x\in \bR^d$ in the inhomogeneous medium. 
A \emph{sound-hard} boundary condition is imposed on the boundary $\partial D_o$ in \eqref{eq:HelmholtzB} for the impenetrable obstacle, where $n$ denotes the outward unit normal vector along $\partial D_o$. Eq.\ \eqref{eq:HelmholtzC} is the \emph{Sommerfeld radiation condition} that guarantees that the scattered wave is outgoing, which is realized by a \emph{perfectly matched layer} (PML) condition \cite{TurkelYefet98}. $r(x) = |x|$ denotes the distance from $x$ to the origin. 

\begin{figure}[!htbp]
	\begin{center}
		\includegraphics[scale=0.5]{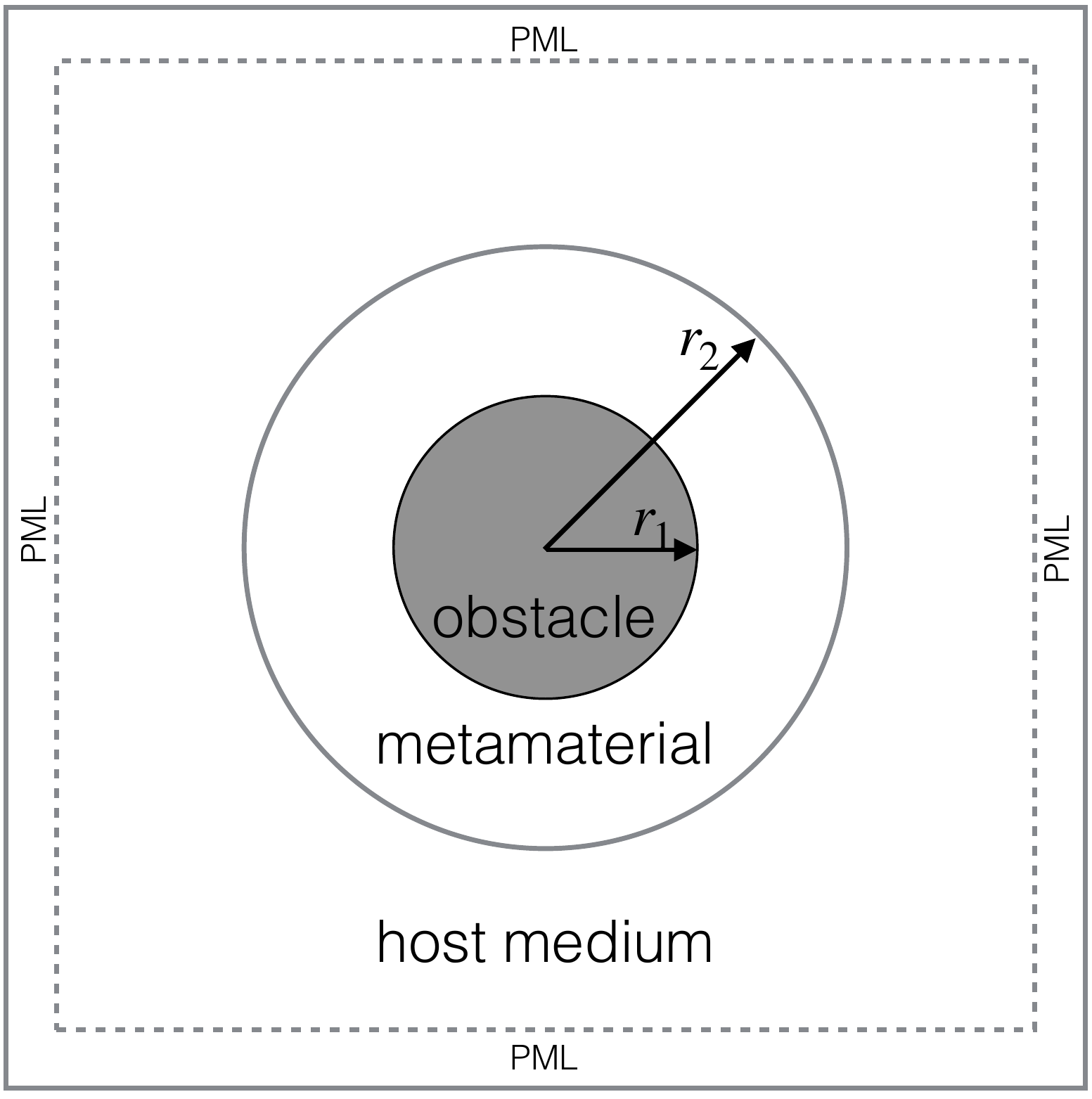}
	\end{center}
	\caption{Sketch of the domain for acoustic wave scattering in an inhomogeneous medium. }\label{fig:SketchAcousticWave}
\end{figure}

\subsection{Weak formulation with PML condition}
To solve the problem numerically, we consider a bounded and square computational domain $D \subset \bR^d \setminus {D}_o$ that includes the inhomogeneous metamaterial medium as shown in Fig.\ \ref{fig:SketchAcousticWave}, and use a PML condition \cite{TurkelYefet98} on its outgoing boundary to prevent reflection of the scattered wave as imposed by \eqref{eq:HelmholtzC}. In two dimensions, i.e., $d = 2$, the PML condition leads to \cite{TurkelYefet98}
\beq\label{eq:HelmholtzPML}
\partial_{x_1} \left(\frac{s_{x_2}}{s_{x_1}} \partial_{x_1} u\right) +  \partial_{x_2} \left(\frac{s_{x_1}}{s_{x_2}} {\partial_{x_2} u}\right) + k^2 s_{x_1} s_{x_2} u =  (k_0^2-k^2)u^{\text{inc}} \quad \text{in } D,
\eeq
where 
\beq
s_{x_1} = 1 + \frac{\sigma_{x_1}}{ik}, \qquad s_{x_2} = 1 + \frac{\sigma_{x_2}}{ik},
\eeq 
where $k(x) = k_0 \sqrt{n(x)} = \omega/c(x)$, $\sigma_{x_1}$ and $\sigma_{x_2}$ are real valued continuous functions in the PML region $D_p$, which depend only on the physical coordinate $x_1$ and $x_2$, respectively. Outside the PML region, i.e., $D \setminus D_{p}$, we have $\sigma_{x_1} = 0$ and $\sigma_{x_2} = 0$, so that \eqref{eq:HelmholtzPML} becomes the same equation as \eqref{eq:HelmholtzA}. The wave function $u$ is complex valued, which can be written as $u = u_1 + i u_2$ with the real and imaginary parts $u_1$ and $u_2$, respectively. Then \eqref{eq:HelmholtzPML} can be written as a set of two equations in $D$ 
with real coefficients as 
\beq\label{eq:HelmholtzSplit}
\begin{split}
	\partial_{x_1} \left( a_1 \partial_{x_1} u_1 - a_2 \partial_{x_1} u_2\right)+ \partial_{x_2} \left( a_3 \partial_{x_2} u_1 - a_4 \partial_{x_2} u_2\right) + b_1 u_1 - b_2 u_2 =  (k_0^2-k^2)u^{\text{inc}}_1\\
	\partial_{x_1} \left( a_1 \partial_{x_1} u_2 + a_2 \partial_{x_1} u_1\right)+ \partial_{x_2} \left( a_3 \partial_{x_2} u_2 + a_4 \partial_{x_2} u_1\right) + b_1 u_2 + b_2 u_1 =  (k_0^2-k^2)u^{\text{inc}}_2 \\
\end{split}
\eeq
where $u^{\text{inc}}_1 = \cos(k_0x\cdot b)$, $u^{\text{inc}}_1 = \sin(k_0x\cdot b)$; the coefficients are given by
\beq
a_1 = \frac{k^2 + \sigma_{x_1}\sigma_{x_2}}{k^2 + \sigma_{x_1}^2}, \; a_2 = \frac{k(\sigma_{x_1} - \sigma_{x_2})}{k^2+\sigma_{x_1}^2}, \;  a_3 = \frac{k^2 + \sigma_{x_1}\sigma_{x_2}}{k^2 + \sigma_{x_2}^2}, \; a_4 = \frac{k(\sigma_{x_2} - \sigma_{x_1})}{k^2+\sigma_{x_2}^2},
\eeq
and 
\beq
b_1 = k^2 - \sigma_{x_1} \sigma_{x_2}, \; b_2 = -k(\sigma_{x_1} + \sigma_{x_2}).
\eeq

To state the weak formulation of the equations \eqref{eq:HelmholtzSplit}, we introduce the following notation.
Let $L^2(D)$ denote the Hilbert space of square-integrable functions in $D$, and let $H^1(D) := \{v \in L^2(D), |\nabla v| \in L^2(D)\}$, $V = H^1(D) \times H^1(D)$.
Then the weak formulation of \eqref{eq:HelmholtzSplit} can be written as: find $u = (u_1, u_2) \in V$, such that 
\beq\label{eq:HelmholtzWeak}
A(u, v) = F(v), \quad \forall v = (v_1, v_2) \in V,
\eeq
where the bilinear form $A: V \times V \to \bR$ is given by
\beq
\begin{split}\label{eq:bilinearHelmholtz}
	A(w, v) & = \int_D (a_1 \partial_{x_1} w_1 - a_2 \partial_{x_1} w_2) \partial_{x_1} v_1 + (a_3 \partial_{x_2} w_1 - a_4 \partial_{x_2} w_2) \partial_{x_2} v_1  \; dx \\
	& + \int_D  (a_1 \partial_{x_1} w_2 + a_2 \partial_{x_1} w_1) \partial_{x_1} v_2 + (a_3 \partial_{x_2} w_2 + a_4 \partial_{x_2} w_1) \partial_{x_2} v_2  \; dx \\
	& - \int_D (b_1 w_1 - b_2 w_2) v_1 + (b_1 w_2 + b_2 w_1) v_2 \; dx 
\end{split} 
\eeq
and the linear form $F: V \to \bR$ is given by 
\beq\label{eq:linearHelmholtz}
F(v) = \int_D  (k_0^2-k^2)u^{\text{inc}}_1 v_1 +  (k_0^2-k^2)u^{\text{inc}}_2 v_2 dx - \int_{\partial D_o} \nabla u^{\text{inc}}_1 \cdot n v_1 +\nabla u^{\text{inc}}_2 \cdot n v_2 ds .
\eeq


\subsection{Uncertain parameter and design variable}
To manufacture the acoustic cloak, additive manufacturing (AM) offers significant promise since it allows the fabrication of complex parts that cannot be readily created using traditional techniques \cite{Mueller12, Matthews15}. Of specific interest here is the potential to construct materials with spatially graded material properties by adjusting process control variables within the build volume. However, this capability is not perfect and errors can be introduced at each manufacturing step. Further, each fabrication approach has some level of uncertainty in the as-built material properties, whose contribution to the final manufacturing accuracy is unclear \cite{PaulAnandGerner14, PintoArrietaAndiaEtAl15}. In this work, we consider an aggregated uncertainty and incorporate it in the sound speed in the cloak, which is represented by 
\beq\label{eq:soundSpeed}
c(x, \omega) = c_0 e^{\zeta(x, \varpi) - \tau(x)}, \quad \forall x \in D_m, a.e.\ \varpi \in \Omega 
\eeq
Here, $\tau$ is the spatially-varying deterministic design variable field of the cloak in the design region $D_m$, which exists in a separable Banach space $Z$. At every $x \in D_m$, $\zeta(x,\cdot)$ is a real valued random variable defined in the probability space $(\Omega, \cF, P)$, with the sample space $\Omega$, a set of events $\cF$, and the probability function $P: \cF \to [0, 1]$. To accommodate the spatial correlation of the random variables $\zeta(x,\cdot)$ at different $x \in D_m$, we consider one of the most popular random fields---\emph{Gaussian random fields} with probability measure $\mu = \cN(\bar{\zeta}, \cC)$ defined in a Hilbert space $X$ with dual $X'$, where $\bar{\zeta} \in X$ is the mean and $\cC$ is the covariance operator that can be viewed as an integral operator with suitable covariance kernel. A very general covariance kernel widely used in spatial statistics, geostatistics, machine learning, etc., is the Mat\'ern covariance, which leads to the Gaussian field $\zeta$ as a solution of the stochastic fractional PDE \cite{LindgrenRueLindstroem11} with homogeneous Neumann boundary condition
\beq\label{eq:fsPDE} 
\begin{split}
	( - \gamma \Delta + \delta I)^{\alpha/2} (\zeta - \bar{\zeta}) & = W \quad \text{ in } D_m, \\
	\quad \nabla \zeta \cdot n & = 0 \quad \text{ on } \partial D_m,
\end{split}
\eeq 
where $W$ represents the spatial Gaussian white noise with unit variance, $I$, $\nabla$, and $\Delta$ are the identity, gradient, and Laplace operators, and $n$ is the outward unit normal vector along $\partial D_m$. Thus  $\cC = (-
\gamma \Delta + \delta I )^{-\alpha}$, with $\alpha > d/2$ controlling the
\emph{regularity}, $\delta$ and $\gamma$ controlling the \emph{variance}, and $\gamma/\delta$ controlling the \emph{correlation length} of $\zeta$. Moreover, $\cC$ is self-adjoint, positive, and of trace class.
Therefore, sampling $\zeta$ involves solution of the elliptic stochastic PDE \eqref{eq:fsPDE}. Generalizations of the stochastic PDE \eqref{eq:fsPDE} may be used to model non-stationary, non-isotropic, complex random fields \cite{LindgrenRueLindstroem11}.

\subsection{Optimal design of acoustic cloak}
Recall that $D_o$ and $D_m$ denote the regions of the obstacle and the metamaterial cloak surrounding the obstacle, respectively; let $D_h = D \setminus (D_o \cup D_m)$ denote the host medium where we can observe the scattered wave. Our goal is to minimize the scattered wave in $D_h$ so that the obstacle becomes ``invisible", i.e., no wave scattering observed outside of the obstacle and its cloak. To achieve this, we define the design objective as 
\beq\label{eq:objective}
Q(u) = \int_{D_h} \left(|u_1|^2 + |u_2|^2\right) dx,
\eeq
which is the scattered wave amplitude measured in the $L^2$-norm. The design objective $Q$ is random and depends on the random (field) variable $\zeta$ through the random scattered wave $u$ as a solution of \eqref{eq:HelmholtzWeak}, where we write $u = u(\zeta, \tau)$ to indicate that the scattered wave depends on the random variable $\zeta$ and the design variable $\tau$. To quantify the randomness of $Q$, we use a mean-variance measure and consider the following objective functional to be minimized 
\beq\label{eq:objectiveFunctional}
J(\tau) = \bE[Q](\tau) + \beta_V \Var[Q](\tau) + \beta_P P(\tau),
\eeq
where the mean and variance of $Q$ are given by 
\beq\label{eq:meanvariance}
\bE[Q](\tau) = \int_X Q(u(\zeta, \tau)) d\mu \text{ and }  \Var[Q](\tau) = \int_X (Q(u(\zeta, \tau)) - \bE[Q](\tau))^2 d\mu,
\eeq
where the integration is taken with respect to the Gaussian measure $\mu = \cN(\bar{\zeta}, \cC)$ in $X$. $P(\tau)$ is a penalty term on the deterministic control $\tau \in Z$. To promote the sparsity of the material in the cloak, we consider an $L^1$-norm for $\tau$, i.e.,
\beq\label{eq:penalization}
P(\tau) = \int_{D_m} |\tau(x)| dx \approx \int_{D_m} (\tau^2(x) + \varepsilon)^{1/2} dx,
\eeq
where we use the approximate form with a small $\varepsilon > 0$ to make $P(\tau)$ differentiable with respect to $\tau$ and thus facilitate gradient based optimization. Further, $\beta_V > 0$ and $\beta_P > 0$ in \eqref{eq:objectiveFunctional} are scalar parameters that weight the importance of the variance and penalty with respect to the mean. The problem of the optimal design of the acoustic cloak under uncertainty is finally formulated as the PDE-constrained stochastic optimization problem
\beq\label{eq:optimization}
\min_{\tau \in  Z} J(\tau), \; \text{ subject to } \eqref{eq:HelmholtzWeak}.
\eeq

\subsection{Optimal design with multiple directions and frequencies}
In the above formulation of the optimal design of the acoustic cloak under uncertainty, we consider only one direction $b$ and one frequency $\omega$ for the incident wave $u^{\text{inc}} = e^{i k_0 x \cdot b}$ where $k_0 = \omega/c_0$. In this section, we extend the formulation to incident waves with multiple directions and multiple frequencies. For notational clarity, for direction $b_i$ and/or frequency $\omega_i$, $i = 1, \dots, I$ for $I \in \bN$, we write the weak formulation \eqref{eq:HelmholtzWeak} as: find $u^i = (u^i_1, u^i_2) \in V$ such that  
\beq\label{eq:HelmholtzWeaki}
A_i(u^i, v^i) = F_i(v^i), \quad \forall v^i = (v^i_1, v^i_2) \in V,
\eeq
and write the design objective \eqref{eq:objective} as  
\beq
Q_i = Q(u^i) = \int_{D_h} \left(|u_1^i|^2 + |u_2^i|^2\right) dx.
\eeq
The objective functional \eqref{eq:objectiveFunctional} then becomes 
\beq\label{eq:objectiveFunctionali}
\cJ(\tau) = \sum_{i = 1}^I \big(\bE[Q_i](\tau) + \beta_V \Var[Q_i](\tau)\big) + \beta_P P(\tau),
\eeq
where the mean, variance, and penalty are given as in \eqref{eq:meanvariance} and \eqref{eq:penalization}. Therefore, the optimal design problem with multiple directions and frequencies becomes
\beq\label{eq:optimizationMultiple}
\min_{\tau \in Z} \cJ(\tau), \text{ subject to } \eqref{eq:HelmholtzWeaki}, \quad i = 1, \dots, I.
\eeq
Note that the approximation and optimization methods developed in the rest of the paper for the optimal design problem \eqref{eq:optimization} with single direction and frequency can be straightforwardly extended to the optimal design problem \eqref{eq:optimizationMultiple} with multiple directions and frequencies. For simplicity, we present methods for only the former case.

\section{Approximation of the mean-variance objective}
\label{sec:approximation}
In this section, we present three classes of approximation methods for the evaluation of the mean and variance in the objective functional: one is a deterministic approximation with the design objective evaluated only at the mean of the random variable $\zeta$, the second is
a classical sample average approximation, and the third is a quadratic Taylor approximation. For notational simplicity, whenever there is no ambiguity, we denote $Q(\zeta)$ for the random objective $Q(u(\zeta, \tau))$ at design $\tau \in Z$, and keep in mind that the dependence of $Q$ on $\zeta$ is implicit through $u$.

\subsection{Deterministic approximation}
\label{sec:deterministicApproximation}
\label{sec:deterministic}
In this approach, we evaluate the design objective at only one fixed sample, e.g., the mean $\bar{\zeta}$ of the random variable $\zeta$, so that the expectation and variance of the design objective are approximated as 
\beq
\bE[Q] \approx Q(\bar{\zeta}) \text{ and } \text{Var}[Q] \approx 0,
\eeq 
which leads to a deterministic optimization problem at $\bar{\zeta}$.

\subsection{Sample average approximation}
\label{sec:SAA}
Let $\zeta_m$, $m = 1, \dots, M$, denote i.i.d.\ random samples drawn from the Gaussian distribution $\cN(\bar{\zeta}, \cC)$, then the mean of $Q$ can be approximated by the average 
\beq
\bE[Q] \approx  \frac{1}{M} \sum_{m=1}^M Q(\zeta_m),
\eeq 
which is known as \emph{sample average approximation} or \emph{Monte Carlo approximation}. The variance can be approximated similarly by the average 
\beq
\Var[Q]  = \bE[Q^2] - (\bE[Q])^2 \approx \frac{1}{M} \sum_{m=1}^M Q^2(\zeta_m) - \left(\frac{1}{M} \sum_{m=1}^M Q(\zeta_m)
\right)^2.
\eeq 
We remark that to balance the approximation errors of the mean and variance, different numbers of i.i.d.\ random samples can be used for the mean and variance evaluation.

\subsection{Taylor approximation} 
\label{sec:Taylor}
Following the previous work \cite{AlexanderianPetraStadlerEtAl17, ChenVillaGhattas19}, we present a Taylor approximation for the design objective $Q$ and the closed form of the mean and variance based on the Taylor approximation. A formal functional Taylor approximation of the objective $Q$ at the mean $\bar{\zeta}$, truncated with $K$ terms, is written as 
\beq\label{eq:Taylor}
T_K Q(\zeta) = \sum_{k=0}^K \partial_\zeta^k Q(\bar{\zeta}) (\zeta - \bar{\zeta})^k,
\eeq
where we assume that $Q$ is $K$-th order Fr\'echet differentiable with respect to $\zeta$. The term $\partial_\zeta^k Q(\bar{\zeta}) (\zeta - \bar{\zeta})^k$ denotes the $k$-th order (tensor) derivative $\partial_\zeta^k Q(\bar{\zeta})$ at $\bar{\zeta}$ acting on $\zeta - \bar{\zeta}$ in each of the $k$ directions, $k = 1, \dots, K$. For $K = 2$, we can write the Taylor approximation \eqref{eq:Taylor} more explicitly as 
\beq\label{eq:TaylorQuadratic}
T_2 Q(\zeta) = \bar{Q} + \langle \bar{g}, \zeta - \bar{\zeta}\rangle + \frac{1}{2} \langle \bar{\cH} (\zeta - \bar{\zeta}), \zeta - \bar{\zeta}\rangle,
\eeq
where $\bar{Q} \in \bR$, $\bar{g} \in X'$, and $\bar{\cH}: X \to X'$ denote the objective and its gradient and Hessian with respect to $\zeta$, evaluated at the mean $\bar{\zeta}$, respectively, and $\langle \cdot , \cdot \rangle = \, _{X'} \langle \cdot , \cdot \rangle_X$ represents the duality pairing in $X'\times X$. Since $\zeta$ is a Gaussian field, the mean and variance of the Taylor approximation of the objective truncated at the quadratic term can be written explicitly as \cite{AlexanderianPetraStadlerEtAl17} 
\beq
\bE[T_2Q] = \bar{Q} + \frac{1}{2} \text{tr} (\cC \bar{\cH}) \quad \text{ and } \quad \Var[T_2Q] = \langle  \bar{g}, \cC \bar{g}\rangle + \frac{1}{2} \text{tr} ((\cC\bar{\cH})^2),
\eeq
where we recall that $\cC: X' \to X$ is the covariance of $\zeta$, and $\text{tr}(\cdot)$ denotes the trace, with 
\beq
\text{tr} (\cC \bar{\cH}) = \sum_{n\geq 1} \lambda_n\text{ and } \text{tr} ((\cC\bar{\cH})^2) = \sum_{n\geq 1} \lambda^2_n. 
\eeq
Here, $(\lambda_n)_{n\geq 1}$ are the eigenvalues of $\cC \bar{\cH}$, which are equivalent to the generalized eigenvalues of $(\bar{\cH}, \cC^{-1})$, i.e., in weak form we can write
\beq\label{eq:generalizedEigenProblem}
\langle \bar{\cH} \psi_n, \phi \rangle = \langle \lambda_n \cC^{-1} \psi_n, \phi \rangle \quad \forall \phi \in X, \quad n\geq 1,
\eeq
where $(\psi_n)_{n\geq 1} \in X$ are the generalized eigenfunctions that satisfy the $\cC^{-1}$-orthonormality condition 
\beq\label{eq:orthonormality}
\langle \cC^{-1} \psi_n, \psi_m \rangle = \delta_{mn}, \quad m\geq 1, n\geq 1.
\eeq

\subsubsection{Randomized algorithm}
It is intractable to solve the generalized eigenvalue problem \eqref{eq:generalizedEigenProblem} for all of the eigenvalues. In practice, these (absolute) eigenvalues decay rapidly as proven for the Hessians of some model problems and numerically demonstrated for many others \cite{BashirWillcoxGhattasEtAl08, FlathWilcoxAkcelikEtAl11,
	Bui-ThanhGhattas12a, Bui-ThanhGhattas13a, Bui-ThanhGhattas12,
	Bui-ThanhBursteddeGhattasEtAl12,
	Bui-ThanhGhattasMartinEtAl13,
	AlexanderianPetraStadlerEtAl16, AlexanderianPetraStadlerEtAl17,
	AlexanderianPetraStadlerEtAl14, CrestelAlexanderianStadlerEtAl17,
	PetraMartinStadlerEtAl14, IsaacPetraStadlerEtAl15,
	MartinWilcoxBursteddeEtAl12, Bui-ThanhGhattas15, ChenVillaGhattas17, ChenVillaGhattas19, ChenGhattas19a, ChenWuChenEtAl19a, ChenGhattas20}. Therefore, we can compute the dominant eigenvalues $\lambda_1, \dots, \lambda_N$, with $|\lambda_1| \geq \cdots \geq |\lambda_N| \geq \lambda_n$ for any $n > N$, and approximate the trace by 
\beq\label{eq:trace}
\text{tr} (\cC \bar{\cH}) \approx \sum_{n\geq 1}^N \lambda_n\text{ and } \text{tr} ((\cC\bar{\cH})^2) \approx \sum_{n\geq 1}^N \lambda^2_n. 
\eeq
To solve the generalized eigenvalue problem \eqref{eq:generalizedEigenProblem} for the $N$ dominant eigenvalues, we apply a randomized algorithm \cite{HalkoMartinssonTropp11, SaibabaLeeKitanidis16} in Algorithm \ref{alg:randomizedEigenSolver}, where  $H, C^{-1}$ of dimension $N_h \times N_h$ denote the discrete approximation of $\bar{\cH}$ and $\cC^{-1}$, e.g., by finite elements. Here, $N_h$ is the number of mesh degrees of freedom representing the discretized field $\zeta$.

\begin{algorithm}[!htb]
	\caption{Randomized algorithm for the generalized eigenvalue problem $(H, C^{-1})$}
	\label{alg:randomizedEigenSolver}
	\begin{algorithmic}
		\STATE{\textbf{Input: } the number of eigenpairs $N$, an oversampling factor $p \leq 10$.}
		\STATE{\textbf{Output: } $(\Lambda_N, \Psi_N)$ with $\Lambda_N = \text{diag}(\lambda_1, \dots, \lambda_N)$ and $\Psi_N = (\psi_1, \dots, \psi_N)$.}
		\STATE{1. Draw a Gaussian random matrix $\Omega \in \bR^{N_h \times (N+p)}$.}
		\STATE{2. Compute $Y = C (H \Omega)$.}
		\STATE{3. Compute $QR$ factorization $Y = QR$ such that $Q^\top C^{-1} Q = I_{N+p}$.}
		\STATE{4. Form $T = Q^\top H Q$ and compute eigendecomposition $T = S \Lambda S^\top$.}
		\STATE{5. Extract $\Lambda_N = \Lambda(1:N, 1:N)$ and $\Psi_N = QS_L$ with $S_N = S(:,1:N)$.}
	\end{algorithmic}
\end{algorithm}

We remark that the computational cost of Algorithm \ref{alg:randomizedEigenSolver} is dominated by the Hessian actions $H \Omega$ and $H Q$, as presented in the next section. These entail $2(N+p)$ forward and adjoint solutions of the Helmholtz equation. The remaining linear algebra in Algorithm \ref{alg:randomizedEigenSolver} is negligible in comparison.
The advantages of Algorithm \ref{alg:randomizedEigenSolver} are \cite{ChenVillaGhattas19, ChenGhattas19a}: (i) the error in the eigenvalues $\lambda_n$, $n = 1, \dots, N$, is bounded by the remaining ones $\lambda_n$, $n > N$, which is small if they decay rapidly;
(ii) the computational cost is dominated by $2(N+p)$ Hessian actions (the application of $C$ on a vector is inexpensive, e.g., it takes only $O(N_h)$ operations by a multigrid solver for $C$ discretized from an elliptic differential operator); (iii) it is scalable as the number of dominant eigenvalues $N$ typically does not depend on the mesh size $N_h$;
(iv) computing the Hessian actions $H\Omega$ and $H Q$ can be asynchronously parallelized.

\subsubsection{$\zeta$-gradient and $\zeta$-Hessian}
\label{sec:zeta-derivative}
The Taylor approximation along with the randomized eigensolver require the computation of the gradient of $Q(\zeta)$ with respect to the random parameter field (the ``$\zeta$-gradient") and the action of the Hessian of $Q(\zeta)$ (the ``$\zeta$-Hessian") in an arbitrary direction, both evaluated at the mean $\bar{\zeta}$. 
To do this, we employ a Lagrangian method as in \cite{ChenVillaGhattas19, ChenGhattas19a}. We begin by forming the Lagrangian 
\beq\label{eq:LagrangeG}
L(u,v; \zeta,\tau) = Q(u) + A(u,v;\zeta, \tau) - F(v;\zeta, \tau),
\eeq
where the bilinear form $A$ and the linear form $F$ defined in \eqref{eq:bilinearHelmholtz} and \eqref{eq:linearHelmholtz} depend on the random parameter and design variable $\zeta, \tau$ through the representation \eqref{eq:soundSpeed}.  The adjoint variable $v$ is a Lagrange multiplier for the forward Helmholtz equation \eqref{eq:HelmholtzWeak}. Then the state $u$ is obtained by setting the variation of the Lagrangian \eqref{eq:LagrangeG} with respect to the adjoint $v$ to zero to obtain the Helmholtz equation evaluated at $\bar{\zeta}$, i.e., find $u\in V$ such that
\beq\label{eq:state}
A(u, \tilde{v}; \bar{\zeta}, \tau) = F(\tilde{v}; \bar{\zeta}, \tau) \quad \forall \tilde{v} \in V,
\eeq
which is the same as \eqref{eq:HelmholtzWeak} evaluated at $\bar{\zeta}$. The adjoint variable $v$ is obtained by setting the variation of \eqref{eq:LagrangeG} with respect to the state $u$ to zero to obtain the adjoint Helmholtz equation evaluated at $\bar{\zeta}$, i.e., find $v\in V$ such that
\beq\label{eq:adjoint}
A(\tilde{u}, v; \bar{\zeta}, \tau) = - \langle \partial_u Q(u), \tilde{u}\rangle \quad \forall \tilde{u} \in V.
\eeq
Then the gradient of the design objective $Q$ with respect to the random variable $\zeta$ evaluated at $\bar{\zeta}$, acting in any direction $\tilde{\zeta} \in X$, is given by the variation of the Lagrangian with respect to $\zeta$, i.e.,
\beq\label{eq:gradient}
\langle \bar{g}, \tilde{\zeta} \rangle =  \langle \partial_\zeta L(u, v; \bar{\zeta}, \tau), \tilde{\zeta} \rangle = \langle \partial_\zeta A(u,v; \bar{\zeta}, \tau) - \partial_\zeta F(v; \bar{\zeta}, \tau), \tilde{\zeta} \rangle.
\eeq
Therefore, the computation of $\zeta$-gradient involves the solution of the Helmholtz equation \eqref{eq:state} for $u$ and the Helmholtz equation \eqref{eq:adjoint} for $v$.

To compute the Hessian of $Q$ at $\bar{\zeta}$ acting in a given direction $\hat{\zeta} \in X$, we form the second Lagrangian $L^H$ by adding the (weak formulation of the) forward and adjoint Helmholtz equations to the (directional) gradient to obtain  
\beq
\begin{split}
	L^H(u, v, \hat{u}, \hat{v}; \zeta, \hat{\zeta}, \tau) & = A(u, \hat{v}; {\zeta}, \tau) - F(\hat{v}; {\zeta}, \tau) \\
	& + A(\hat{u}, v; {\zeta}, \tau) + \langle \partial_u Q(u), \hat{u}\rangle \\
	& + \langle \partial_\zeta A(u,v; {\zeta}, \tau) - \partial_\zeta F(v; {\zeta}, \tau), \hat{\zeta} \rangle,
\end{split}
\eeq
where $\hat{v}, \hat{u}, \hat{\zeta}$ are the Lagrange multipliers for the forward Helmholtz equation \eqref{eq:state}, the adjoint Helmholtz equation \eqref{eq:adjoint}, and the gradient \eqref{eq:gradient}.
Proceeding as with the gradient derivation, we set the variation of $L^H$ with respect to $v$ and $u$ to obtain the incremental state variable $\hat{u}$ as the solution of the ``incremental forward Helmholtz equation" (evaluated at $\bar{\zeta}$)
\beq\label{eq:incrementalState}
A(\hat{u}, \tilde{v}; \bar{\zeta}, \tau) = - \langle \partial_{\zeta} A(u,\tilde{v}; \bar{\zeta}, \tau)  - \partial_{\zeta} F(\tilde{v}; \bar{\zeta}, \tau), \hat{\zeta} \rangle \quad \forall \tilde{v} \in V,
\eeq
and the incremental adjoint variable $\hat{v}$ as the solution of  the ``incremental adjoint Helmholtz equation" (evaluated at $\bar{\zeta}$)
\beq\label{eq:incrementalAdjoint}
A(\tilde{u}, \hat{v}; \bar{\zeta}, \tau) = -  \langle \partial_{uu} Q(u) \hat{u} , \tilde{u} \rangle - \langle \partial_{\zeta} A(\tilde{u},v; \bar{\zeta}, \tau), \hat{\zeta} \rangle \quad \forall \tilde{u} \in V.
\eeq
Finally, the Hessian action at $\bar{\zeta}$ in direction $\hat{\zeta}$, tested again $\tilde{\zeta}$, can be evaluated as 
\beq\label{eq:hessian}
\begin{split}
	\langle \bar{\cH} \hat{\zeta}, \tilde{\zeta} \rangle  = \langle \partial_\zeta L^H, \tilde{\zeta} \rangle & = \langle \partial_\zeta A(u, \hat{v}; \bar{\zeta}, \tau) - \partial_\zeta F(\hat{v}; \bar{\zeta}, \tau), \tilde{\zeta} \rangle\\
	& + \langle \partial_\zeta A(\hat{u}, v; \bar{\zeta}, \tau), \tilde{\zeta} \rangle\\
	& + \langle \partial_{\zeta\zeta} A(u,v; \bar{\zeta}, \tau) \hat{\zeta} - \partial_{\zeta\zeta} F(v; \bar{\zeta}, \tau)\hat{\zeta}, \tilde{\zeta} \rangle.
\end{split}
\eeq
Therefore, each Hessian action involves the solution of the incremental forward Helmholtz equation \eqref{eq:incrementalState} and the incremental adjoint Helmholtz equation \eqref{eq:incrementalAdjoint}. To compute the objective functional \eqref{eq:objectiveFunctional} with the quadratic Taylor approximation \eqref{eq:TaylorQuadratic} and the randomized algorithm for trace estimation (Algorithm \ref{alg:randomizedEigenSolver}), we need to solve one forward Helmholtz equation \eqref{eq:state}, one adjoint Helmholtz equation \eqref{eq:adjoint}, and $2(N+p)$ pairs of incremental forward and adjoint Helmholtz equations \eqref{eq:incrementalState} and \eqref{eq:incrementalAdjoint}. 

\section{Optimization}
\label{sec:optimization}
In the PDE-constrained optimization problem \eqref{eq:optimization}, the design variable field is a function over the cloaking region, and is thus high-dimensional after discretization by finite elements. To solve the resulting high-dimensional optimization problem, we propose an approximate Newton method with backtracking line search for globalization, where the Hessian of the objective functional with respect to the design variable, denoted as the $\tau$-Hessian, is approximated by the Hessian evaluated at the mean of the random field, while the gradient, denoted as the $\tau$-gradient, is computed accurately. The Newton system is solved inexactly in matrix-free fashion by a preconditioned conjugate gradient method.
In this section, we present the computation of the $\tau$-Hessian at the mean as well as the $\tau$-gradient of the objective functional \eqref{eq:objectiveFunctional} for both the sample average approximation of Section \ref{sec:SAA} and the Taylor approximation of Section \ref{sec:Taylor}.

\subsection{$\tau$-gradient and $\tau$-Hessian for the deterministic approximation}\label{sec:gradDeterministic}
Using the deterministic approximation of Section \ref{sec:deterministic}, we obtain the deterministic optimization problem:
\beq\label{eq:optimizationDeter}
\begin{split}
	& \min_{\tau \in Z} J_{\bar{\zeta}}(\tau) \quad \text{ where } J_{\bar{\zeta}}(\tau) = Q(u) + \beta_P P(\tau), \\
	& \text{subject to } A(u, v; \bar{\zeta}, \tau) = F(v; \bar{\zeta}, \tau) \quad \forall v \in V.
\end{split}
\eeq
To compute the gradient and Hessian of the objective functional with respect to the design variable $\tau$, we use a Lagrangian method akin to that presented in Section \ref{sec:zeta-derivative} for the gradient and Hessian of the design objective with respect to the random variable. Specifically, we first form the Lagrangian 
\beq
L_{\bar{\zeta}}(u, v; \bar{\zeta}, \tau) = Q(u) + \beta_P P(\tau) + A(u, v; \bar{\zeta}, \tau) - F(v; \bar{\zeta}, \tau).
\eeq
The state variable $u$ and the adjoint variable $v$ are obtained by setting the variation of this Lagrangian with respect to the adjoint $v$ and the state $u$ to zero and solving the forward and adjoint Helmholtz equations, which leads to the same problems as in \eqref{eq:state} and \eqref{eq:adjoint}.
The $\tau$-gradient (the Fr\'echet derivative of the objective in a direction $\tilde{\tau}$) is then given by 
\beq\label{eq:gradienttau}
\begin{split}
	\langle \nabla_\tau J_{\bar{\zeta}}(\tau), \tilde{\tau} \rangle &= \langle \partial_\tau L_{\bar{\zeta}}(u, v; \bar{\zeta}, \tau), \tilde{\tau} \rangle \\
	&= \langle \beta_P \nabla_\tau P(\tau) + \partial_\tau A(u, v; \bar{\zeta}, \tau) - \partial_\tau F(v; \bar{\zeta}, \tau), \tilde{\tau} \rangle.
\end{split}
\eeq
To compute the $\tau$-Hessian acting in a direction $\hat{\tau} \in Z$, we form the second Lagrangian
\beq
\begin{split}
	L^H_{\bar{\zeta}}(u, v, \hat{u}, \hat{v}; \bar{\zeta}, \tau, \hat{\tau}) & = A(u, \hat{v}; \bar{\zeta}, \tau) - F(\hat{v}; \bar{\zeta}, \tau) \\
	& + A(\hat{u}, v; \bar{\zeta}, \tau) + \langle \partial_u Q(u), \hat{u}\rangle \\
	& + \langle \beta_P \nabla_\tau P(\tau) + \partial_\tau A(u,v; \bar{\zeta}, \tau) - \partial_\tau F(v; \bar{\zeta}, \tau), \hat{\tau} \rangle,
\end{split}
\eeq
where $\hat{v}, \hat{u}, \hat{\tau}$ are the Lagrange multipliers for the forward Helmholtz equation \eqref{eq:state}, the adjoint Helmholtz equation \eqref{eq:adjoint}, and the gradient \eqref{eq:gradienttau}, respectively.
Once again, by setting the variation of $L^H_{\bar{\zeta}}$ with respect to $v$ and $u$ to zero, we obtain the incremental state variable $\hat{u}$ as the solution of the incremental forward Helmholtz equation 
\beq\label{eq:incrementalStatetau}
A(\hat{u}, \tilde{v}; \bar{\zeta}, \tau) = - \langle \partial_{\tau} A(u,\tilde{v}; \bar{\zeta}, \tau)  - \partial_{\tau} F(\tilde{v}; \bar{\zeta}, \tau), \hat{\tau} \rangle \quad \forall \tilde{v} \in V,
\eeq
and the incremental adjoint variable $\hat{v}$ as the solution of the incremental adjoint Helmholtz equation 
\beq\label{eq:incrementalAdjointtau}
A(\tilde{u}, \hat{v}; \bar{\zeta}, \tau) = -  \langle \partial_{uu} Q(u) \hat{u} , \tilde{u} \rangle - \langle \partial_{\tau} A(\tilde{u},v; \bar{\zeta}, \tau), \hat{\tau} \rangle \quad \forall \tilde{u} \in V.
\eeq
Then the $\tau$-Hessian action at $\tau$ in a direction $\hat{\tau}$, tested against $\tilde{\tau}$, can be evaluated as 
\beq\label{eq:hessiantau}
\begin{split}
	\langle \nabla_{\tau\tau} J_{\bar{\zeta}} \hat{\tau}, \tilde{\tau} \rangle & = \langle \partial_\tau L^H_{\bar{\zeta}}, \tilde{\tau} \rangle \\
	& = \langle \partial_\tau A(u, \hat{v}; \bar{\zeta}, \tau) - \partial_\tau F(\hat{v}; \bar{\zeta}, \tau), \tilde{\tau} \rangle\\
	& + \langle \partial_\tau A(\hat{u}, v; \bar{\zeta}, \tau), \tilde{\tau} \rangle\\
	& + \langle \beta_P \nabla_{\tau \tau} P(\tau) \hat{\tau}+ \partial_{\tau\tau} A(u,v; \bar{\zeta}, \tau) \hat{\tau} - \partial_{\tau\tau} F(v; \bar{\zeta}, \tau)\hat{\tau}, \tilde{\tau} \rangle.
\end{split}
\eeq
Therefore, at each $\tau$, after solving the forward Helmholtz equation \eqref{eq:state} and the adjoint Helmholtz equation \eqref{eq:adjoint}, to compute the $\tau$-Hessian action in each direction $\hat{\tau}$, we need to solve two PDEs---one incremental forward Helmholtz equation \eqref{eq:incrementalStatetau} and one incremental adjoint Helmholtz equation \eqref{eq:incrementalAdjointtau}. In Section \ref{sec:approximateHessianalg}, we derive how this capability for computing the action of the $\tau$-Hessian in an arbitrary direction can be used to solve the (approximate) Newton system by conjugate gradients.

\subsection{$\tau$-gradient for the sample average approximation}\label{sec:gradSAA}
With the sample average approximation (SAA), the optimization problem \eqref{eq:optimization} becomes 
\beq\label{eq:optimizationSAA}
\begin{split}
	&\min_{\tau \in Z} J_{\text{SAA}}(\tau)\\
	&\text{subject to } A(u_m, v; \zeta_m, \tau) = F(v; \zeta_m, \tau) \quad \forall v \in V, \quad m = 1, \dots, M,
\end{split}
\eeq
where $u_m = u(\zeta_m, \tau)$ represents the solution at $\zeta_m$ and $\tau$, and the SAA of the objective functional, $J_{\text{SAA}}(\tau)$, is given by 
\beq\label{eq:J_SAA}
J_{\text{SAA}}(\tau) = 
\frac{1}{M} \sum_{m=1}^M Q(u_m) + \frac{\beta_V}{M} \sum_{m=1}^M Q^2(u_m) - \beta_V \left(\frac{1}{M} \sum_{m=1}^M Q(u_m)\right)^2 + \beta_P P(\tau).
\eeq
To compute the $\tau$-gradient of $J_{\text{SAA}}$, we form the Lagrangian 
\beq
\begin{split}
	&L_{\text{SAA}}((u_m)_{m=1}^M, (v_m)_{m=1}^M; (\zeta_m)_{m=1}^M, \tau) \\
	&= J_{\text{SAA}}(\tau) + \sum_{m=1}^M A(u_m, v_m; \zeta_m, \tau) - F(v_m; \zeta_m, \tau),
\end{split}
\eeq
where $v_m$, $m = 1, \dots, M$, are the adjoint variables or the Lagrange multipliers. By setting the variation of the Lagrangian with respect to the state $u_m$ to zero for each $m = 1, \dots, M$, we obtain: find $v_m \in V$ such that 
\beq
\begin{split}\label{eq:adjointM}
	A(\tilde{u}, v_m; \zeta_m, \tau) = C_m \langle \partial_u Q(u_m), \tilde{u} \rangle \quad \forall \tilde{u} \in V, \quad m = 1, \dots, M,
\end{split}
\eeq 
where the constant $C_m$ is given by 
\beq
C_m = - \frac{1}{M} \left(1+ 2\beta_V Q(u_m) - 2\beta_V \left(\frac{1}{M} \sum_{m=1}^M Q(u_m)\right)\right), \quad m = 1, \dots, M.
\eeq
The $\tau$-gradient of $J_{\text{SAA}}(\tau)$ in \eqref{eq:J_SAA} can be computed as 
\beq
\begin{split}
	\nabla_\tau J_{\text{SAA}}(\tau) &= \partial_\tau L_{\text{SAA}}((u_m)_{m=1}^M, (v_m)_{m=1}^M; (\zeta_m)_{m=1}^M, \tau)  \\
	&= \beta_P \nabla_\tau P(\tau) + \sum_{m=1}^M \partial_\tau A(u_m, v_m; \zeta_m, \tau) - \partial_\tau F(v_m; \zeta_m, \tau).
\end{split}
\eeq
Hence, $M$ forward Helmholtz problems in \eqref{eq:optimizationSAA} are solved to compute $J_{\text{SAA}}(\tau)$, and $M$ adjoint problems \eqref{eq:adjointM} are solved to compute its $\tau$-gradient.

\subsection{$\tau$-gradient for the quadratic Taylor approximation}\label{sec:gradTaylor}
With the quadratic Taylor approximation of the design objective $T_2Q$, the objective functional \eqref{eq:objectiveFunctional} becomes 
\beq
J_{T_2}(\tau) = Q(u) + \frac{1}{2}\sum_{n=1}^N \lambda_n  + \beta_V \left(\langle \bar{g}, \cC\bar{g}\rangle  +  \frac{1}{2}\sum_{n=1}^N \lambda_n^2\right) + \beta_P P(\tau),
\eeq
where the $\tau$-gradient $\bar{g}$ is given by \eqref{eq:gradient}. Then the optimization problem \eqref{eq:optimization} reads 
\beq\label{eq:T2J}
\min_{\tau\in Z} J_{T_2}(\tau)
\eeq
subject to 
\beq\label{eq:T2constraints}
\begin{split}
	A(u, \tilde{v}; \bar{\zeta}, \tau) &= F(\tilde{v}; \bar{\zeta}, \tau) \; \forall \tilde{v} \in V, \\[4pt]
	A(\tilde{u}, v; \bar{\zeta}, \tau) &= - \langle \partial_u Q(u), \tilde{u}\rangle \; \forall \tilde{u} \in V, \\[4pt]
	A(\hat{u}_n, \tilde{v}; \bar{\zeta}, \tau) &= -  \langle \partial_{\zeta} A(u,\tilde{v}; \bar{\zeta}, \tau)  + \partial_{\zeta} F(\tilde{v}; \bar{\zeta}, \tau), \psi_n \rangle,  \; \forall \tilde{v} \in V, n = 1, \dots, N, \\[4pt]
	A(\tilde{u}, \hat{v}_n; \bar{\zeta}, \tau) &= -  \langle \partial_{uu} Q(u) \hat{u}_n , \tilde{u} \rangle -  \langle \partial_{\zeta} A(\tilde{u},v; \bar{\zeta}, \tau), \psi_n \rangle \; \forall \tilde{u} \in V, n = 1, \dots, N,\\[4pt]
	\langle \bar{\cH} \psi_n, \phi \rangle &= \langle \lambda_n \cC^{-1} \psi_n, \phi \rangle \; \forall \phi \in X,  n= 1, \dots, N,\\[4pt]
	\langle \cC^{-1} \psi_n, \psi_m \rangle &= 1, \; m, n = 1, \dots, N,
\end{split}
\eeq 
which correspond to the forward Helmholtz equation \eqref{eq:state}, the adjoint Helmholtz equation \eqref{eq:adjoint}, the incremental forward Helmholtz equation \eqref{eq:incrementalState} for $\hat{\zeta} = \psi_n$, $n = 1, \dots, N$, the incremental adjoint Helmholtz equation \eqref{eq:incrementalAdjoint} for $\hat{\zeta} = \psi_n$, $n = 1, \dots, N$, the generalized eigenvalue problem \eqref{eq:generalizedEigenProblem} for the eigenpairs $(\lambda_n, \psi_n)$, where the $\tau$-Hessian action $\bar{\cH}\psi_n$ is given by \eqref{eq:hessian}, $n = 1, \dots, N$, and the orthonormality condition \eqref{eq:orthonormality} for the eigenfunctions $\psi_n$, $n = 1, \dots, N$. As can be seen, the dominant cost of computing the objective functional $J_{T_2}(\tau)$ is $N$ pairs of (incremental) forward/adjont Helmholtz equations. This is in contrast with the $M$ forward Helmholtz equations which must be solved to compute the SAA objective $J_{\text{SAA}}$.

To compute the $\tau$-gradient of the approximate objective functional \eqref{eq:T2J} with the PDE constraints \eqref{eq:T2constraints}, we form the Lagrangian 
\beq
\begin{split}
	L_{T_2}(&u, v, (\hat{u}_n)_{n=1}^N, (\hat{v}_n)_{n=1}^N, (\lambda_n)_{n=1}^N, (\psi_n)_{n=1}^N, \\
	&u^*, v^*, (\hat{u}^*_n)_{n=1}^N, (\hat{v}^*_n)_{n=1}^N, (\lambda^*_{m,n})_{m,n=1}^N, (\psi^*_n)_{n=1}^N, \tau)\\
	& = J_{T_2}(\tau) \\
	& + A(u, v^*; \bar{\zeta}, \tau) - F(v^*; \bar{\zeta}, \tau)\\
	& + A(u^*, v; \bar{\zeta}, \tau) + \langle \partial_u Q(u), u^*\rangle \\
	& + \sum_{n=1}^N A(\hat{u}_n, \hat{v}^*_n; \bar{\zeta}, \tau)  + \langle \partial_{\zeta} A(u,\hat{v}^*_n; \bar{\zeta}, \tau)   + \partial_{\zeta} F(\hat{v}^*_n; \bar{\zeta}, \tau), \psi_n \rangle\\
	& + \sum_{n=1}^N A(\hat{u}^*_n, \hat{v}_n; \bar{\zeta}, \tau)+ \langle \partial_{uu} Q(u) \hat{u}_n , \hat{u}^*_n \rangle + \langle \partial_{\zeta} A(\hat{u}^*_n, v; \bar{\zeta}, \tau), \psi_n \rangle \\
	& + \sum_{n=1}^N \langle \bar{\cH} \psi_n, \psi_n^* \rangle - \langle \lambda_n \cC^{-1} \psi_n, \psi_n^* \rangle\\
	& + \sum_{m,n=1}^N \lambda_{m,n}^*  \left(\langle \cC^{-1} \psi_n, \psi_m \rangle - \delta_{mn}\right).
\end{split}
\eeq
By setting the variation of this Lagrangian with respect to $\lambda_n$ to zero, we obtain 
\beq\label{eq:psinstar}
\psi_n^* = \frac{1+2\beta_V \lambda_n}{2} \psi_n, \quad n = 1, \dots, N.
\eeq
By setting the variation with respect to $\hat{v}_n$ to zero, we have: find $\hat{u}_n^* \in V$ such that 
\beq
A(\hat{u}_n^*, \tilde{v}; \bar{\zeta}, \tau) = - \langle \partial_\zeta A(u, \tilde{v}; \bar{\zeta}, \tau) - \partial_\zeta F(\tilde{v};\bar{\zeta}, \tau ), \psi_n^* \rangle \quad \forall \tilde{v} \in V,
\eeq
which has the same form as the incremental forward Helmholtz equation \eqref{eq:incrementalState}, so that by \eqref{eq:psinstar} we have
\beq\label{eq:unhatstar}
\hat{u}_n^* = \frac{1+2\beta_V \lambda_n}{2} \hat{u}_n, \quad n = 1, \dots, N.
\eeq
Similarly, by setting the variation of $L_{T_2}$ with respect to $\hat{u}_n$ to zero, we have: find $\hat{v}_n^* \in V$ such that 
\beq
A(\tilde{u}, \hat{v}_n^*; \bar{\zeta}, \tau) = -  \langle \partial_{uu} Q(u) \tilde{u}, \hat{u}_n^* \rangle - \langle \partial_{\zeta} A(\tilde{u},v; \bar{\zeta}, \tau), \psi_n^* \rangle \quad \forall \tilde{u} \in V,
\eeq
which has the same form as the incremental adjoint Helmholtz equation \eqref{eq:incrementalAdjoint}, so that by \eqref{eq:psinstar} and \eqref{eq:unhatstar} we have
\beq
\hat{v}_n^* = \frac{1+2\beta_V \lambda_n}{2} \hat{v}_n, \quad n = 1, \dots, N.
\eeq
Then, by setting the variation of $L_{T_2}$ with respect to $v$ to zero, we obtain: find $u^* \in V$ such that 
\beq
\begin{split}
	A(u^*, \tilde{v}; \bar{\zeta}, \tau)  = &- 2\beta_V  \langle \partial_\zeta A(u, \tilde{v}; \bar{\zeta}, \tau) - \partial_\zeta F(\tilde{v}; \bar{\zeta}, \tau), \cC \bar{g} \rangle \\
	& - \langle \partial_{\zeta} A(\hat{u}^*_n, \tilde{v}; \bar{\zeta}, \tau), \psi_n \rangle - \langle \partial_\zeta A(\hat{u}_n, \tilde{v}; \bar{\zeta}, \tau), \psi_n^* \rangle  \\
	& - \langle \partial_{\zeta\zeta} A(u,\tilde{v}; \bar{\zeta}, \tau) \psi_n - \partial_{\zeta\zeta} F(\tilde{v}; \bar{\zeta}, \tau)\psi_n, \psi_n^* \rangle \quad \forall \tilde{v} \in V.
\end{split}
\eeq
Finally, by setting the variation of $L_{T_2}$ with respect to $u$ to zero, we obtain: find $v^* \in V$ such that 
\beq
\begin{split}
	A(\tilde{u}, v^*; \bar{\zeta}, \tau) = &- \langle \partial_u Q(u), \tilde{u} \rangle - 2\beta_V \langle \partial_\zeta A(\tilde{u}, v;  \bar{\zeta}, \tau), \cC \bar{g} \rangle \\
	& - \langle \partial_{uu} Q(u) u^*, \tilde{u} \rangle - \langle \partial_\zeta A(\tilde{u}, \hat{v}_n^*; \bar{\zeta}, \tau), \psi_n \rangle \\
	& - \langle \partial_\zeta A(\tilde{u}, \hat{v}_n;  \bar{\zeta}, \tau)+  \partial_{\zeta\zeta} A(\tilde{u}, v;  \bar{\zeta}, \tau) \psi_n, \psi_n^* \rangle \quad \forall\tilde{u}\in V.
\end{split}
\eeq
Note that the design variable $\tau$ is not involved in the orthonormality condition of the eigenfunctions, so there is no need to compute $\lambda_{m,n}^*$ in the Lagrangian. 
With all the other Lagrange multipliers available, we can compute the $\tau$-gradient as 
\beq
\begin{split}
	\nabla_\tau J_{T_2}(\tau) = \partial_\tau L(&u, v, (\hat{u}_n)_{n=1}^N, (\hat{v}_n)_{n=1}^N, (\lambda_n)_{n=1}^N, (\psi_n)_{n=1}^N, \\
	&u^*, v^*, (\hat{u}^*_n)_{n=1}^N, (\hat{v}^*_n)_{n=1}^N, (\lambda^*_{m,n})_{m,n=1}^N, (\psi^*_n)_{n=1}^N, \tau).
\end{split}
\eeq

\subsection{The approximate Newton algorithm}
\label{sec:approximateHessianalg}
Once the $\tau$-gradient is computed for the different approximations, and the $\tau$-Hessian action is computed for the deterministic approximation, we can solve the optimization problem by an approximate Newton algorithm with backtracking line search to guarantee monotonic convergence, where the $\tau$-Hessian is computed or approximated by the $\tau$-Hessian of the deterministic approximation, and the resulting linear system is solved by inexact preconditioned conjugate gradient method with Steihaug's stopping criteria.
\begin{algorithm}[!htb]
	\caption{Line search inexact approximate Newton--pCG algorithm}
	\label{alg:InexactNewtonCG}
	\begin{algorithmic}
		\STATE{\textbf{Input: }  the maximum numbers of approximate Newton, CG, and line search iterations $N_{qn}$, $N_{cg}$, $N_{ls}$, and the convergence tolerance $\varepsilon_{qn}$ for the approximate Newton.}
		\STATE{\textbf{Output: } solution of the optimization problem $\tau^*$.}
		\STATE{1. Initialize a design variable $\tau_0$, set $n_{qn}, n_{cg}, n_{ls} = 0$, set the tolerance $\epsilon_{qn} = 2 \varepsilon_{qn}$, set the tolerance for CG convergence to $\varepsilon_{cg} = \varepsilon_{cg}^0$.}
		\WHILE{ $n_{qn} < N_{qn}$ and $\epsilon_{qn} < \varepsilon_{qn}$}
		\STATE{2. Solve the Newton linear system: find the update direction $\delta\tau \in Z$ by solving
			\beq\label{eq:quasiNewton}
			\nabla_\tau^2 J_{\bar{\zeta}}(\tau_{n_{qn}}) \; \delta \tau = - \nabla_\tau J_{a}(\tau_{n_{qn}})
			\eeq
			using a CG method with preconditioner $\beta_P \nabla_\tau^2 P(\tau)$,  terminated by Steihaug's criteria, namely, when $n_{cg} \geq N_{cg}$, or $\epsilon_{cg} \geq \varepsilon_{cg}$, or $\langle \nabla_\tau^2 J_{\bar{\zeta}}(\tau_{n_{qn}}) \; \delta \tau, \delta \tau \rangle < 0$ (i.e., when a direction of negative curvature is encountered).}
		\vspace*{0.2cm}
		\WHILE{$J_a(\tau_{n_{qn}} + \alpha \delta \tau) > J_a(\tau_{n_{qn}}) +c_{\text{AG}}  \alpha  \delta \tau$ and $n_{ls} < N_{ls}$}
		\STATE{3. Set $\alpha = 2^{-n_{ls}}$ and compute $J_a(\tau_{n_{qn}} + \alpha \delta \tau) $. Set $n_{ls} \leftarrow n_{ls} + 1$. }
		\ENDWHILE
		\vspace*{0.2cm}
		\STATE{4. Break the while loop if $n_{ls} \geq N_{ls}$.}
		\vspace*{0.2cm}
		\STATE{5. Set $\tau_{n_{qn}+1} = \tau_{n_{qn}} + \alpha \delta \tau$, $n_{qn} \leftarrow n_{qn} + 1$, $n_{cg}, n_{ls} = 0$, compute $\epsilon_{qn}$, and update the tolerance for CG convergence at $\varepsilon_{cg} = \min\{\varepsilon_{cg}^0, ||\nabla J_a(\tau_{n_{cg}})||/||\nabla J_a(\tau_0)||\}$.}
		\ENDWHILE
		\STATE{6. Set $\tau^* = \tau_{n_{qn}}$.}
	\end{algorithmic}
\end{algorithm}

The method is summarized in Algorithm Algorithm \ref{alg:InexactNewtonCG}. In \eqref{eq:quasiNewton} of Algorithm \ref{alg:InexactNewtonCG}, $J_a$ represents the approximation of the objective functional $J \approx J_a$, with the deterministic approximation $J_a = J_{\bar{\zeta}}$, the sample average approximation $J_a = J_{\text{SAA}}$, and the Taylor approximation $J_a = J_{T_2}$. For the sample average approximation and Taylor approximation, the $\tau$-Hessian of $J_a$ is approximated by the $\tau$-Hessian of the deterministic approximation $\nabla_\tau^2 J_{\bar{\zeta}}$ given by \eqref{eq:hessiantau}, while the gradients are computed as in Sections \ref{sec:gradDeterministic}, \ref{sec:gradSAA}, and \ref{sec:gradTaylor} for the deterministic, SAA, and Taylor approximations, respectively. For the termination condition in step $n_{qn}+1$ of the approximate Newton iteration, we can use a quantity related to the norm of the gradient $||\nabla_\tau J_a(\tau_{n_{qn}})||$ and/or $\langle \nabla_\tau J_a(\tau_{n_{qn}}), \delta \tau \rangle$. $c_{\text{AG}}$ is a small constant for Armijo--Goldstein conditions, e.g., $c_{\text{AG}} =  10^{-4}$.

In each of the approximate Newton iteration, we have to compute once the $\tau$-gradient $\nabla_\tau J_a$, perform $n_{cg}$ $\tau$-Hessian actions, i.e., the actions of $\nabla_\tau^2 J_{\bar{\zeta}}$ in given CG directions while solving \eqref{eq:quasiNewton}, which requires solution of a pair of incremental forward/adjoint Helmholtz equations \eqref{eq:incrementalStatetau} and \eqref{eq:incrementalAdjointtau} for each Hessian action, as well as $n_{ls}$ backtracking line search iterations, which requires $n_{ls}$ evaluations of $J_a$. For relatively small uncertainty, i.e., small signal-to-noise level, we expect that the Hessian $\nabla_\tau^2 J_{\bar{\zeta}}$ is a good approximation of $\nabla_\tau^2 J_{a}$, and the total number of Newton iterations $n_{qn}$ is independent of the dimension of the discretized design variable field. Moreover, the number of preconditioned CG iterations $n_{cg}$ is also expected to be independent of the design variable dimension when the $\tau$-Hessian of the approximation for $\bE[Q] + \beta_V \text{Var}[Q]$ is low-rank. Therefore, the inexact approximate Newton--pCG algorithm is expected to be scalable with respect to the dimension of the design variable $\tau$, in the sense that the number of Helmholtz solves will be independent of the design variable dimension. This will be demonstrated numerically in Section \ref{sec:scalability}.

\section{Numerical experiments}
\label{sec:numerics}

In this section, we present several numerical experiments to: (1) demonstrate the effectiveness of the optimization strategy in a deterministic setting, (2) compare the difference between various approximation methods for the optimization under uncertainty, (3) illustrate the scalability of the Taylor approximation and the approximate Newton-pCG algorithm with respect to the dimension of the discretized random variable and design variable fields, respectively, (4) show that the proposed optimization strategy can achieve cloaking for incident waves with multiple directions and multiple frequencies under uncertainty, and finally (5) elucidate the applicability of the proposed optimization strategy to more complex geometries beyond disks.   

\subsection{Cloaking in a deterministic setting}
In this experiment, we aim to demonstrate the effectiveness of the optimization strategy in designing a cloak that makes the obstacle invisible to acoustic waves. In what follows, we use normalized units for all quantities. The configuration of the design problem is displayed in Fig.\ \ref{fig:SketchAcousticWave}, where the obstacle is a disk of radius $r_1 = 1$, which is surrounded by the cloaking region with radius $r_2 = 3$, and immersed in a host square medium of size $[-6,6]^2$ with PML boundaries of length $1$ on all sides. The incident wave is a plane wave propagating from the left side to the right side, given by $u^{\text{inc}} = e^{i k_0 x \cdot b}$ with direction $b = (1,0)$ and wavenumber $k_0 = \omega/c_0$ with frequency $\omega = 2\pi$ and sound speed $c_0 = 1$. For this experiment we do not consider uncertainty in the optimal design and fix the random variable at its mean $\zeta = \bar{\zeta} = 0$ in \eqref{eq:soundSpeed}. This approach is equivalent to the deterministic approximation presented in Section \ref{sec:deterministicApproximation}. For the regularization of the design variable, we set $\beta_P = 10^{-2}$ in the objective functional \eqref{eq:optimizationDeter}. A finite element method is used to solve the scattering problem, with mesh of triangles with 172,803 vertices, leading to 345,606, 34,217, and 57,462 degrees of freedom for the discrete state variable (using piecewise linear elements), the discrete random variable (piecewise linear elements), and the discrete design variable (piecewise constant elements), respectively.

\begin{figure}[!htb]
	\begin{center}
		\includegraphics[scale=0.45]{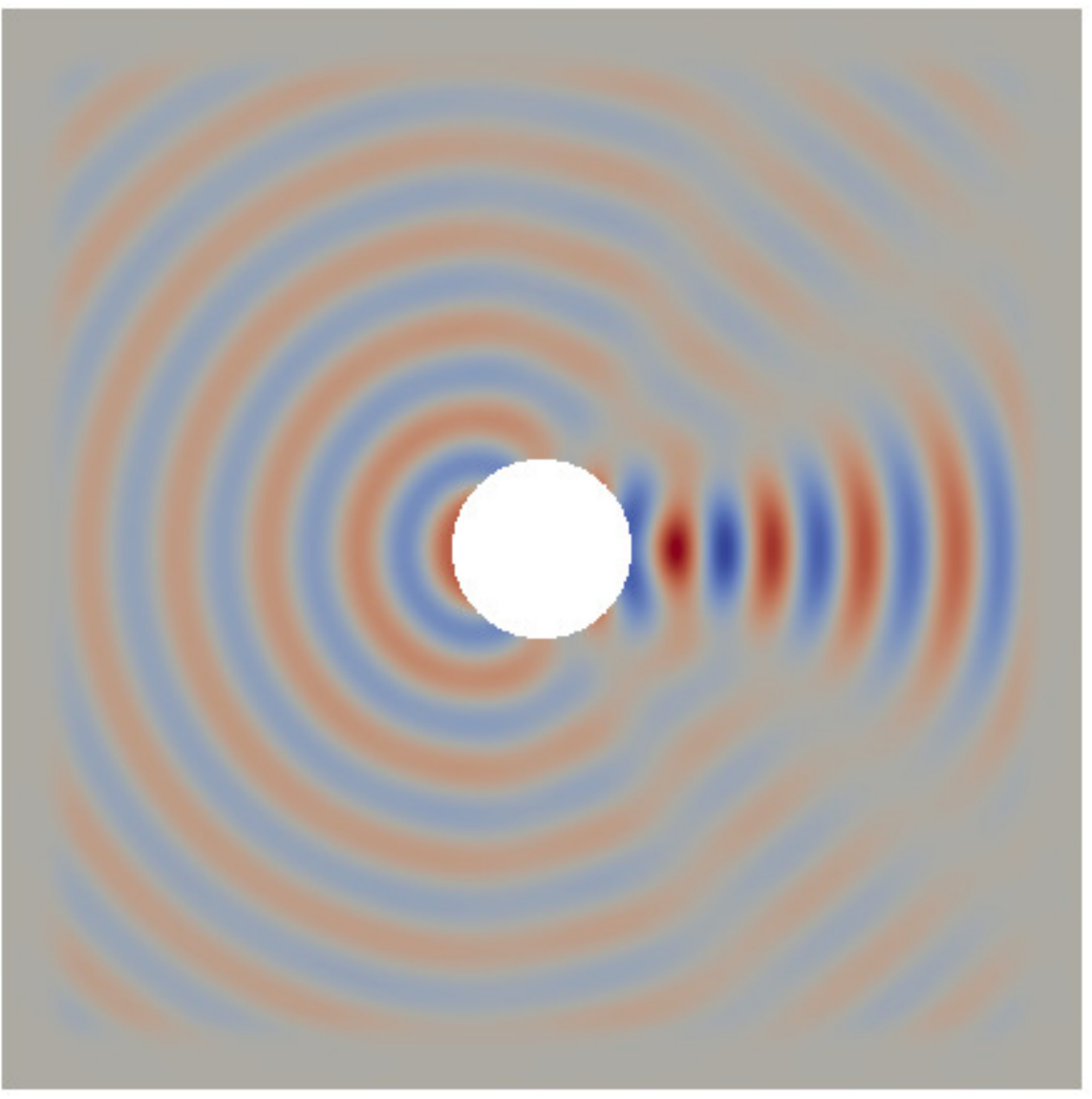}
		\includegraphics[scale=0.45]{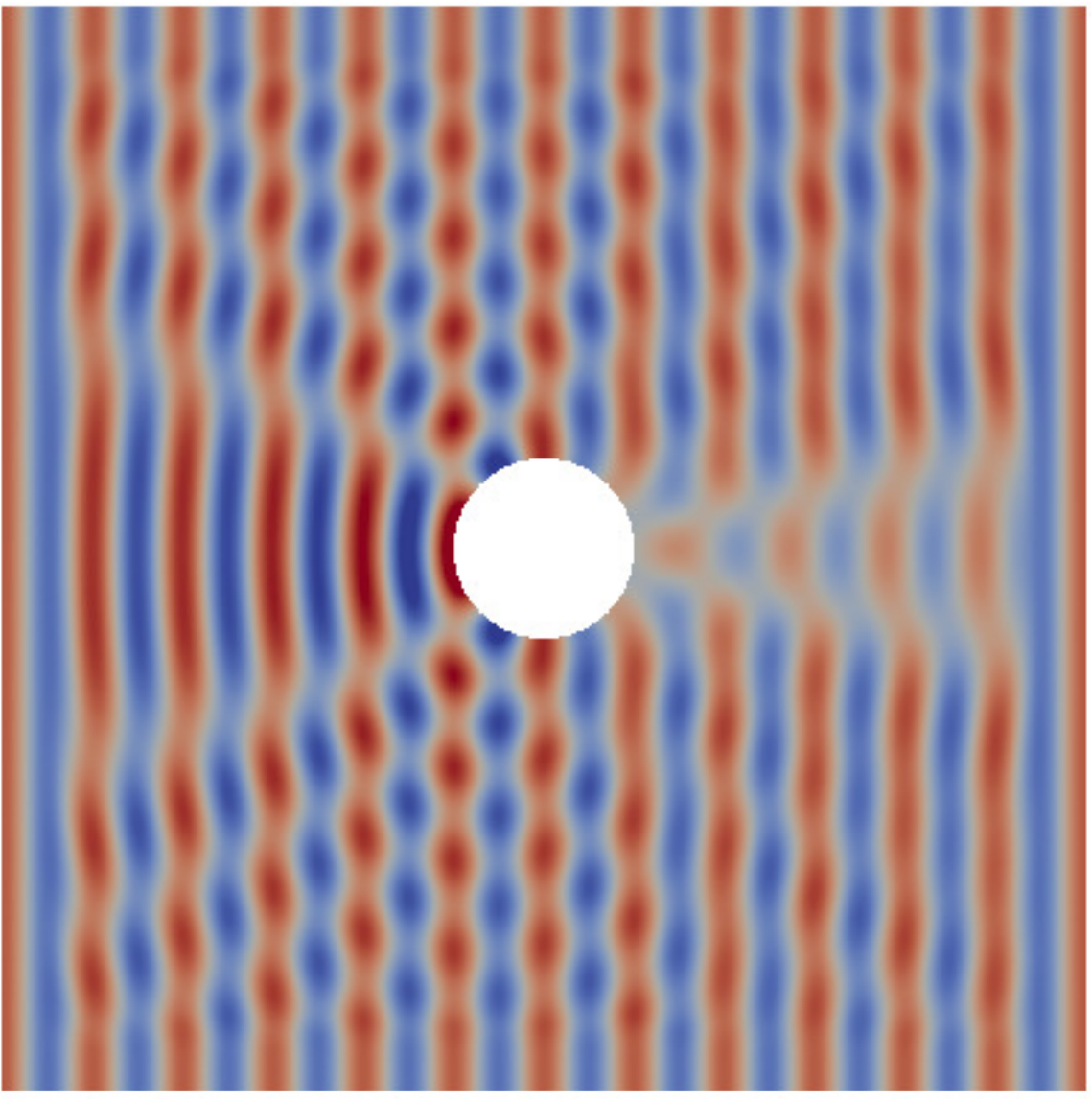}
		\includegraphics[scale=0.45]{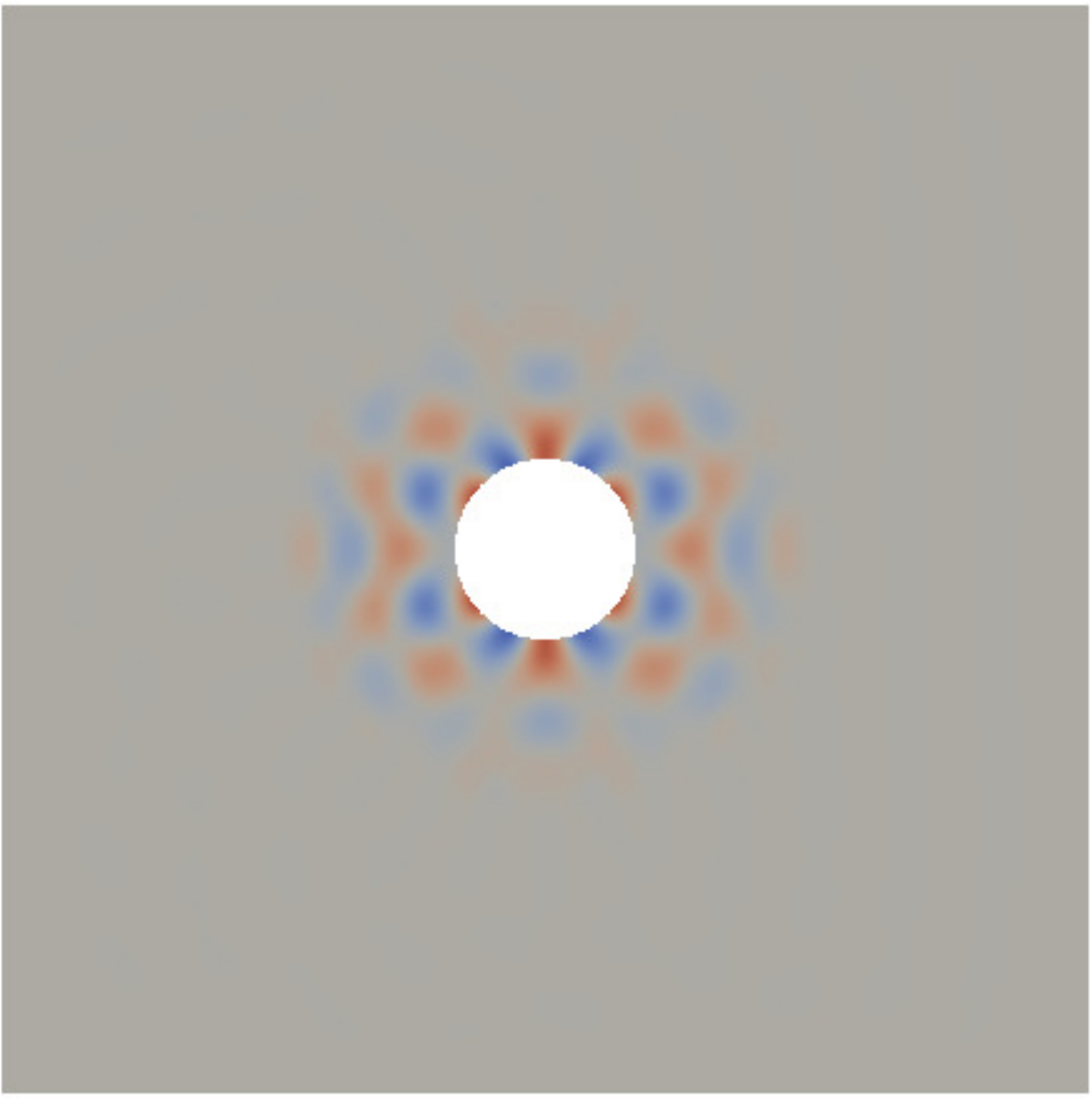}
		\includegraphics[scale=0.45]{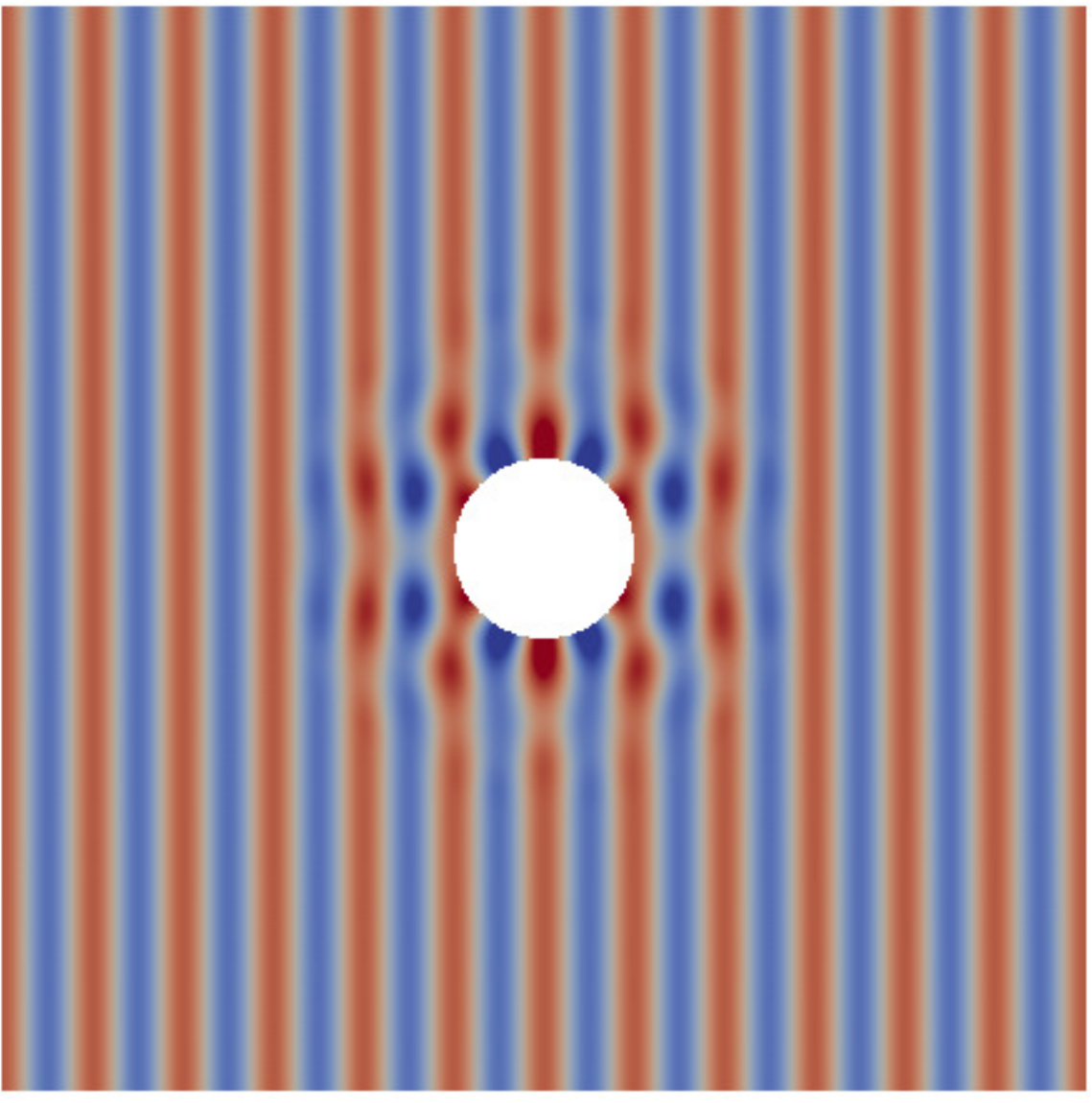}
		\includegraphics[scale=0.45]{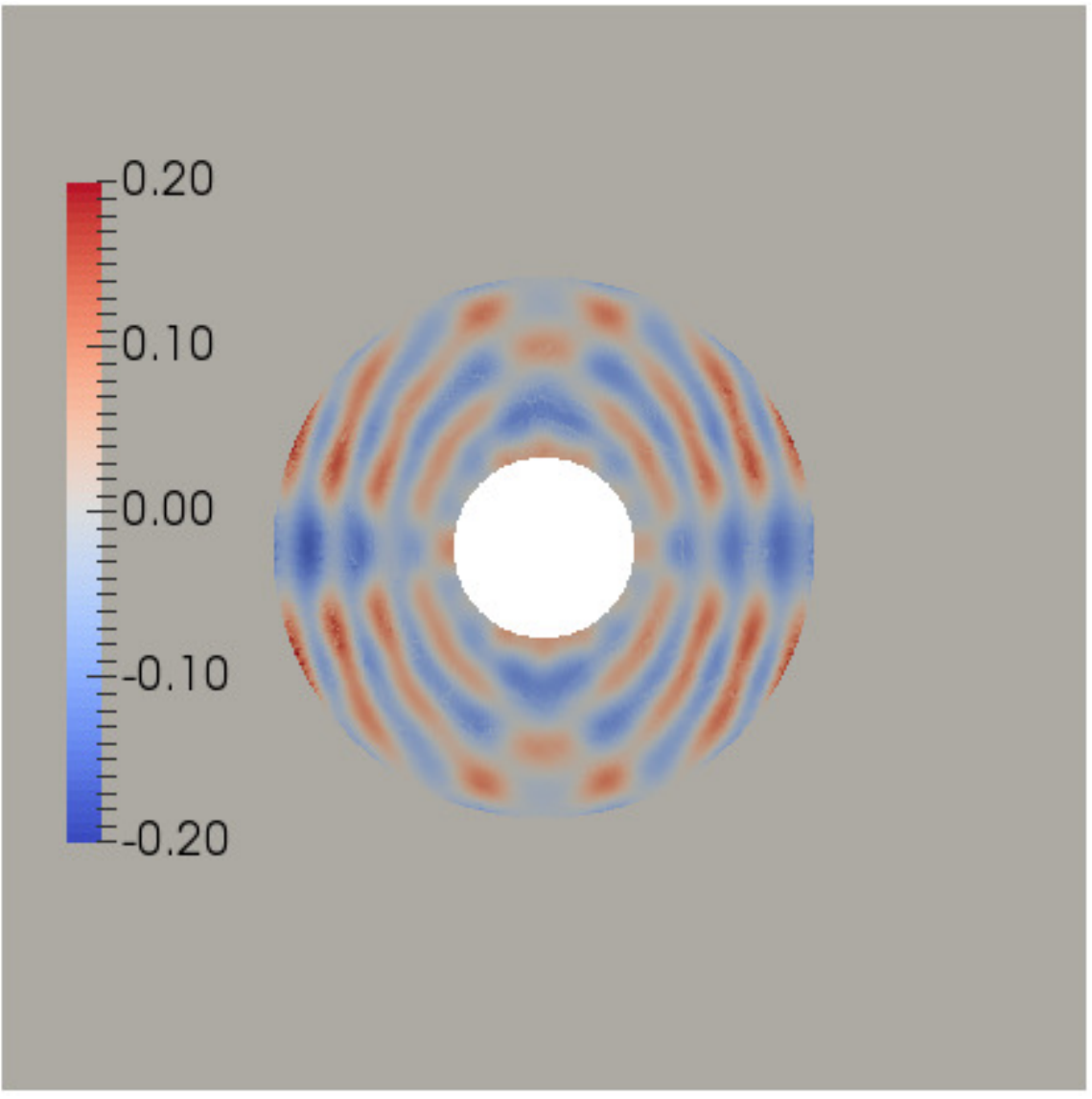}
		\includegraphics[scale=0.45]{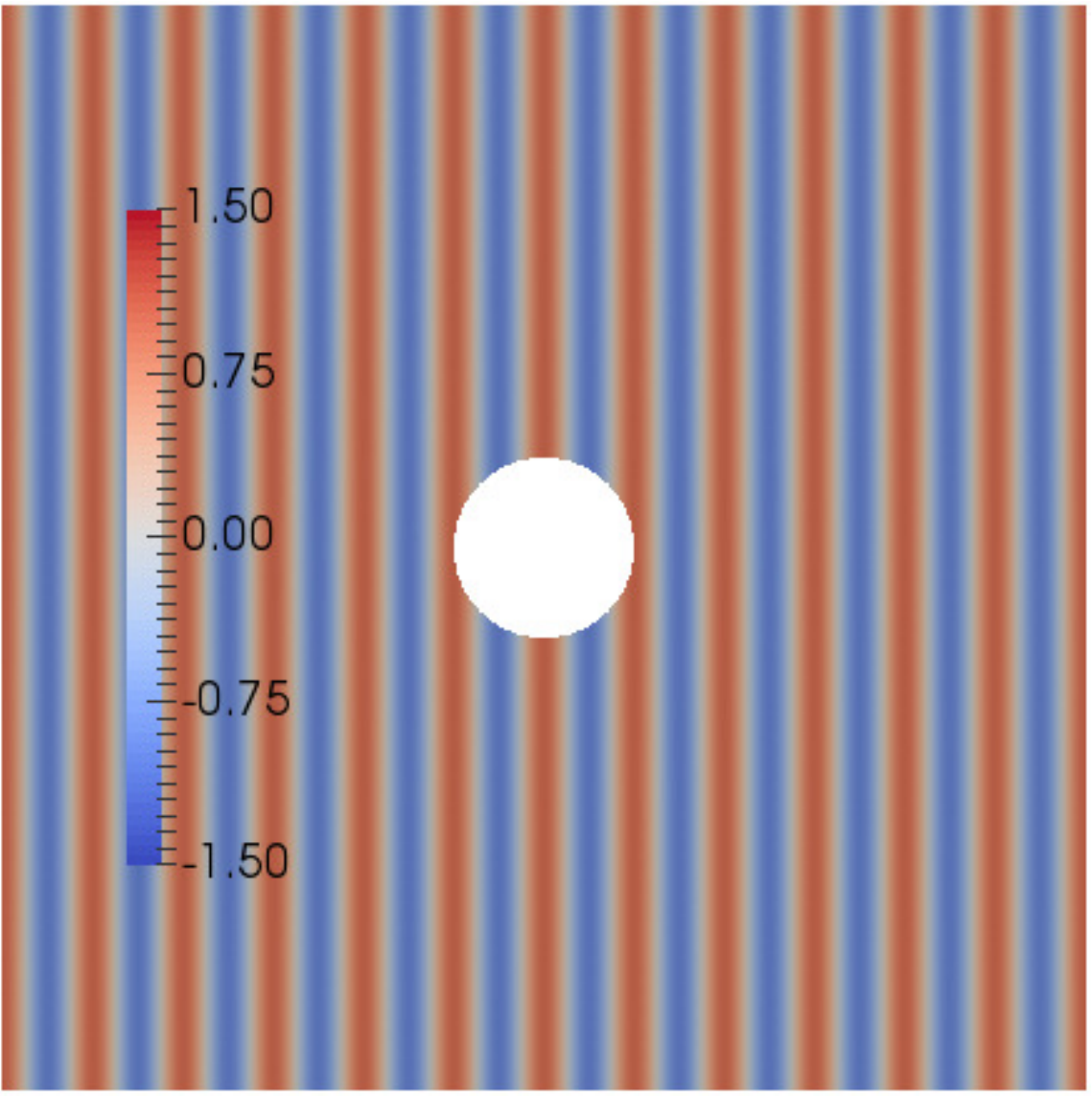}
	\end{center}
	\caption{Top: wave scattering from an impenetrable obstacle; left: scattered wave field; right: total wave field. Middle: wave scattering with the optimized cloak; left: scattered wave field; right: total wave field. Bottom, left: the optimal design variable field $\tau^*$ obtained by the deterministic optimization; right: the incident wave field, i.e., total wave field in homogeneous medium. The real part of all wave fields are shown.}\label{fig:cloaking}
\end{figure}

We initialize the design variable $\tau = 0$ in \eqref{eq:soundSpeed} and run the approximate Newton algorithm as presented in Algorithm \ref{alg:InexactNewtonCG} to minimize the objective functional \eqref{eq:optimizationDeter} with respect to the design variable $\tau$, with $N_{qn}=10$, $N_{cg}=10$, $N_{ls} = 10$, and $\varepsilon_{qn} = 10^{-2}$. The algorithm converged in 6 iterations. The results are shown in Fig.\ \ref{fig:cloaking} with the real part of the scattered and total wave fields shown in the top two images, in which the reflection of the incident wave from the impenetrable obstacle without  the cloak is evident. In the middle two images, the scattered and total wave fields are displayed with  the cloak at the optimal design. From the middle-left image, we can see a clear reduction of the scattered wave in the observation region---which is essentially invisible outside the cloak region. Inside the cloak region, the scattered wave fields is significantly altered from that without the cloak. From the middle-right images of the total wave field, we can observe an effective cloaking of the obstacle, i.e., the total field coincides with the incident field outside the design region as shown in the bottom-right image. All wave fields are scaled to the range $[-1.5, 1.5]$ for the sake of visual comparison. The optimal design variable $\tau^*$ is shown in the bottom-left image, from which we can see a sub-wavelength structure within the cloak, which effectively eliminates the scattered wave in the background medium rendering it undetectable to an external observer. It is worth noting some similarity in the cloaking structure between this approach, which permits continuously varying material properties are possible, and the cloaks constructed from distributions of discrete scatterers reported in references \cite{SanchisGarcia-ChocanoLlopis-PontiverosEtAl13,AndkjerSigmund13,LuSanchisWenEtAl18}. 

\subsection{Cloaking under uncertainty}
\begin{figure}[!htb]
	\begin{center}
		\includegraphics[scale=0.5]{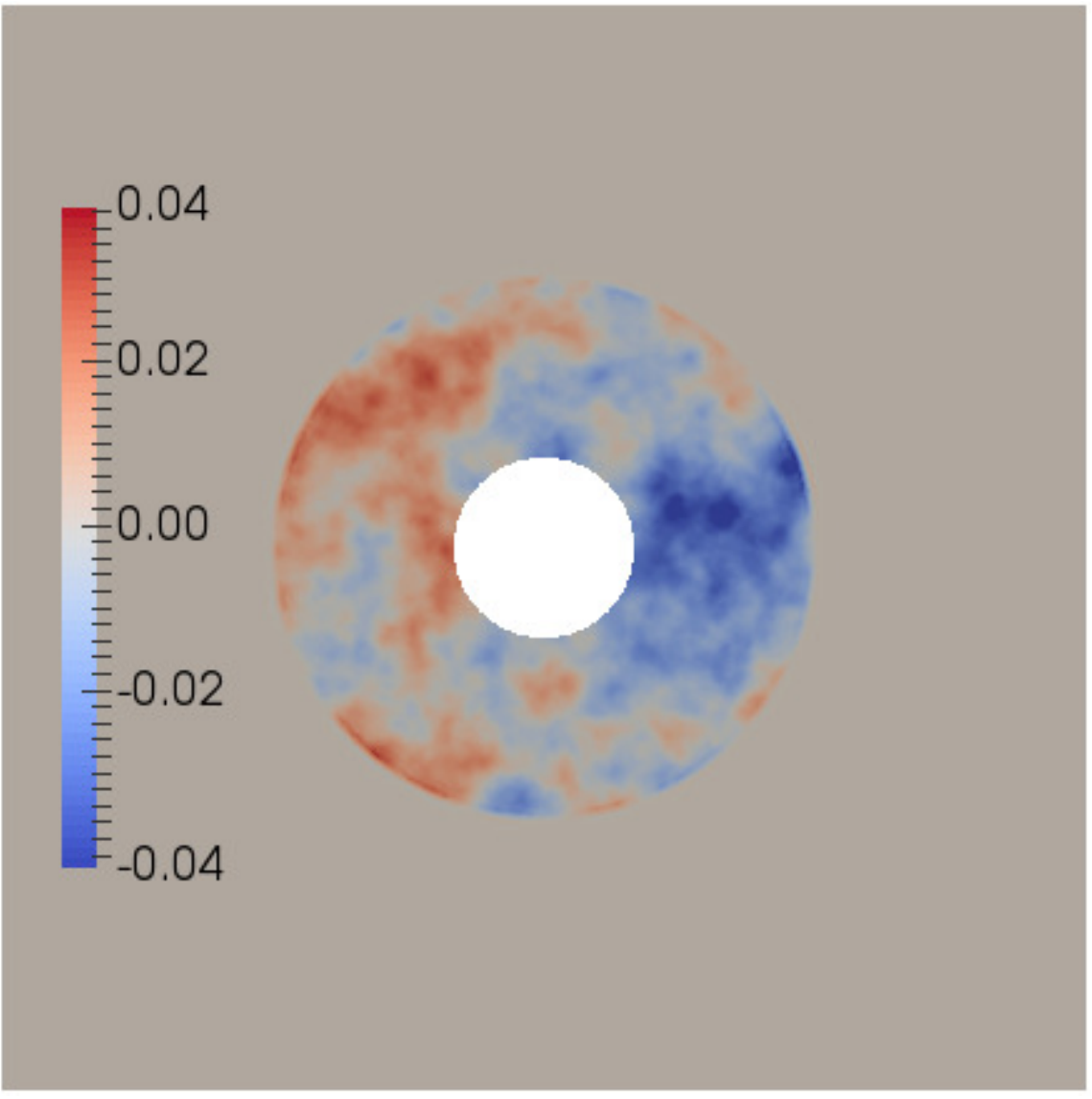}
		\includegraphics[scale=0.5]{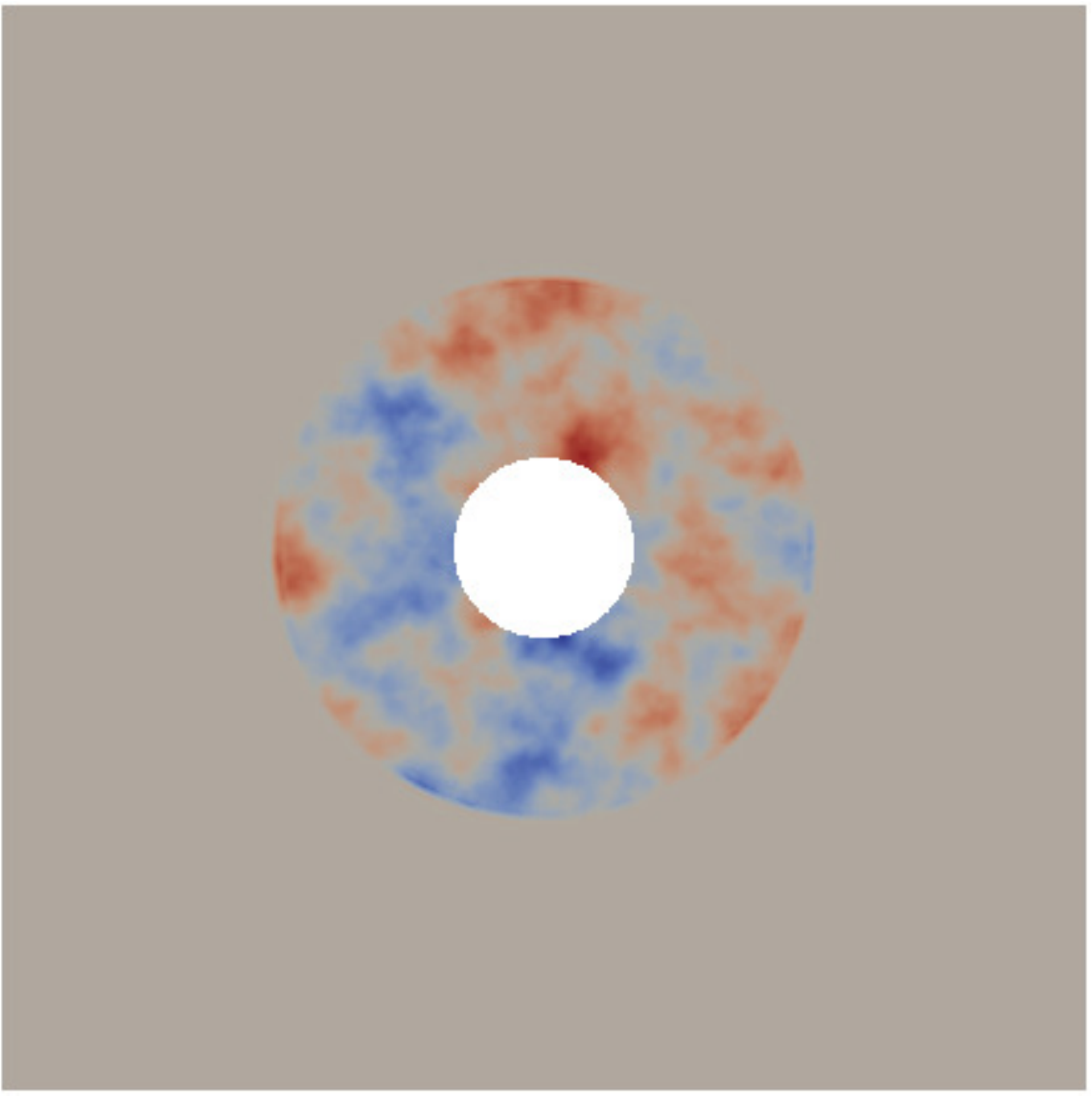}
	\end{center}
	\caption{Random samples of $\zeta \sim \cN(\bar{\zeta}, \cC)$ with $\bar{\zeta} = 0$, and $\gamma = 10, \delta = 50$ for $\cC = (-\gamma \Delta + \delta I)^{-2}$. }\label{fig:samples}
\end{figure}

In this experiment, we compare the optimal cloaking performance under uncertainty by the three approximation methods presented in Section \ref{sec:approximation}. This uncertainty, due to manufacturing errors or variability in material properties, is modeled as an additive Gaussian random field $\cN(\bar{\zeta},\cC)$ with the covariance operator $\cC = (-\gamma \Delta + \delta I)^{-2}$. We take $\gamma = 10$ and $\delta = 50$ such that the noise-to-signal ratio of the random variable is about $20\%$ of the design variable. Two samples of the random (field) variable are shown in Fig.\ \ref{fig:samples}. 

\begin{figure}[!htb]
	\begin{center}
		\includegraphics[scale=0.28]{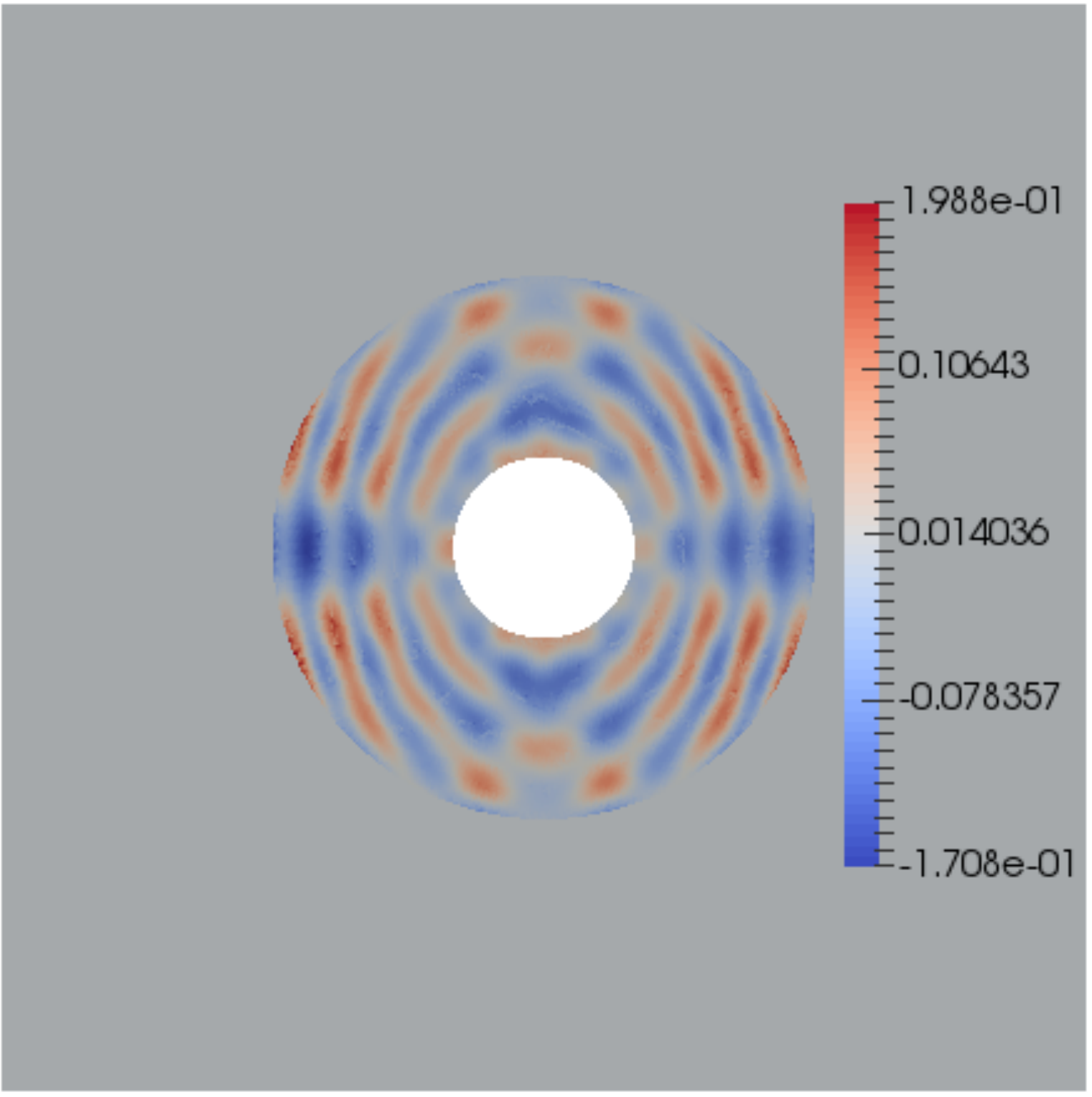}
		\includegraphics[scale=0.28]{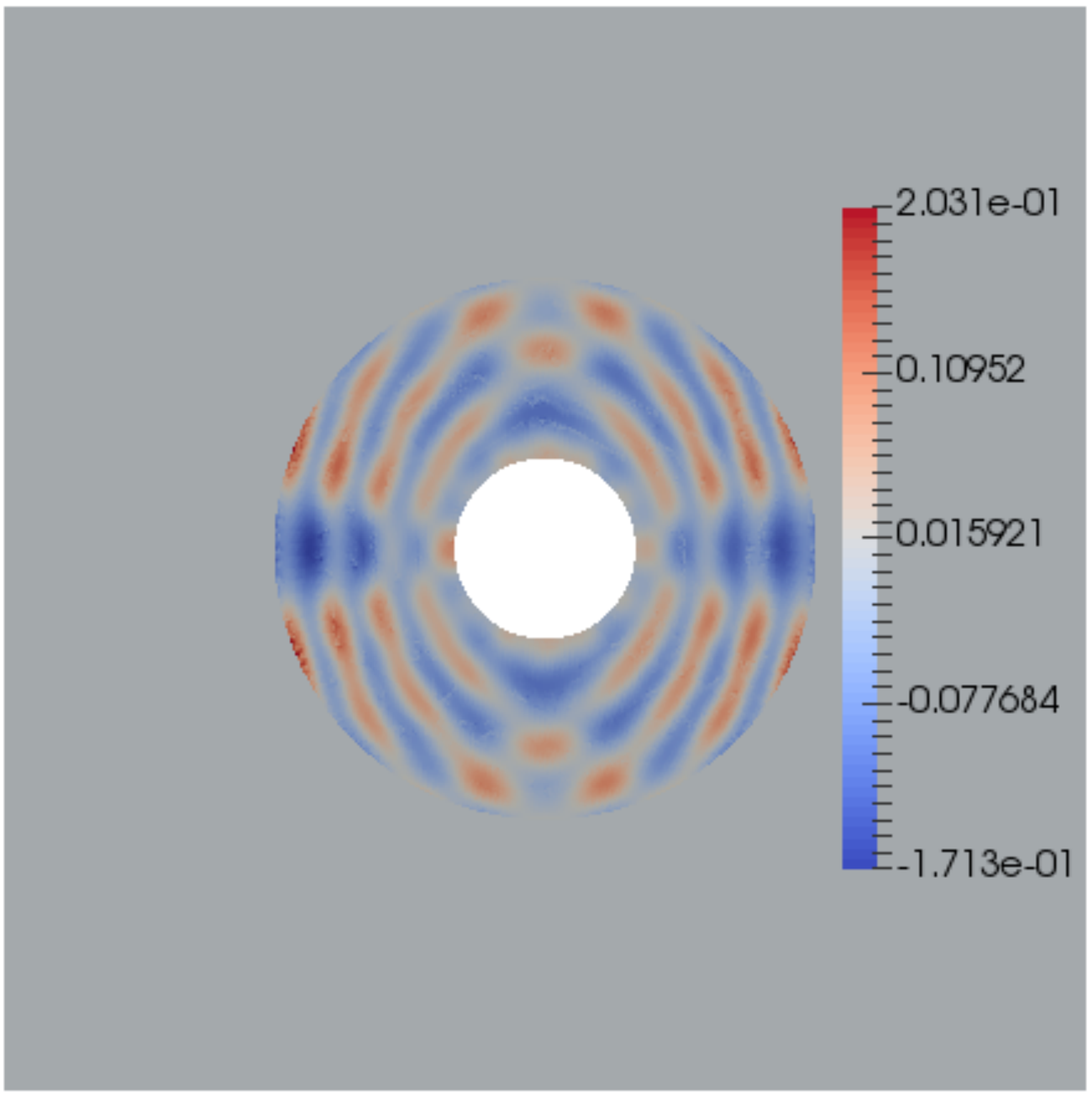}
		\includegraphics[scale=0.28]{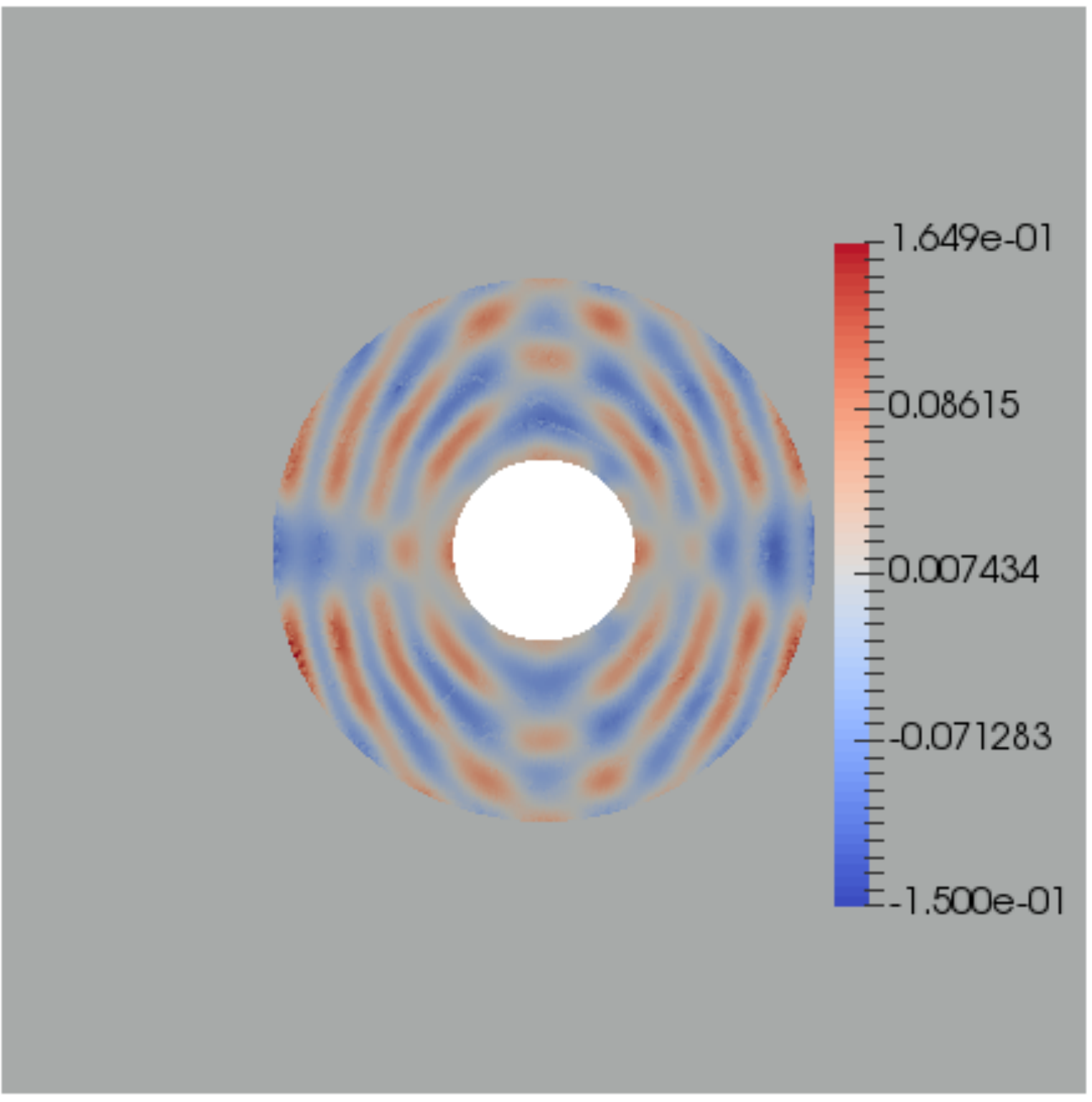}
	\end{center}
	\caption{Optimal design variable field $\tau^*$ obtained by using  deterministic (left), quadratic (middle), and sample average (right) approximations of the objective functional.}\label{fig:designs}
\end{figure}

The optimal design variables obtained by using different approximations of the objective functional are shown in Fig.\ \ref{fig:designs}. We use 50 eigenvalues in the trace estimate \eqref{eq:trace} for the quadratic approximation, which achieves about $99\%$ accuracy (shown in the next section). One hundred samples are used for the sample average approximation, which requires similar computational cost as the quadratic approximation. Slight differences can be noticed even though they share the same topological structure. 

\begin{figure}[!htb]
	\begin{center}
		\includegraphics[scale=0.35]{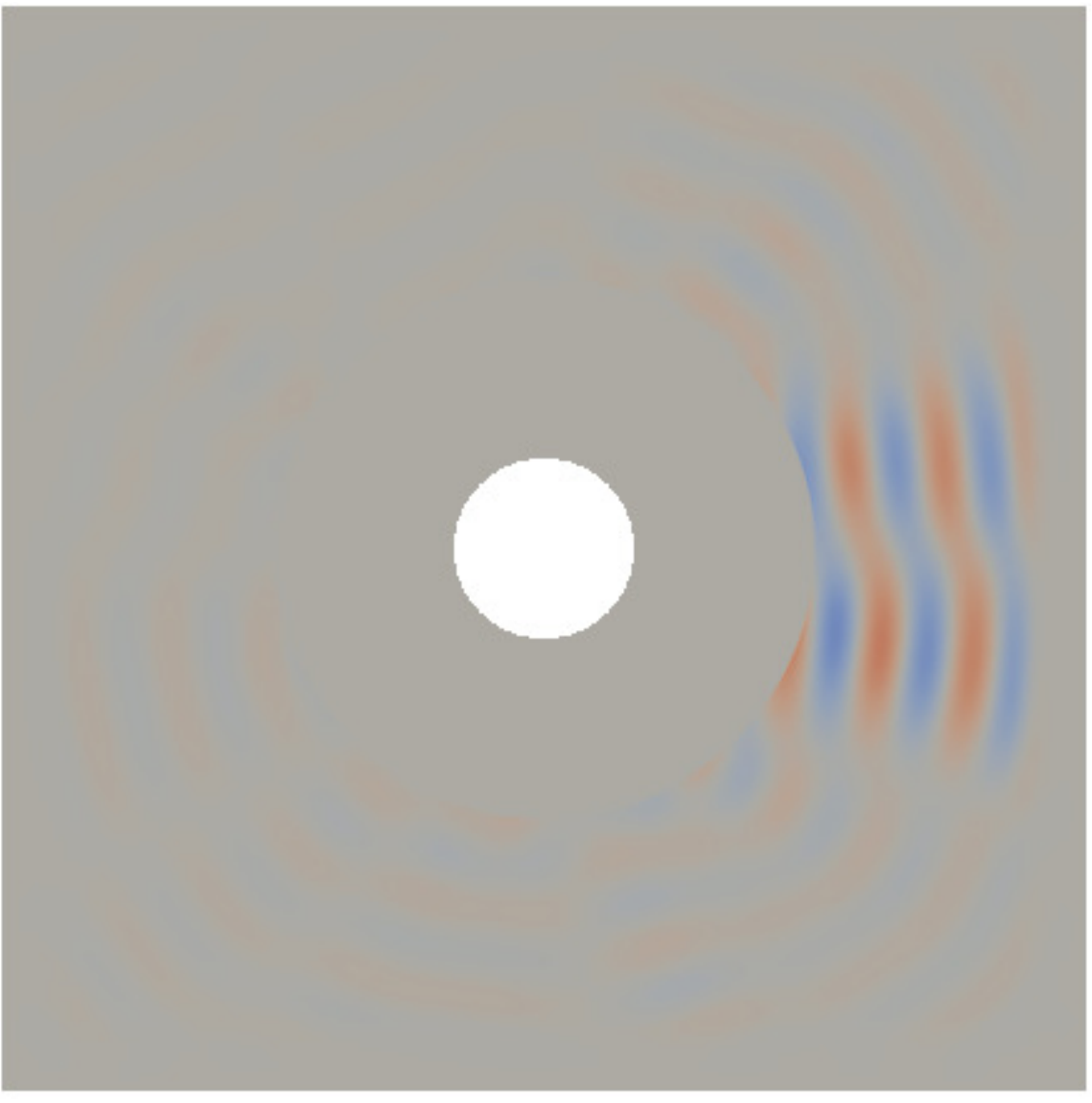}
		\includegraphics[scale=0.35]{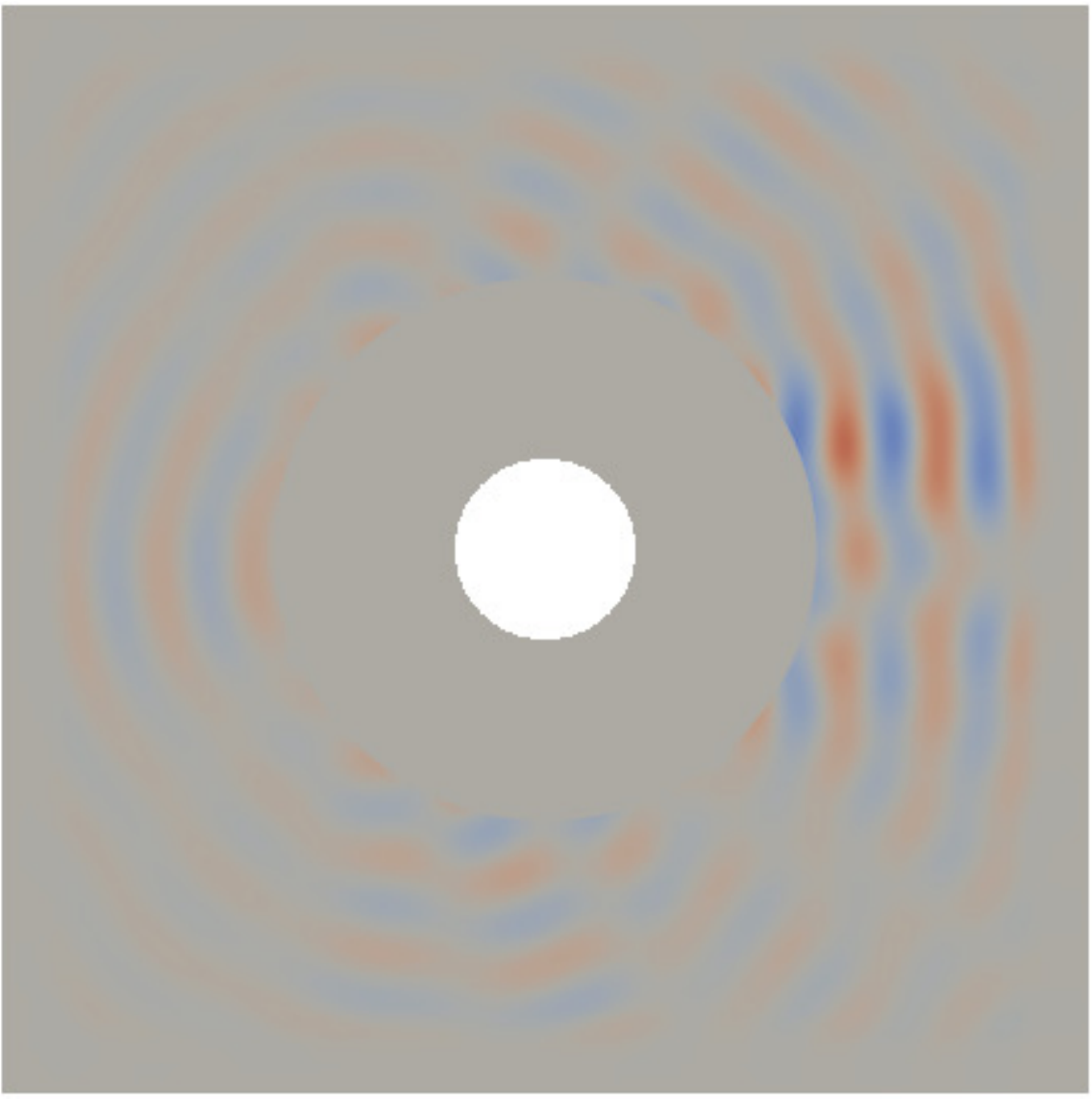}
		\includegraphics[scale=0.35]{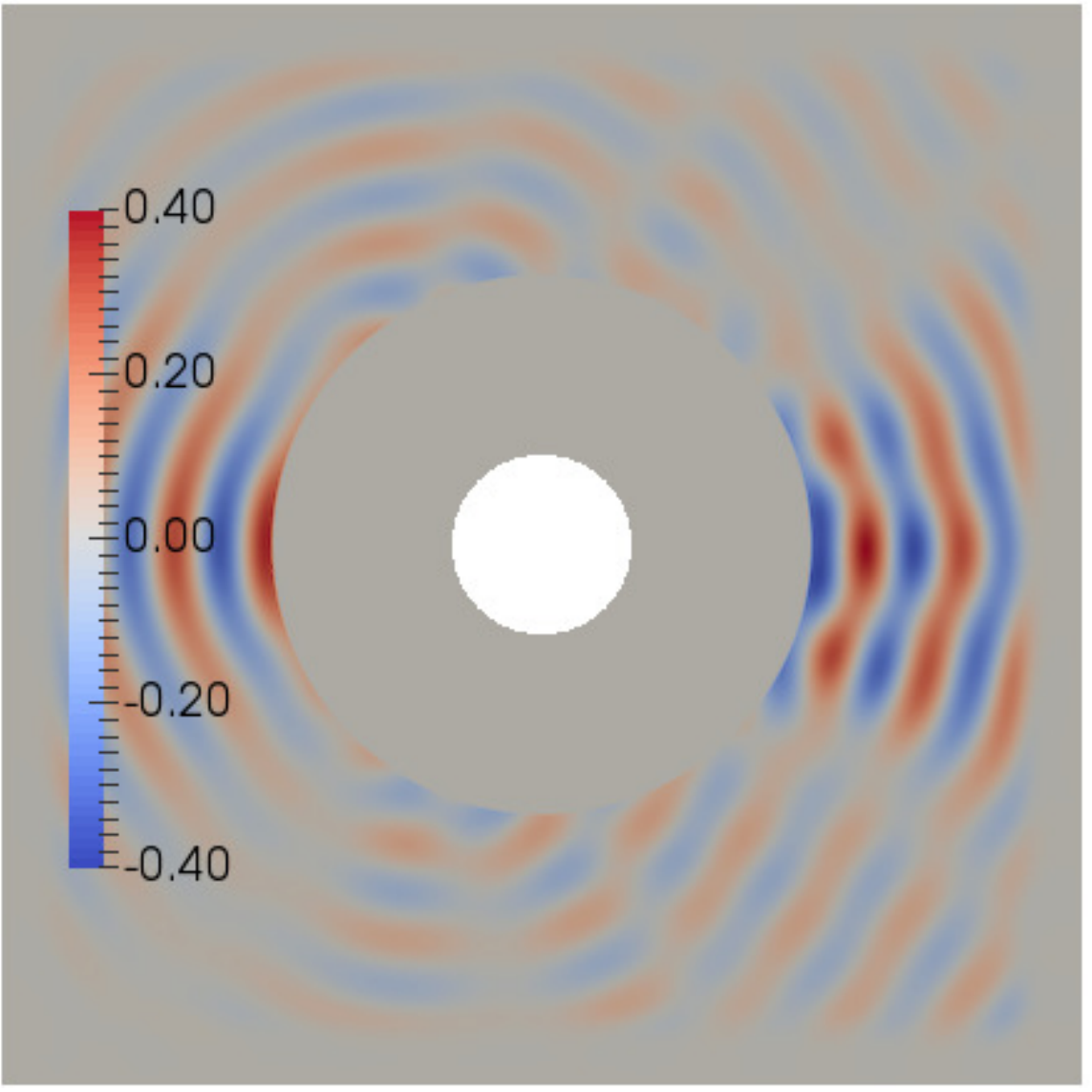}
		
		\includegraphics[scale=0.35]{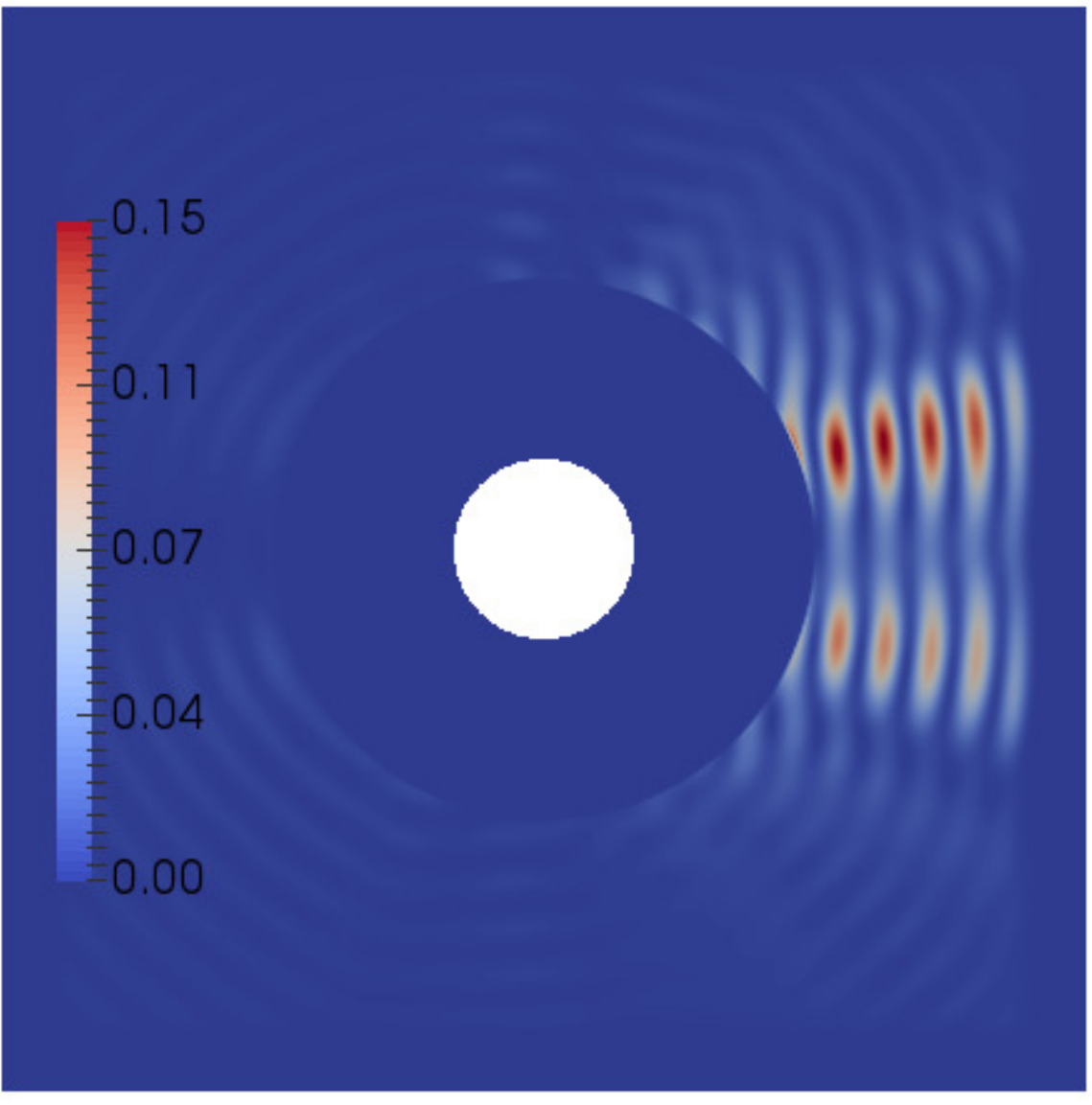}
		\includegraphics[scale=0.35]{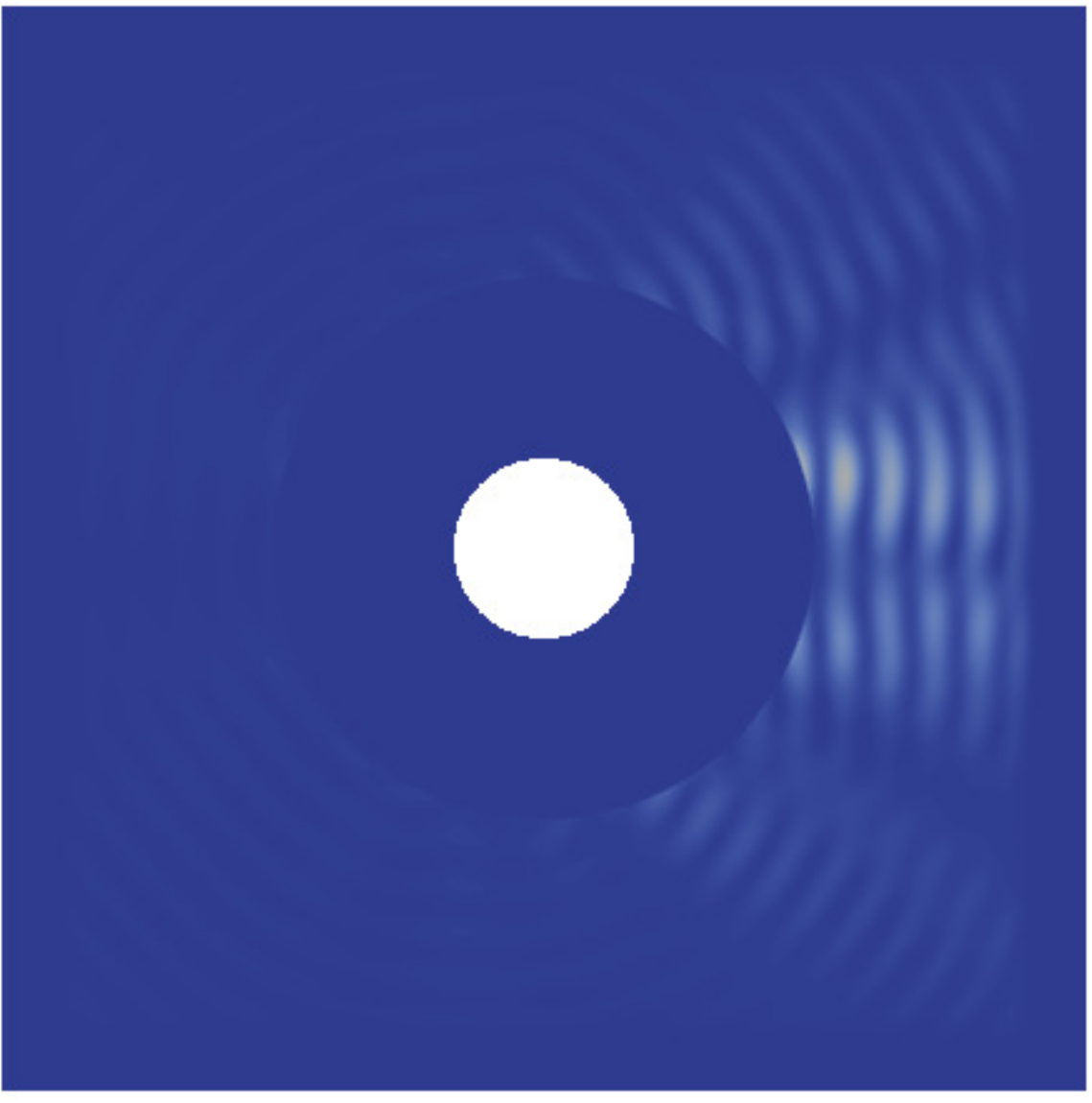}
		\includegraphics[scale=0.35]{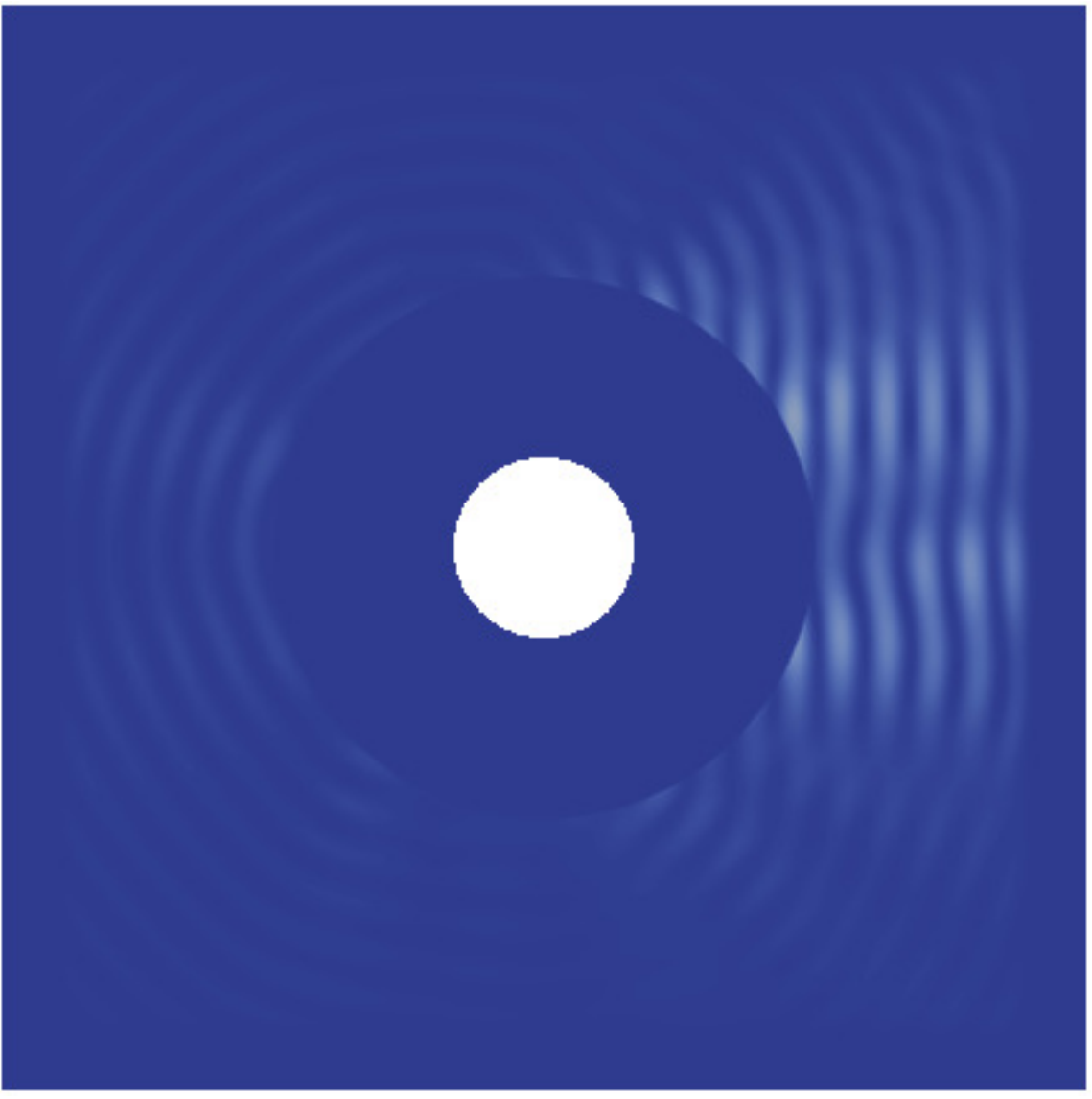}
	\end{center}
	\caption{Mean (top) and standard deviation (bottom) of the scattered wave field at the optimal design $\tau^*$ obtained by using deterministic (left), quadratic (middle), and sample average (right) approximations of the objective functional.}\label{fig:meanstd}
\end{figure}

We next draw 10 random samples of the random variable $\zeta$, and solve the scattering problem for each optimal design field. The mean and standard deviation of the scattered fields for the 10 random samples are shown in Fig.\ \ref{fig:meanstd}. We can observe that the sample average approximation leads to a more biased scattered field (as seen from its large mean), while the deterministic approximation gives rise to large variation of the scattered field (as seen from its large standard deviation).  

%

\begin{table}[!htb]
	\caption{Estimates $\hat{Q}$ of the design objective $Q$ and mean squared errors (MSE) for $\hat{Q}$, $Q - T_1 Q$, and $Q-T_2 Q$, based on 10 samples for a random design $\tau_{\text{random}}$, the deterministic optimal design $\tau_{\text{deter}}$, and the optimal design under uncertainty using the quadratic $\tau_{\text{quad}}$ and the sample average $\tau_{\text{saa}}$ approximations.}\label{tab:MSEQ}
	\begin{center}
		\begin{tabular}{|c|c|c|c|c|}
			\hline
			design   &     $\hat{Q}$       &     MSE($\hat{Q}$)  &     MSE($Q-T_1Q$) &     MSE$(Q-T_2Q$) \\
			\hline
			$\tau_{\text{random}} $ &  1.19E+01        & 8.34E$-$02     &     4.50E$-$03        &     4.89E$-$05 \\
			\hline 
			$\tau_{\text{deter}}$  &   1.39E+00        &     5.47E-02        &     5.47E$-$02        &     1.48E$-$04  \\
			\hline
			$\tau_{\text{quad}}$  & 8.28E$-$01        &     2.37E$-$02        &     1.62E$-$02        &     3.56E$-$05 \\
			\hline
			$\tau_{\text{saa}}$  &  2.00E+00        &     8.40E$-$03        &     2.30E$-$02        &     5.66E$-$05 \\
			\hline
		\end{tabular}
	\end{center}
\end{table}

\begin{table}[!htb]
	\caption{Estimates $\hat{q}$ of $q=(Q-Q(\bar{\zeta}))^2$ and mean squared errors (MSE) for $\hat{q}$, $q - T_1 q$, and $q - T_2 q$ based on 10 samples for a random design $\tau_{\text{random}}$, the deterministic optimal design $\tau_{\text{deter}}$, and the optimal design under uncertainty using the quadratic $\tau_{\text{quad}}$ and the sample average $\tau_{\text{saa}}$ approximations.}\label{tab:MSEq}
	\begin{center}
		\begin{tabular}{|c|c|c|c|c|}
			\hline
			design   &     $\hat{q}$       &     MSE($\hat{q}$)  &     MSE($q-T_1 q$) &     MSE$(q-T_2q$) \\
			\hline
			$\tau_{\text{random}} $ & 1.42E+02        &     4.83E+01        &     2.29E+00        &     3.24E$-$02 \\
			\hline 
			$\tau_{\text{deter}}$  &  2.48E+00        &     7.49E$-$01        &     7.49E$-$01        &     3.92E$-$03  \\
			\hline
			$\tau_{\text{quad}}$  &  9.22E$-$01        &     1.07E$-$01        &     9.04E$-$02        &     2.53E$-$04 \\
			\hline
			$\tau_{\text{saa}}$  &  4.08E+00        &     1.46E$-$01        &     2.48E$-$01        &     1.06E$-$03 \\
			\hline
		\end{tabular}
	\end{center}
\end{table}

To assess the accuracy of the Taylor approximation, we compute the mean squared errors (MSE) of the the design objective $Q$ and its residual using the linear and quadratic ($T_2$) Taylor approximations, as well as the quantity $q = (Q- Q(\bar{\zeta}))^2$ in the evaluation of the variance. The results are obtained at a random design, and the optimal design with deterministic, quadratic Taylor, and sample average approximations, and are shown in Table \ref{tab:MSEQ} and \ref{tab:MSEq}. These results indicate that the quadratic approximation is much more accurate than the linear approximation, both achieving errors smaller than $1\%$. We further remark that if higher accuracy is required, we can use the quadratic approximation as a control variate to reduce the variance in a sample average approximation, as introduced in \cite{ChenVillaGhattas19}.

\subsection{Scalability of the approximation and optimization methods}
\label{sec:scalability}
The random variable and the design variable are spatially distributed functions, whose dimensions can be very high after discretization. It is therefore crucial that the approximation and optimization are scalable with respect to both random and design variables. 
To illustrate the scalability of the approximation and optimization methods, we use a sequence of refined meshes as reported in Table \ref{tab:DOF}, which correspond to a sequence of increased dimensions for the random and design variables. 

\begin{table}[!htb]
	\caption{Degrees of freedom (DOF) for finite element discretization of the state variable $u$ and random variable $\zeta$ with piecewise linear elements (P1), and design variable $\tau$ with piecewise constant elements (P0), at a sequence of (uniformly refined) meshes, denoted by {meth1, meth2, meth3, meth4, meth5}.}\label{tab:DOF}
	\begin{center}
		\begin{tabular}{|c|c|c|c|c|c|}
			\hline
			DOF    & mesh1 & mesh2 & mesh3 & mesh4 & mesh5 \\
			\hline
			$u$(P1)	  &22,110& 86,788& 345,606& 1,373,814&5,488,216\\
			\hline
			$\zeta$(P1)	 & 2,347& 8,795& 34,217& 134,796& 535,321\\
			\hline
			$\tau$(P0)	 & 4,454& 17,114& 67,462& 267,640& 1,066,761\\
			\hline
		\end{tabular}
		\label{table:discretization}
	\end{center}
\end{table}

\begin{figure}[!htb]
	\begin{center}
		\includegraphics[scale=0.38]{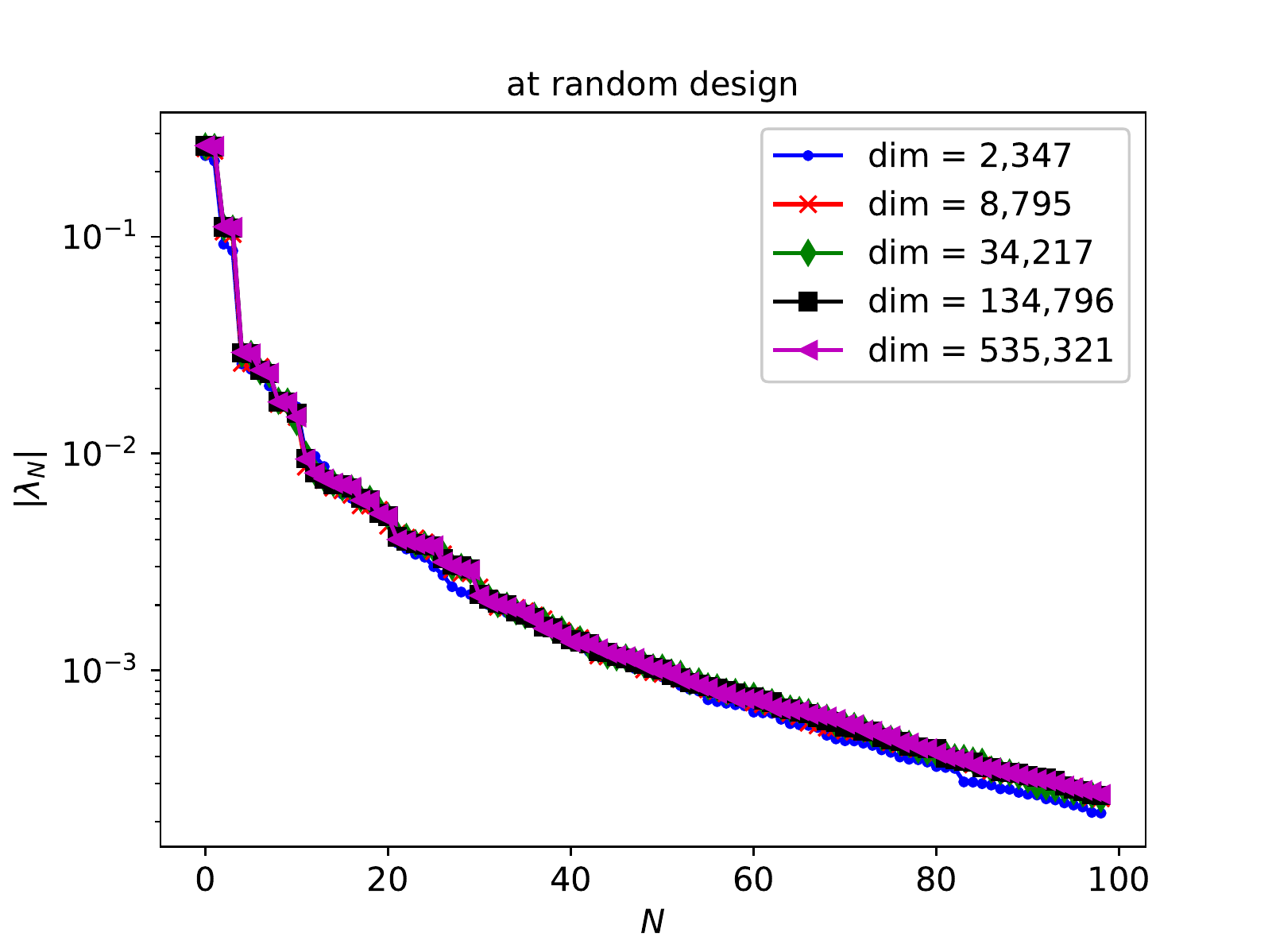}
		\includegraphics[scale=0.38]{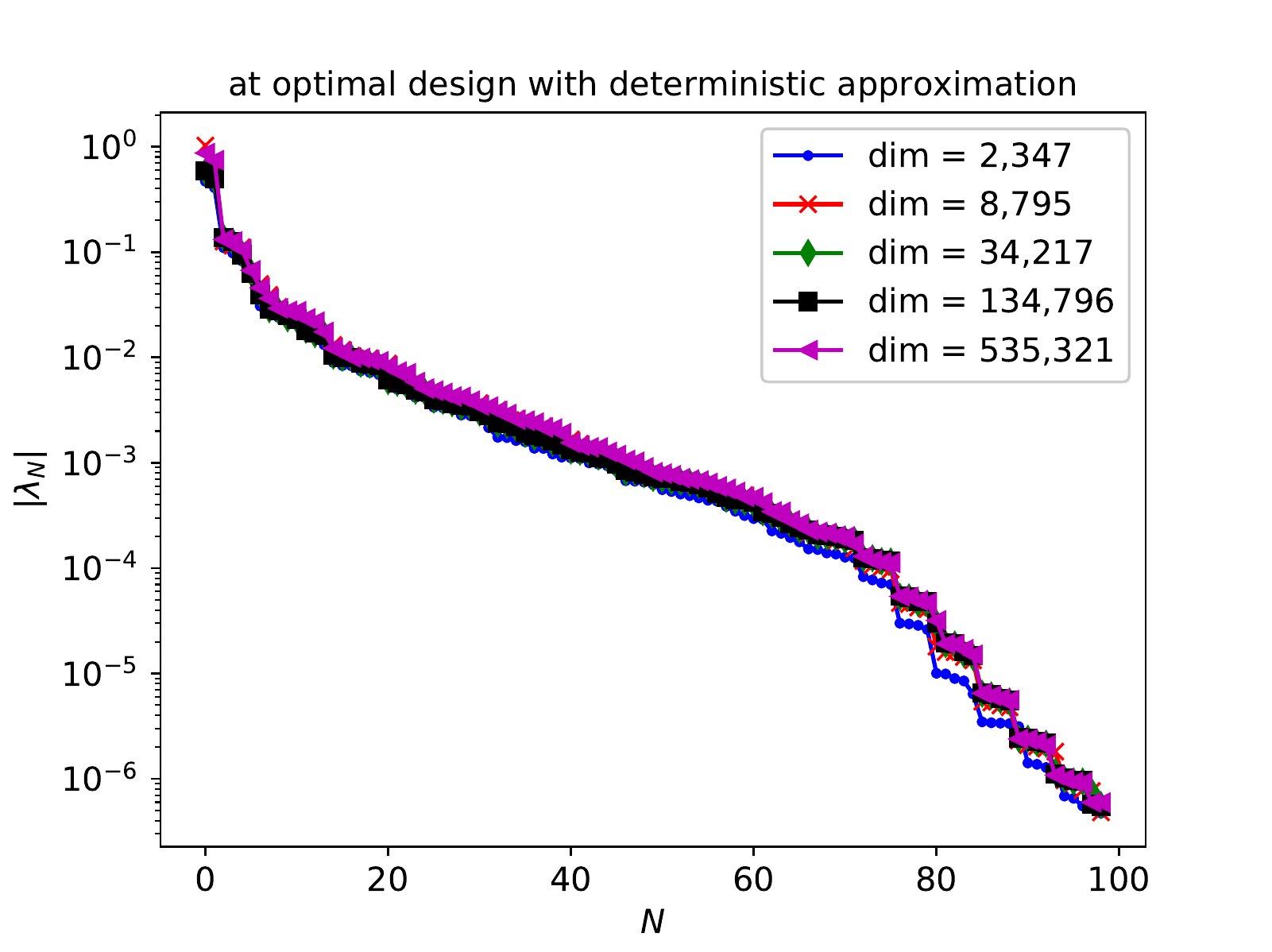}
		
		\includegraphics[scale=0.38]{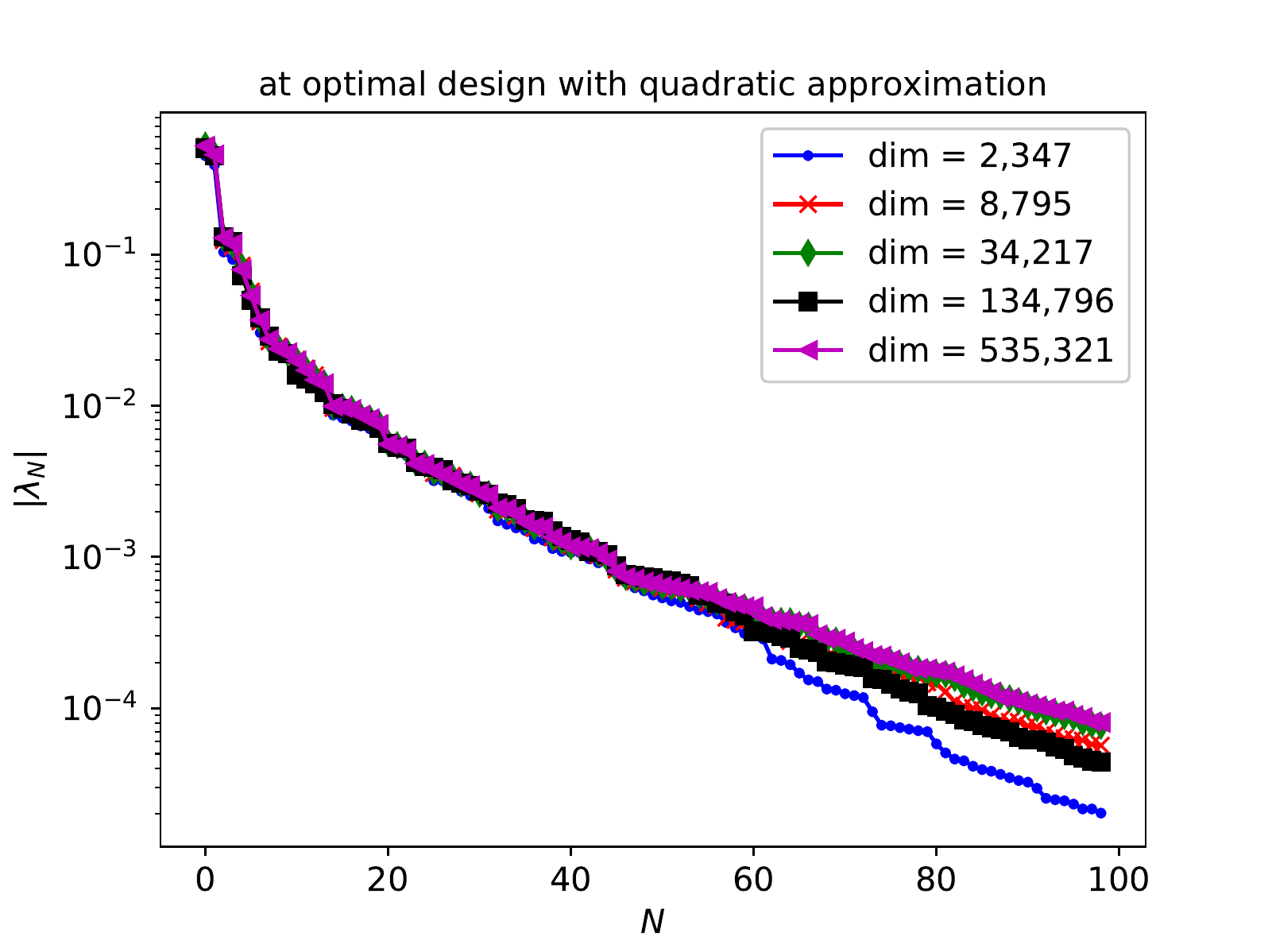}
		\includegraphics[scale=0.38]{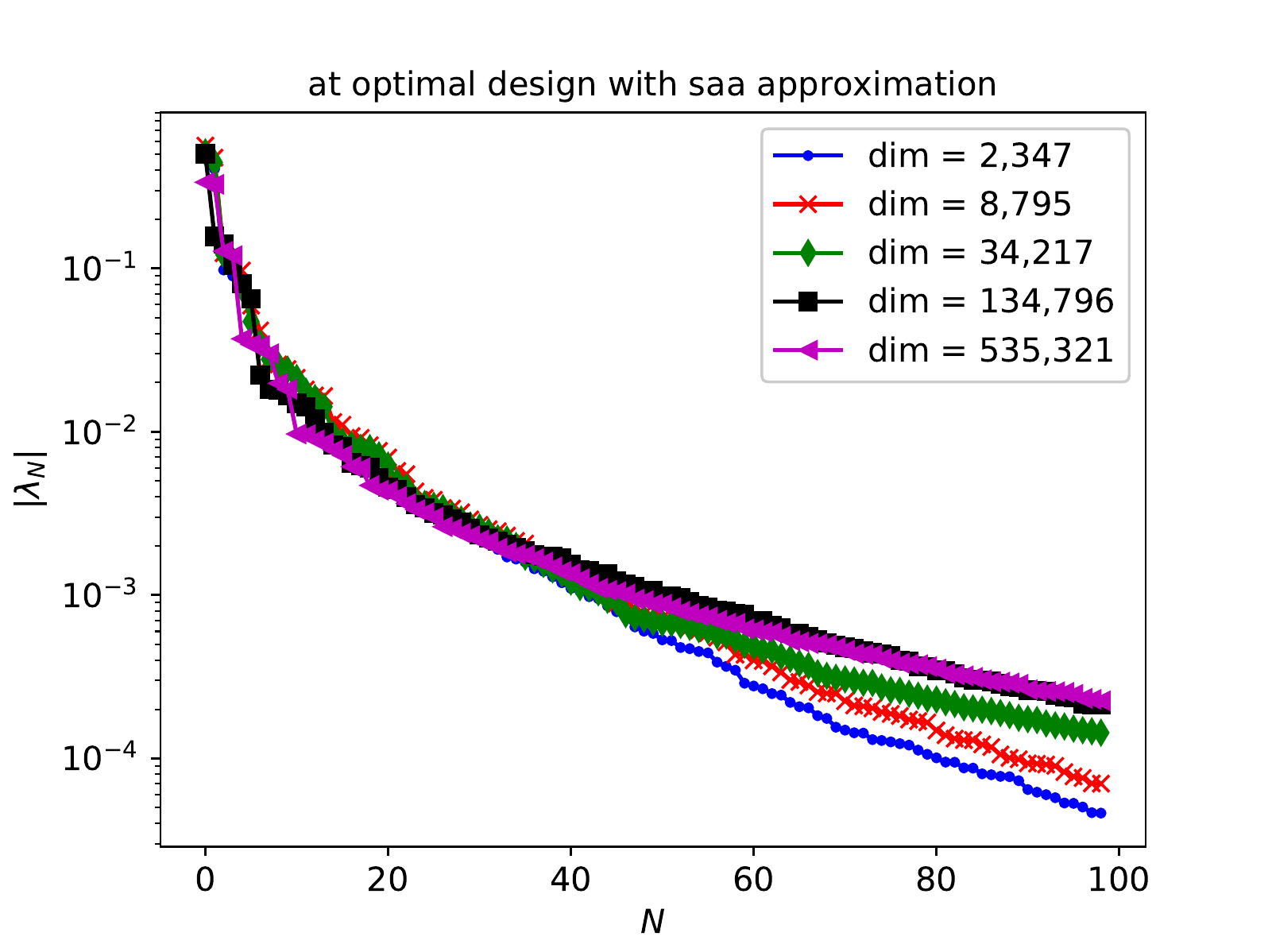}
	\end{center}
	\caption{Decay of the absolute generalized eigenvalues of the covariance preconditioned Hessian in \eqref{eq:generalizedEigenProblem} at different designs. A design at a realization of space white noise (top-left); the optimal design obtained with deterministic (top-right), quadratic (bottom-left), and sample average approximations (bottom-right).}\label{fig:trace}
\end{figure}

As shown in Fig.\ \ref{fig:trace}, the scalability with respect to the complexity of the quadratic approximation is implied by the similar decay pattern of the absolute eigenvalues of the generalized eigenvalue problem \eqref{eq:generalizedEigenProblem} across the refined meshes, which determines the accuracy of the trace estimate. Moreover, the accuracy of the quadratic approximation measured by the mean squared errors is reported in Table \ref{tab:Qmse} and \ref{tab:qmse}, which remains about $1\%$  with increasing dimensions, and indicates that the accuracy of the quadratic approximation is also scalable.

\begin{table}[!htb]
	\caption{Estimates $\hat{Q}$ of the design objective $Q$ and mean squared errors (MSE) for $\hat{Q}$, $Q - T_1 Q$, and $Q-T_2 Q$, based on 10 samples for different parameter dimensions. }\label{tab:Qmse}
	\begin{center}
		\begin{tabular}{|c|c|c|c|c|}
			\hline
			dimension   &     $\hat{Q}$       &     MSE($\hat{Q}$)  &     MSE($Q-T_1Q$) &     MSE$(Q-T_2Q$) \\
			\hline
			2,347 &  6.49E$-$01        &     1.28E$-$02        &     8.92E$-$03        &     1.01E$-$04 \\
			\hline 
			8,795  &  7.66E$-$01        &     1.66E$-$02        &     1.07E$-$02        &     1.54E$-$04  \\
			\hline
			34,217 & 8.28E$-$01        &     2.37E$-$02        &     1.62E$-$02        &     3.56E$-$05 \\
			\hline
		\end{tabular}
	\end{center}
\end{table}

\begin{table}[!htb]
	\caption{Estimates $\hat{q}$ of $q=(Q-Q(\bar{\zeta}))^2$ and mean squares errors (MSE) $\hat{q}$, $q - T_1 q$, and $q - T_2 q$, based on 10 samples for different parameter dimensions.}\label{tab:qmse}
	\begin{center}
		\begin{tabular}{|c|c|c|c|c|}
			\hline
			dimension   &     $\hat{q}$       &     MSE($\hat{q}$)  &     MSE($q-T_1q$) &     MSE$(q-T_2q$) \\
			\hline
			2,347 & 5.49E$-$01        &     3.42E$-$02        &     3.28E$-$02        &     5.61E$-$04    \\
			\hline 
			8,795  &  7.54E$-$01        &     9.00E$-$02        &     7.39E$-$02        &     1.64E$-$03  \\
			\hline
			34,217  &  9.22E$-$01        &     1.07E$-$01        &     9.04E$-$02        &     2.53E$-$04   \\
			\hline
		\end{tabular}
	\end{center}
\end{table}

\begin{figure}[!htb]
	\begin{center}
		\includegraphics[scale=0.38]{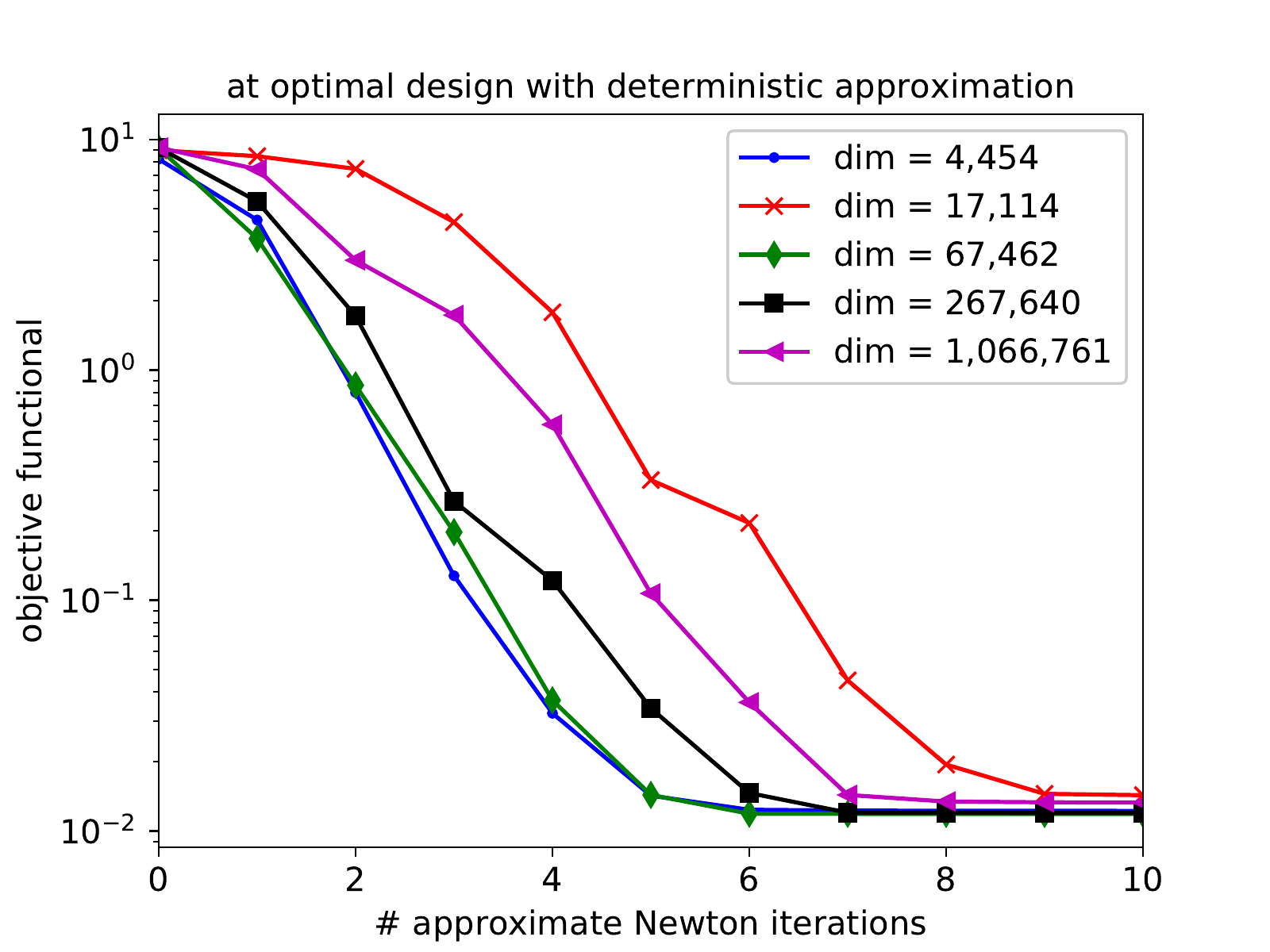}
		\includegraphics[scale=0.38]{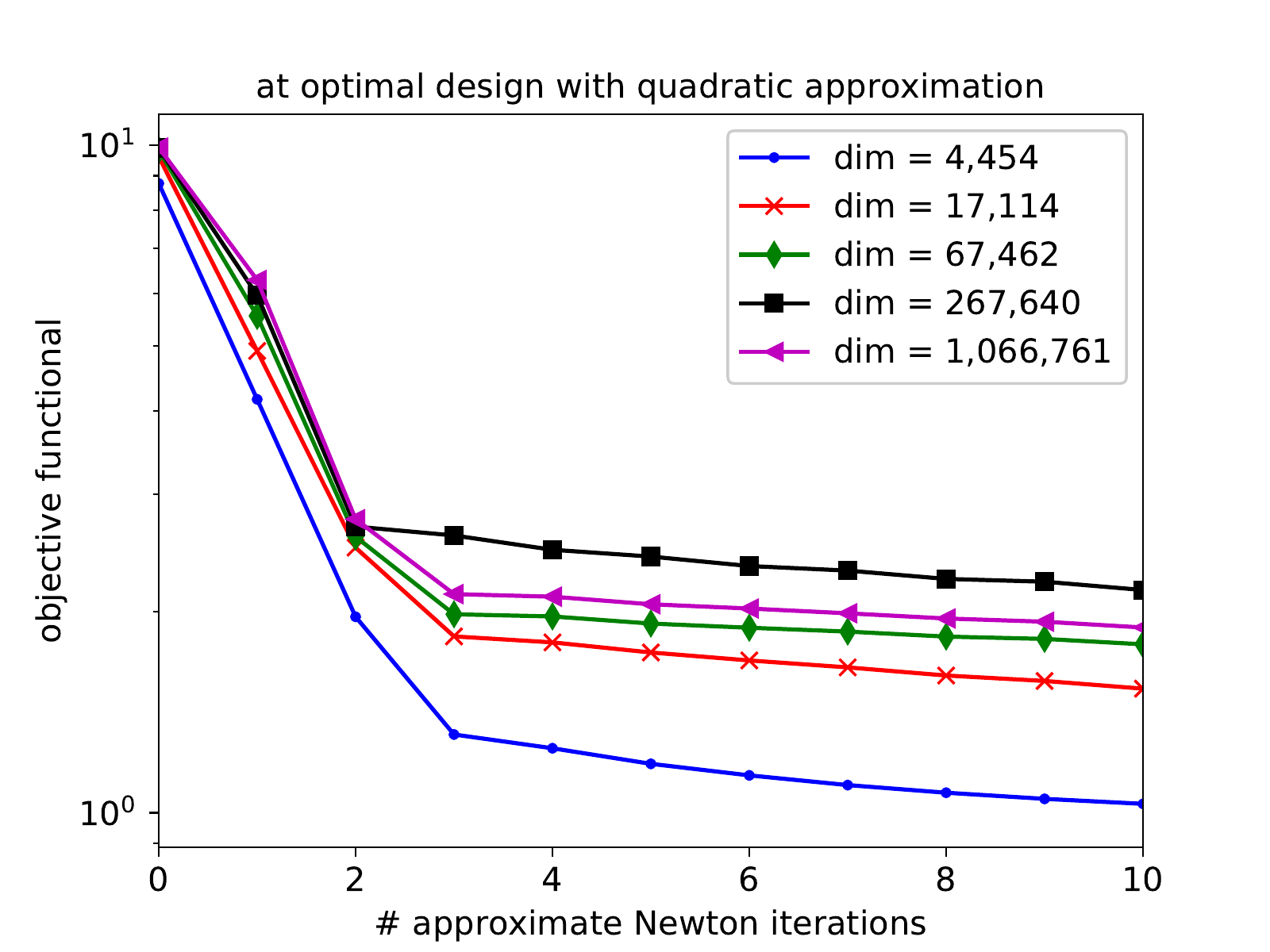}
		\includegraphics[scale=0.38]{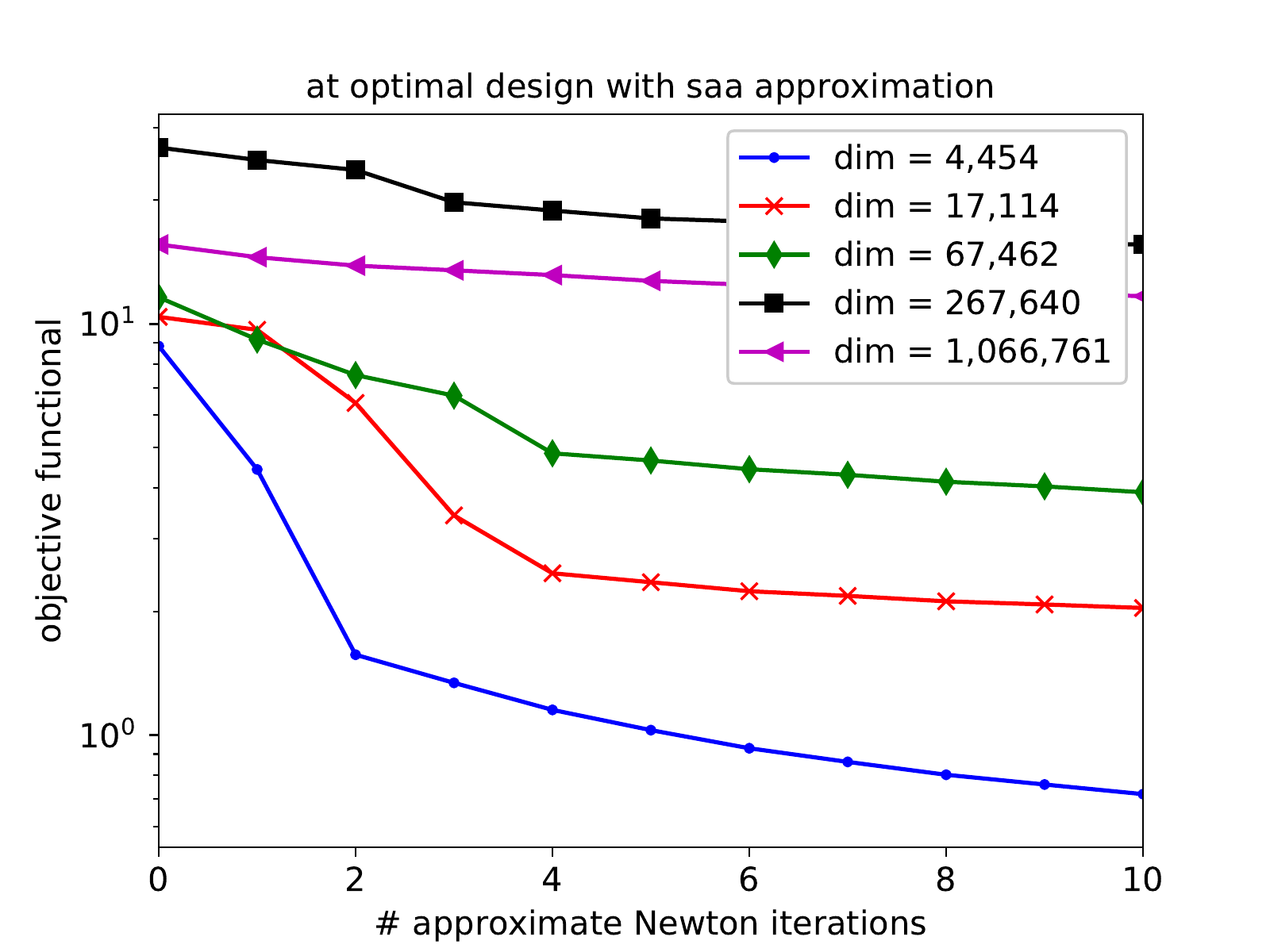}
	\end{center}
	\caption{Decay of the objective functional with the number of approximate Newton optimization steps for deterministic approximation (top-left), quadratic approximation (top-right), and sample average approximation (bottom).}\label{fig:quasiNewton}
\end{figure}

As for the scalability of the approximate Newton optimization
algorithm, we plot the decay of the objective functional against the
number of optimization iterations in Fig.\ \ref{fig:quasiNewton}.
Fast and relatively mesh-independent decay of the objective functional
can be observed for the deterministic approximation, which is
understandable since the Hessian approximation in Section
\ref{sec:gradDeterministic} is in fact exact in this case, so that the
method is a proper Newton method. For the quadratic approximation,
convergence is only weakly dependent on the discretization (with
sufficient mesh resolution); thus the use of the deterministic Hessian
in place of the true Hessian still results in a relatively scalable
number of optimization iterations. In contrast, the use of the
deterministic Hessian for the sample average approximation does not
yield a scalable method in this case, as shown by the dependence of
the iterations on mesh size and the resulting problem dimension.  

\subsection{Multiple directions and multiple frequencies}

In this numerical experiment, we access the ability of the optimal cloak to hide the obstacle from the incident wave from multiple attack angles and multiple frequencies. In the test, for the incident wave $e^{i k x \cdot b}$ we choose four attack angles, $b = (1,0), (0,1), (0, -1), (-1, 0)$, and four frequencies $k= k_0/2, 2k_0/3, 5k_0/6, k_0$, and set three test trials. In the first trial, we use four directions at one frequency $k = k_0$; in the second trial, we use four frequencies at one direction $b = (1, 0)$; in the third trial, we use four directions at four frequencies $b = (1,0), (0,1), (0, -1), (-1, 0)$ and $k= k_0/2, 2k_0/3, 5k_0/6, k_0$.

\begin{figure}
	\begin{center}
		\includegraphics[scale=0.45]{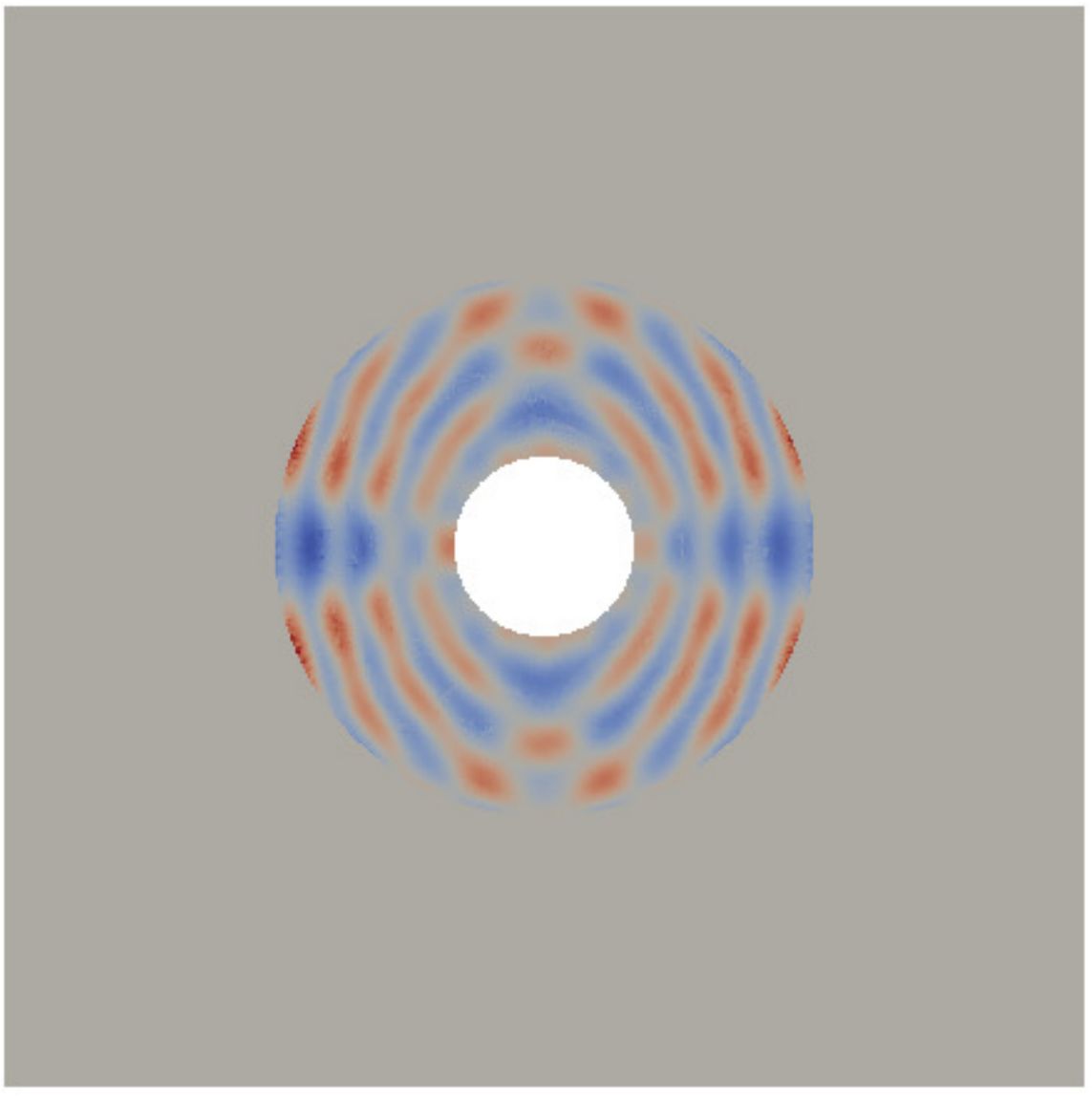}
		\includegraphics[scale=0.45]{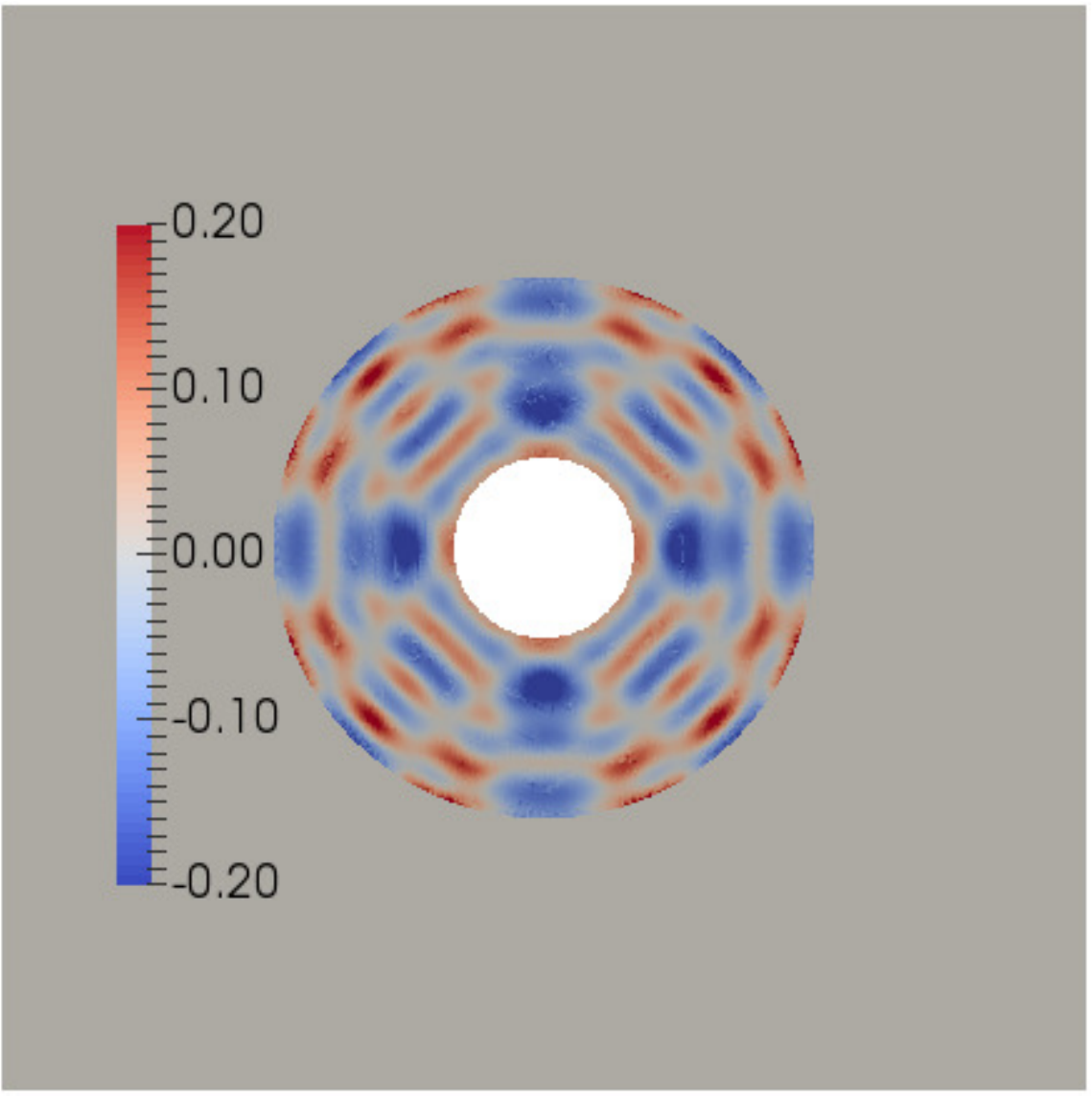}
		\includegraphics[scale=0.45]{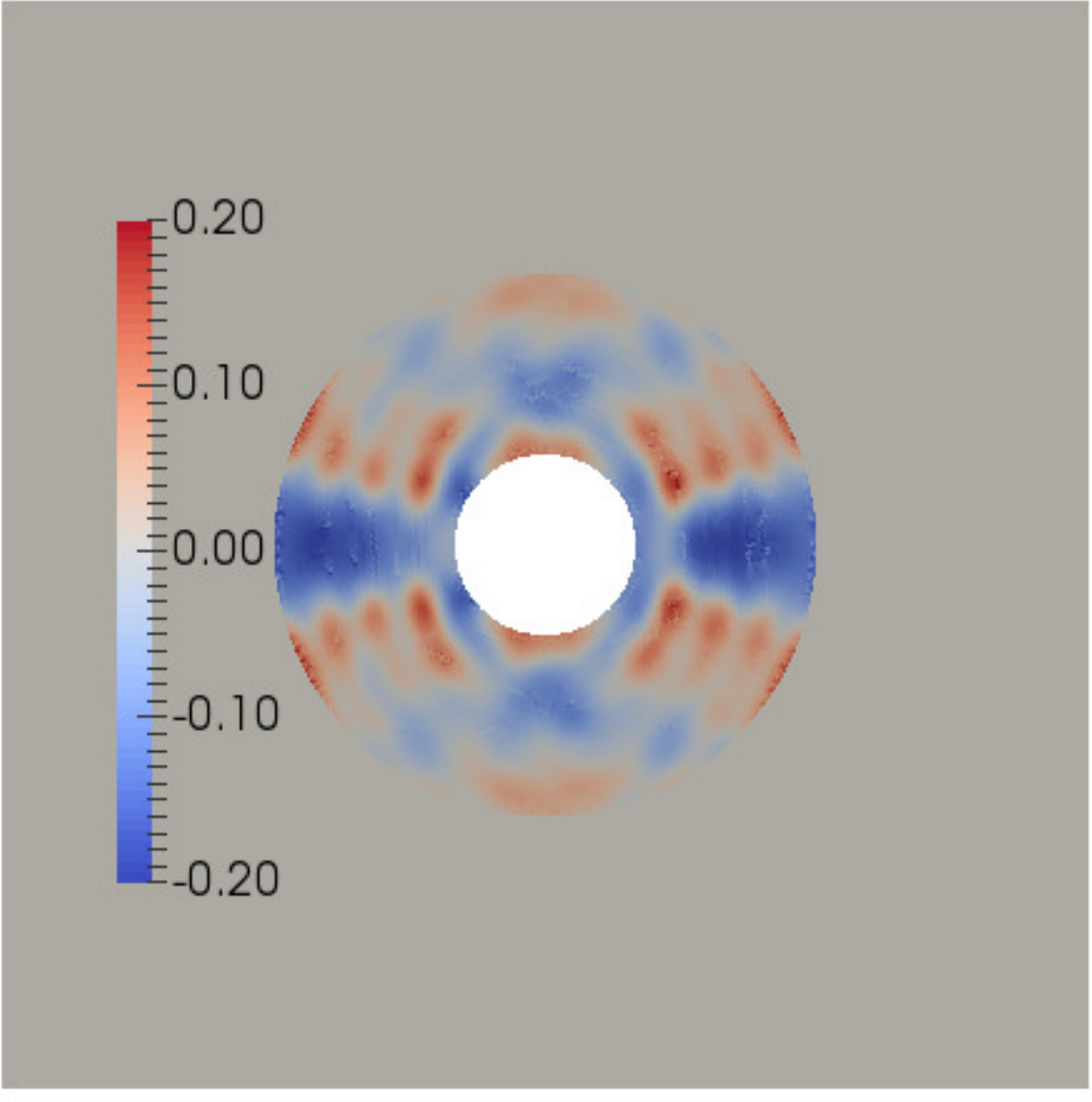}
		\includegraphics[scale=0.336]{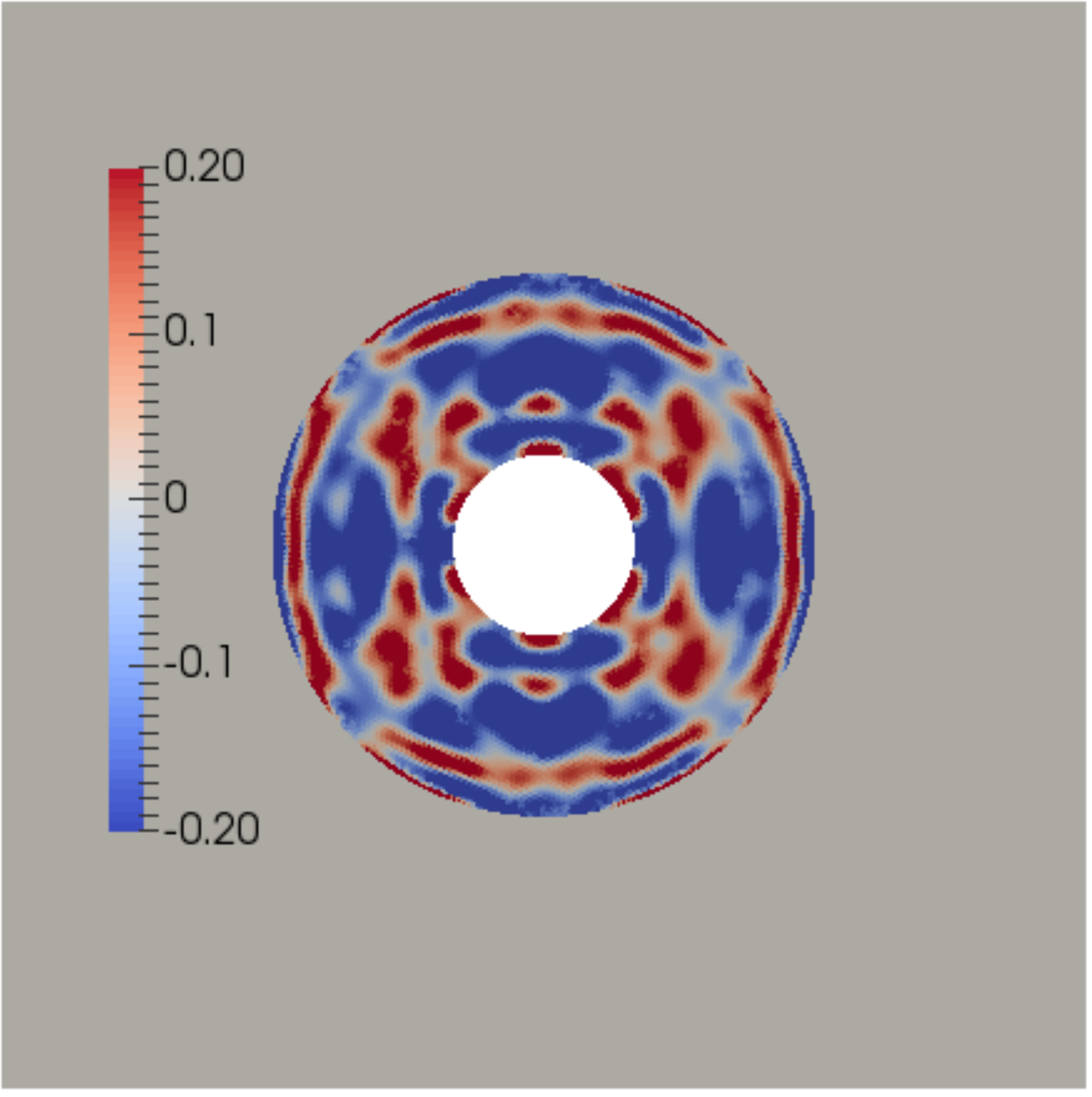}
	\end{center}
	\caption{Optimal designs under uncertainty using quadratic approximation. One direction and one frequency (top-left), which is the same as in Fig.\ \ref{fig:designs}, four directions and one frequency (top-right), one direction and four frequencies (bottom-left), and four directions and four frequencies (bottom-right).}\label{fig:4designs}
\end{figure}

\begin{figure}
	\begin{center}
		\includegraphics[scale=0.28]{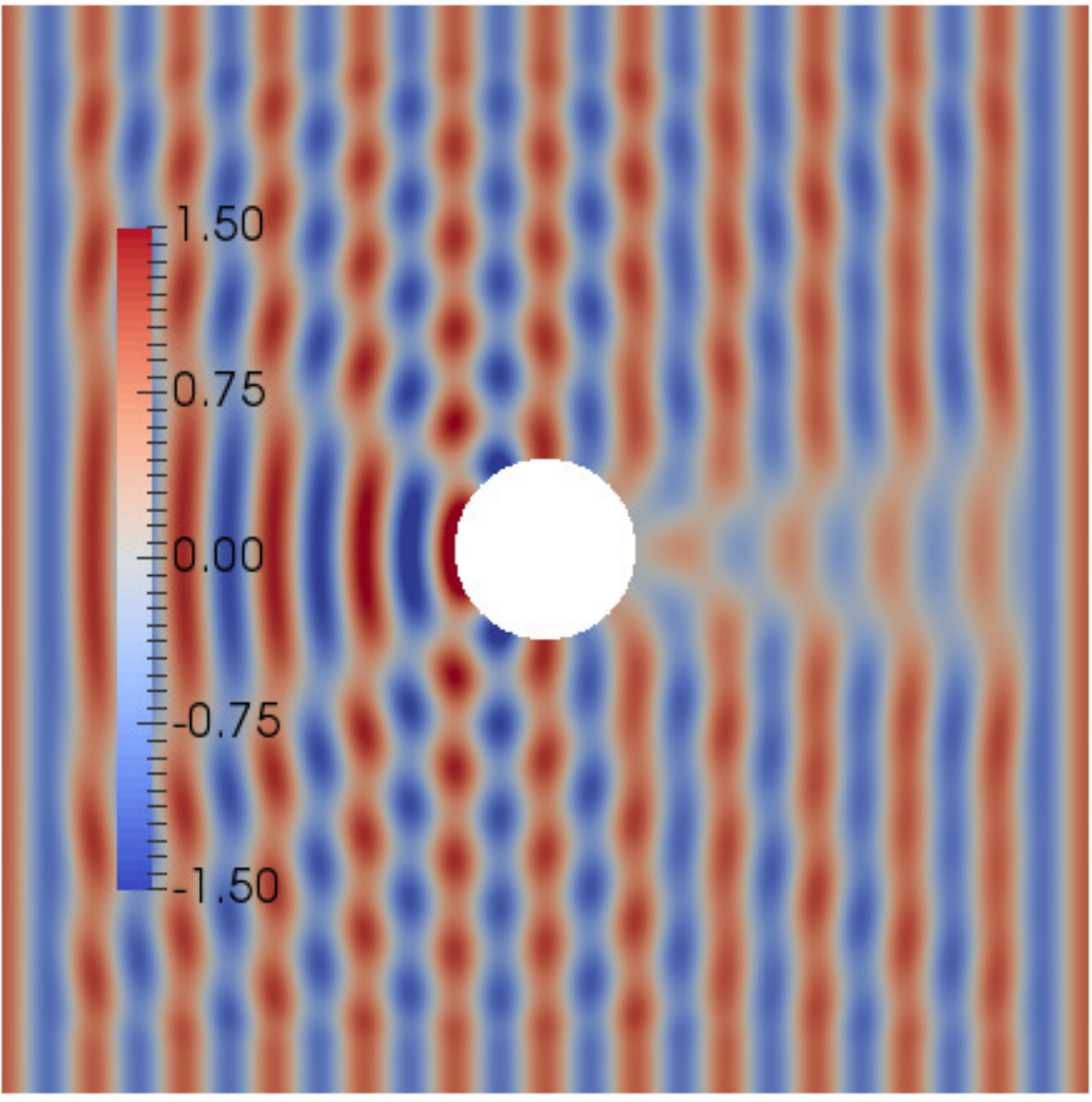}
		\includegraphics[scale=0.28]{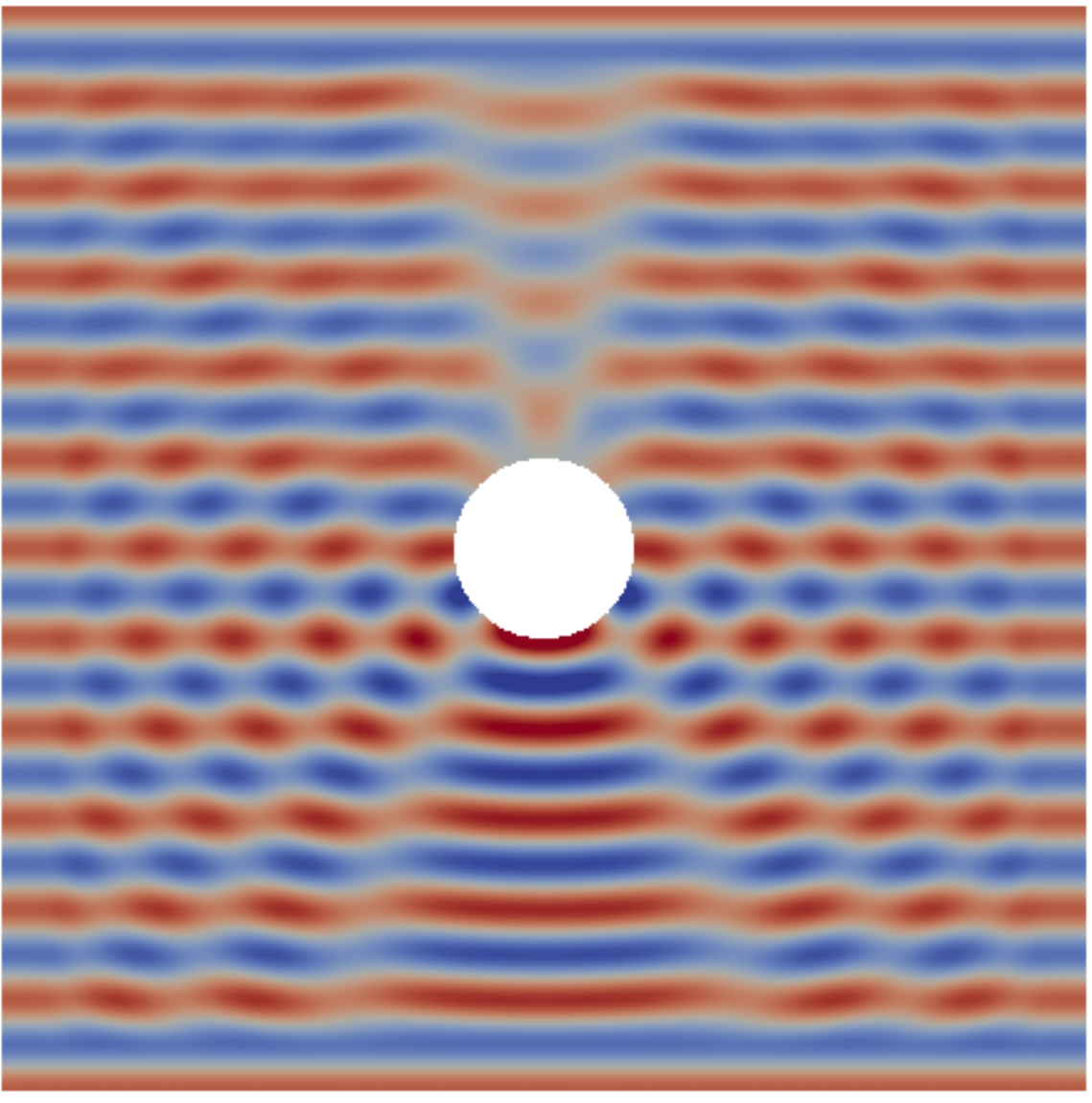}
		\includegraphics[scale=0.28]{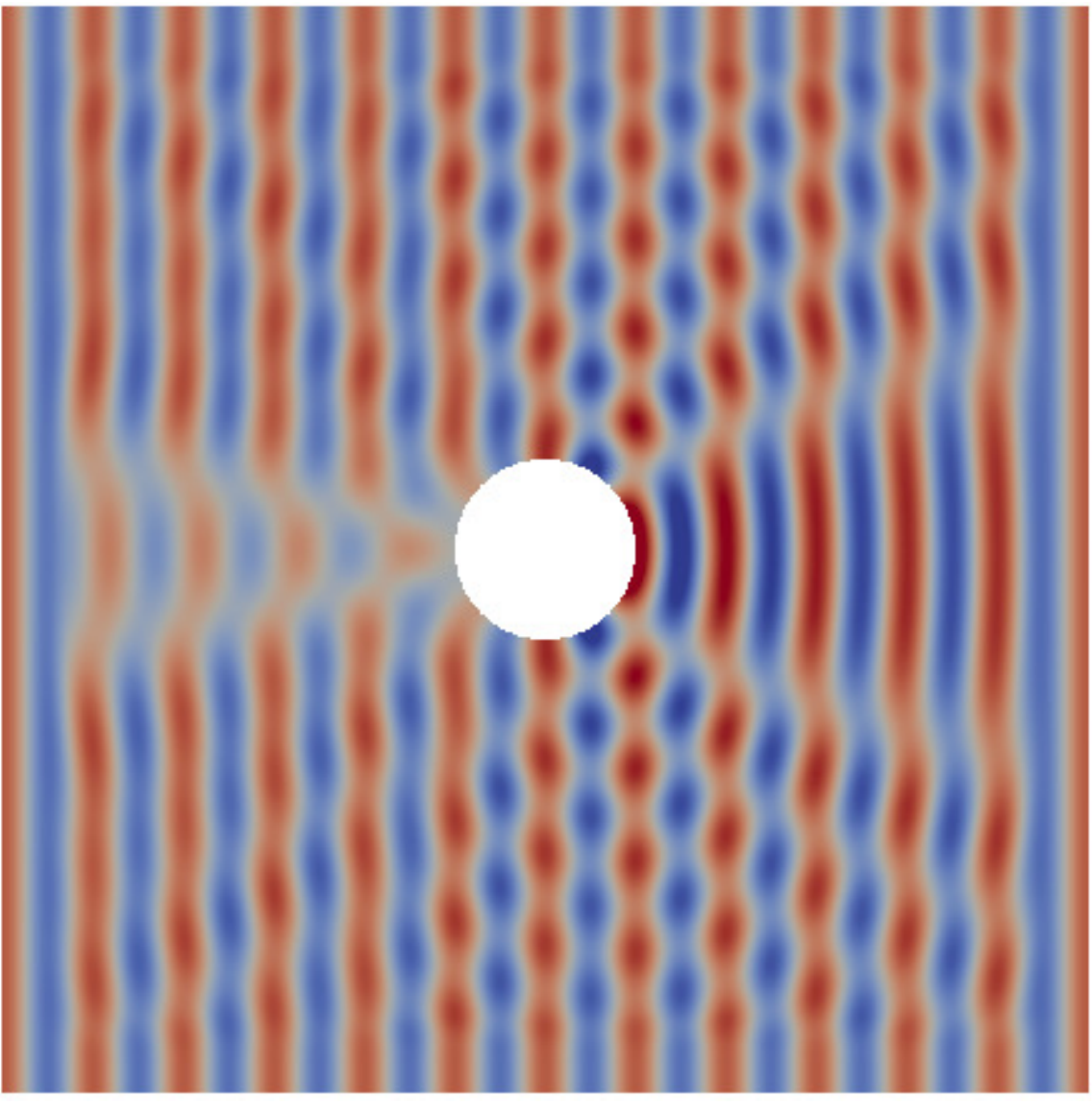}
		\includegraphics[scale=0.28]{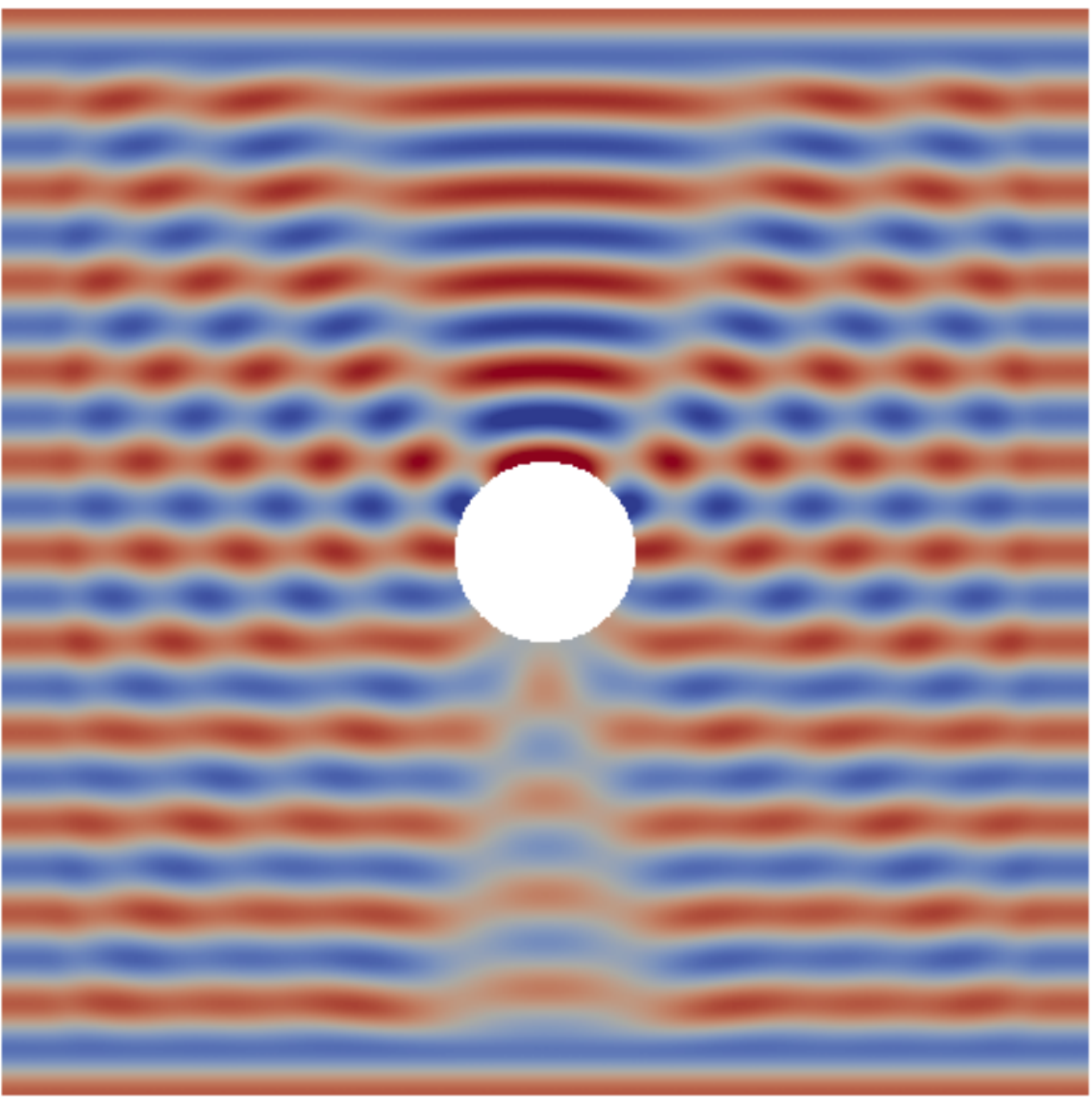}
		\includegraphics[scale=0.28]{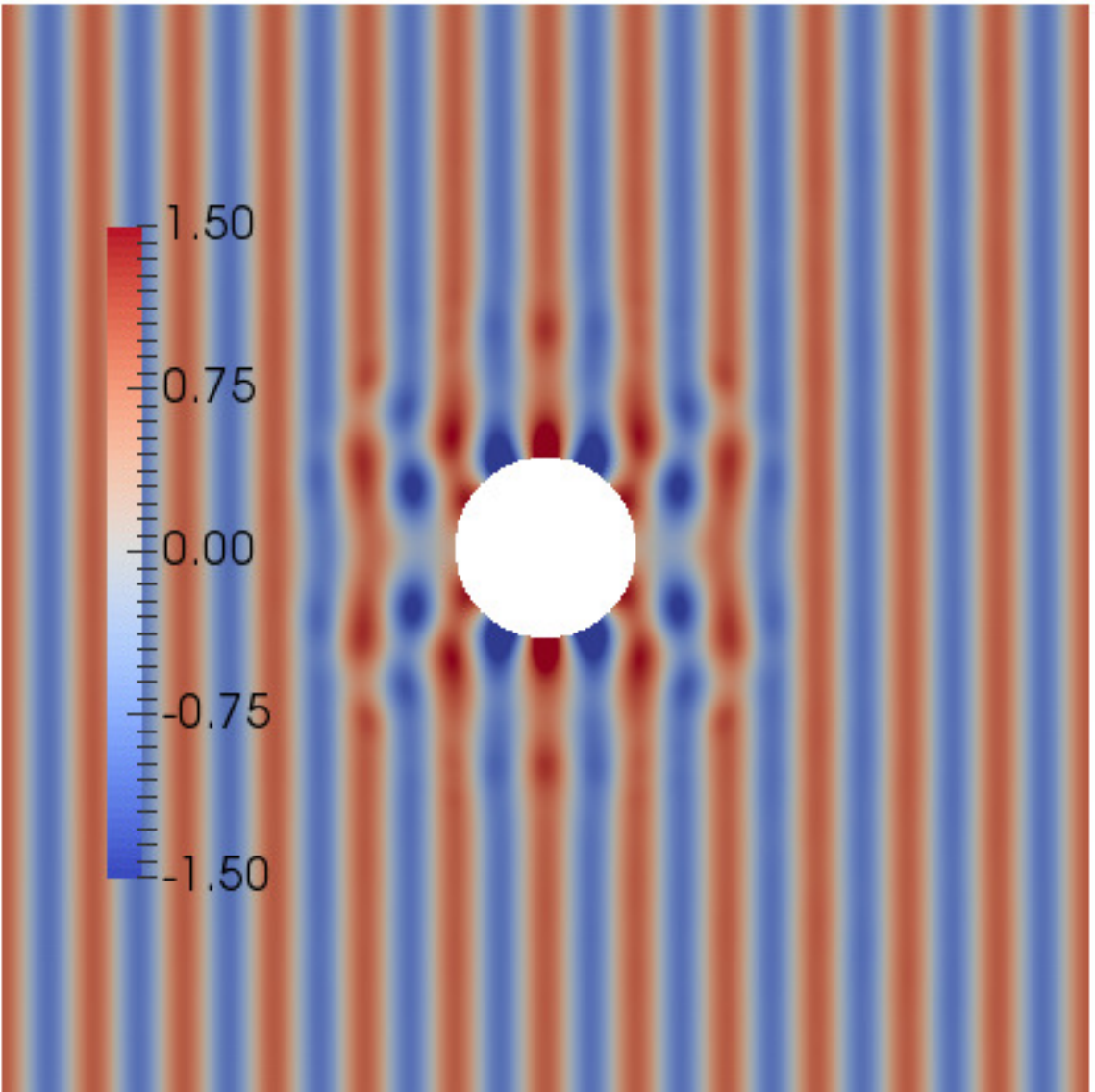}
		\includegraphics[scale=0.28]{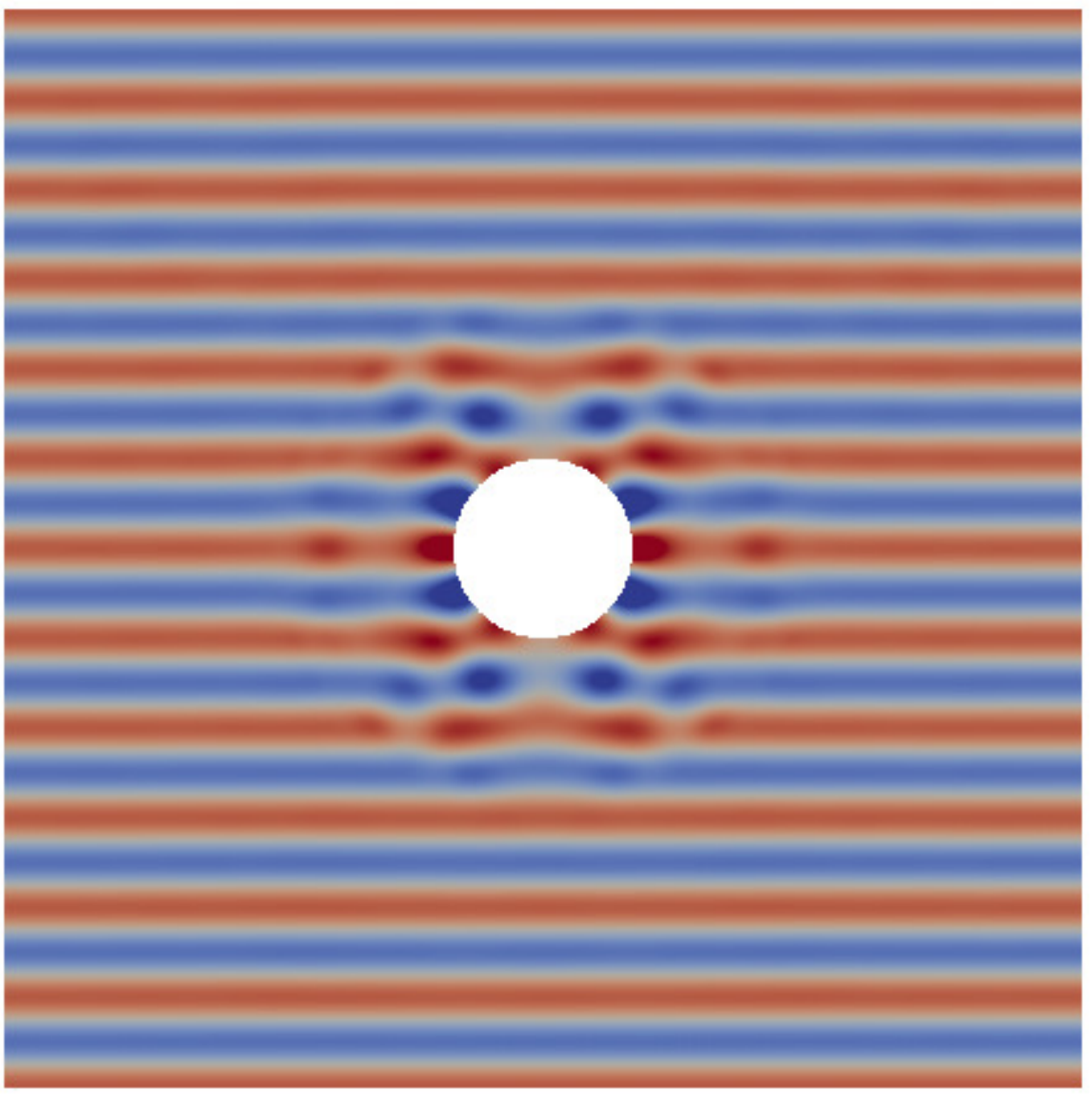}
		\includegraphics[scale=0.28]{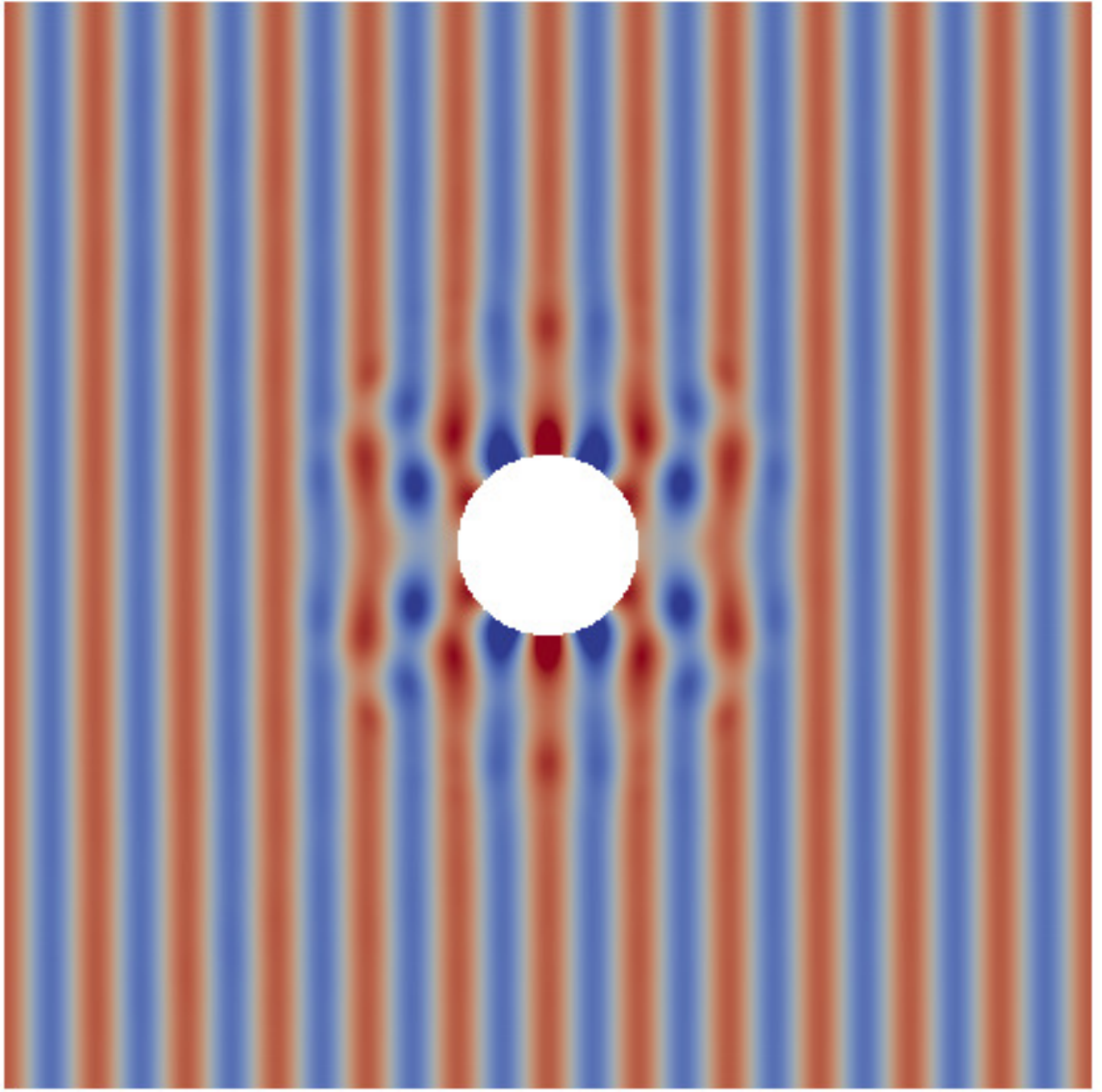}
		\includegraphics[scale=0.28]{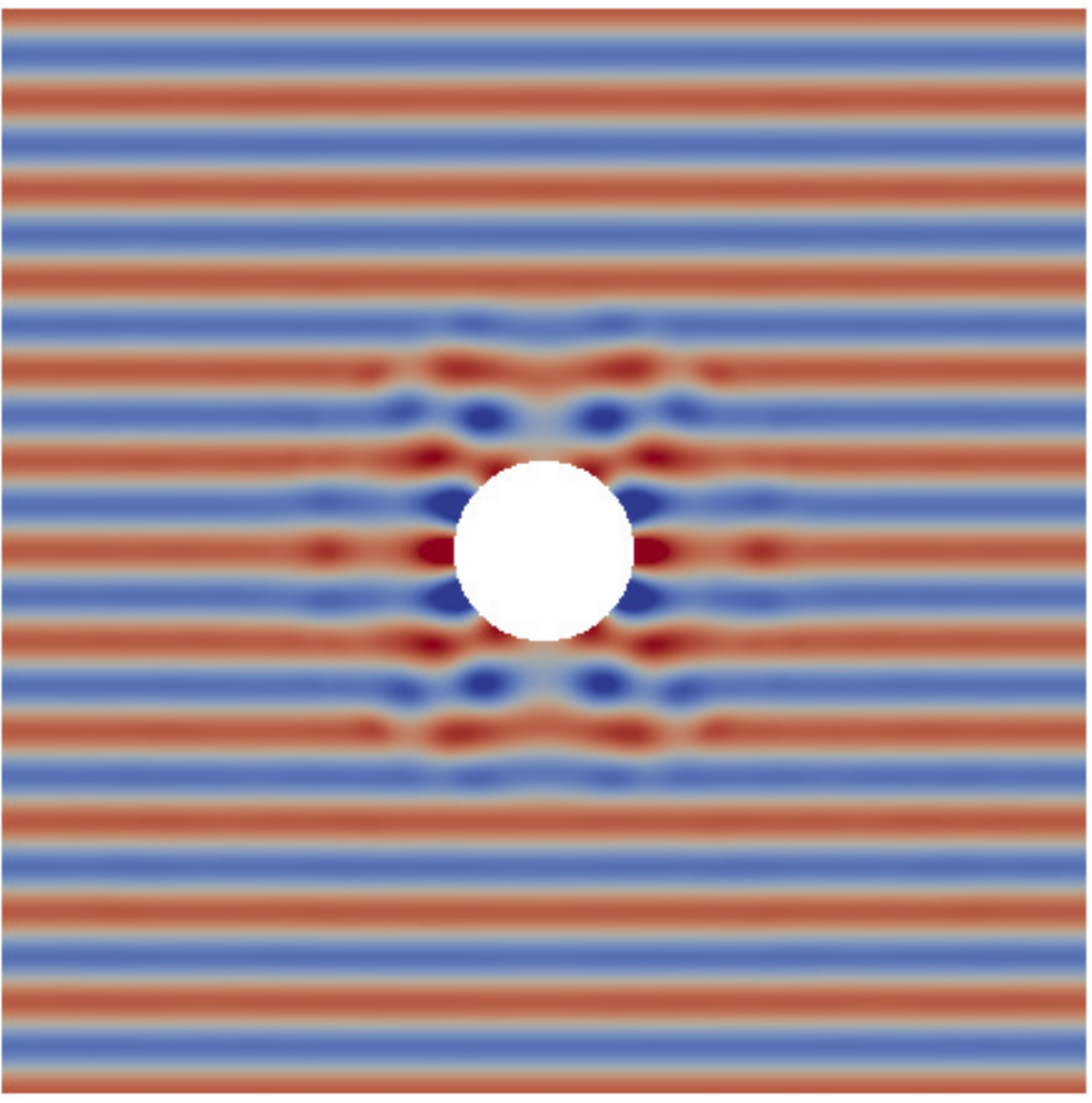}
	\end{center}
	\caption{The real part of the total wave fields without (top) and with (bottom) the cloak designed under uncertainty for the case of four directions and one frequency.}\label{fig:4sources}
\end{figure}

\begin{figure}
	\begin{center}
		\includegraphics[scale=0.28]{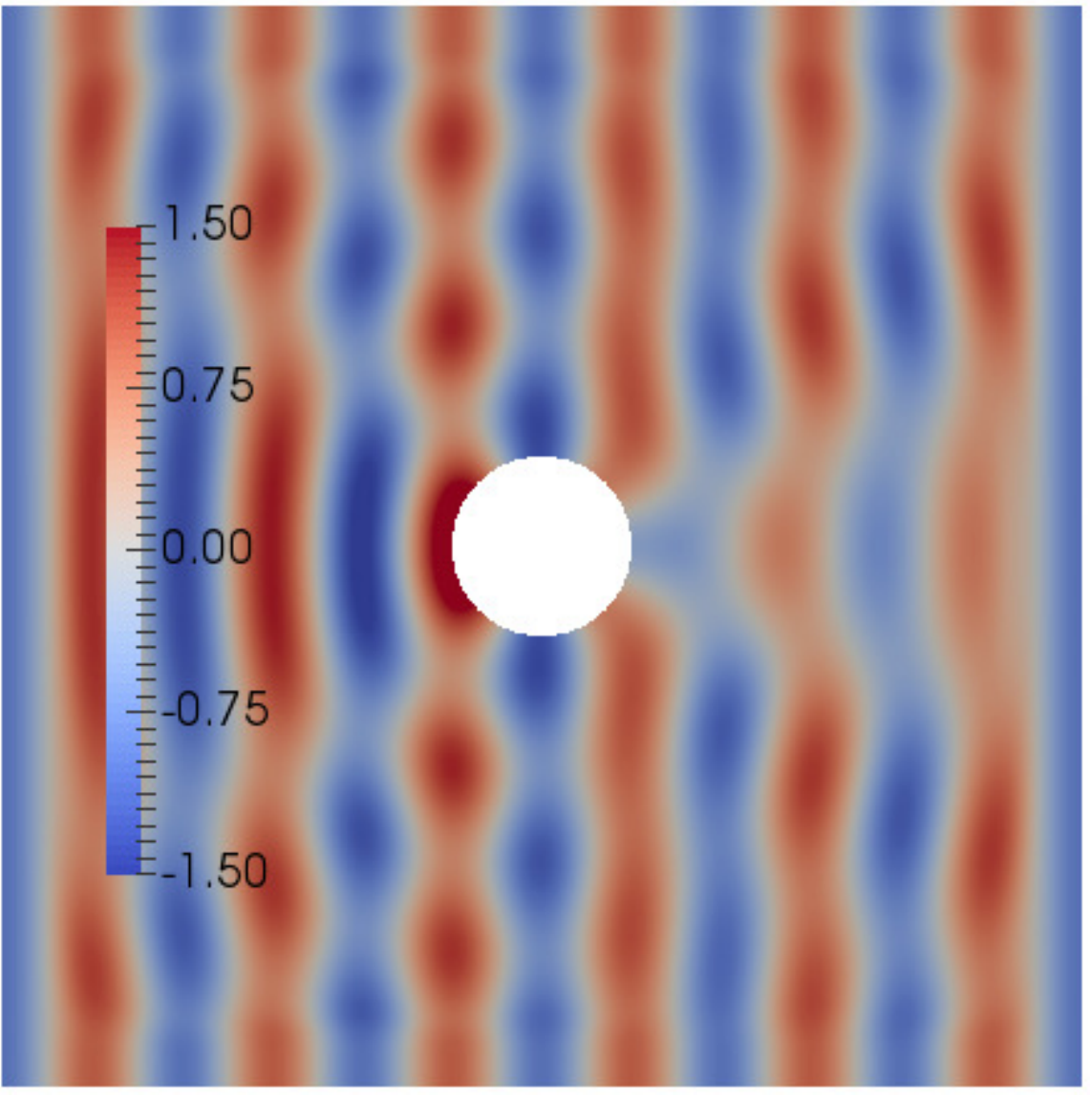}
		\includegraphics[scale=0.28]{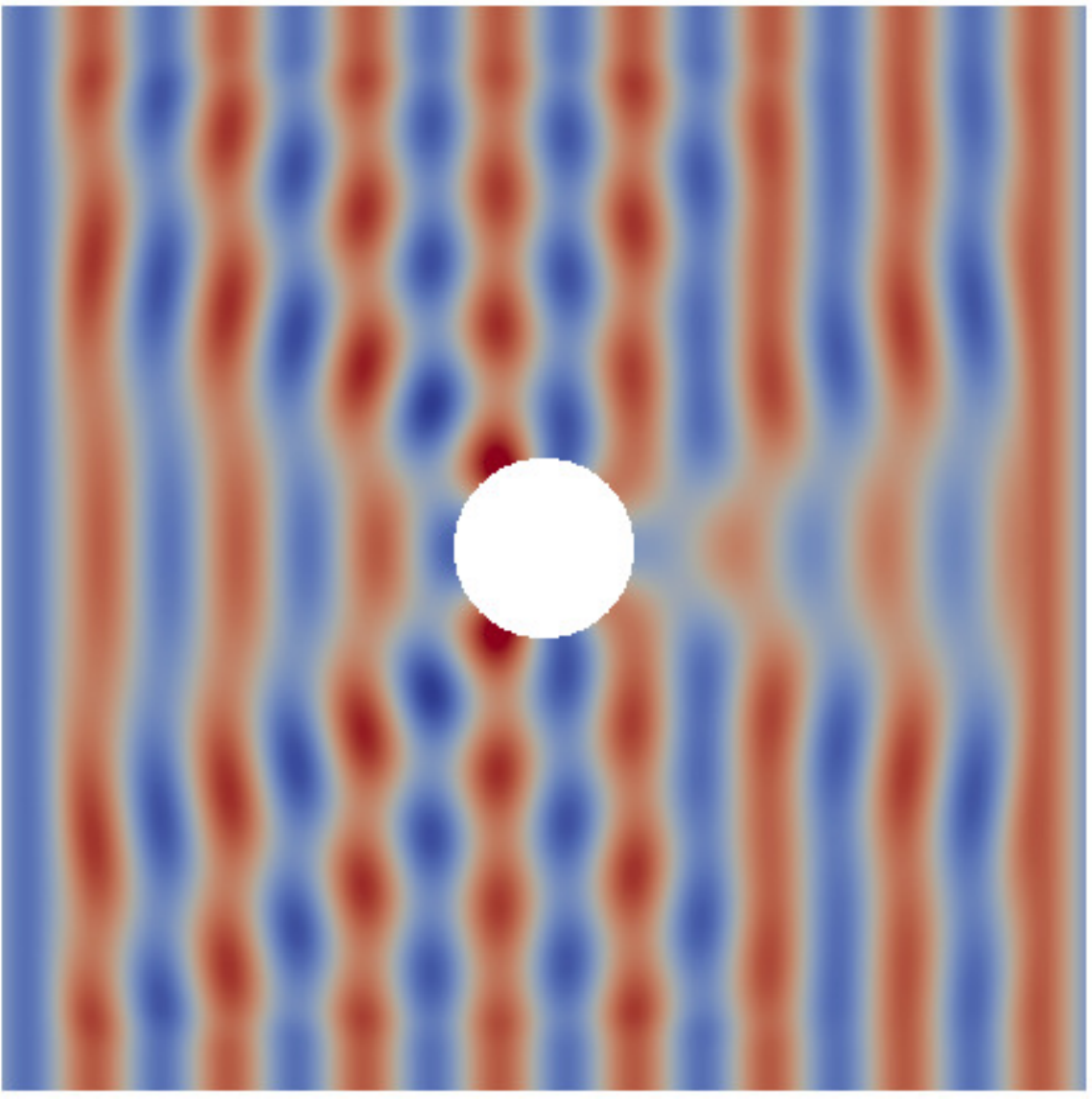}
		\includegraphics[scale=0.28]{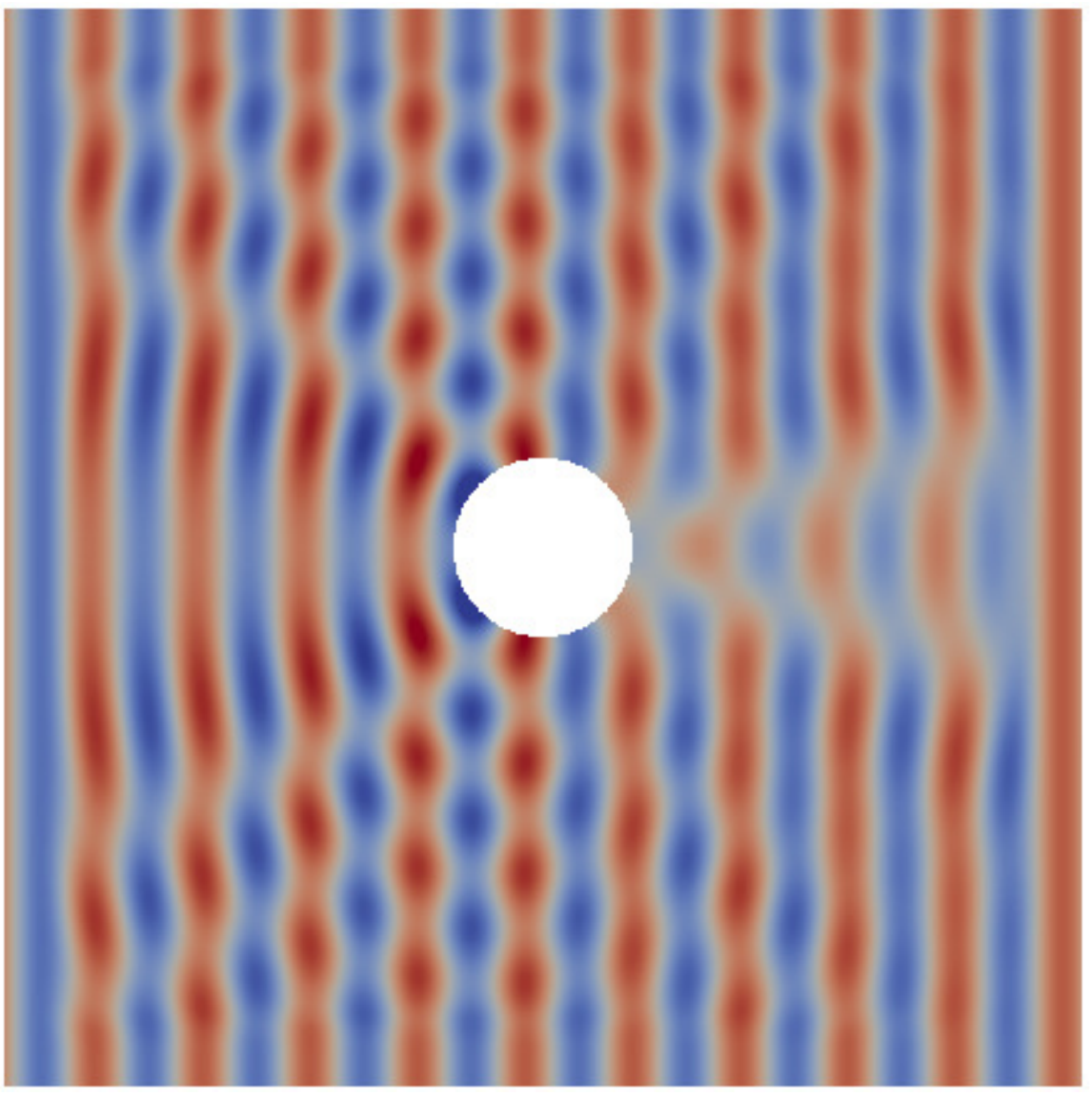}
		\includegraphics[scale=0.28]{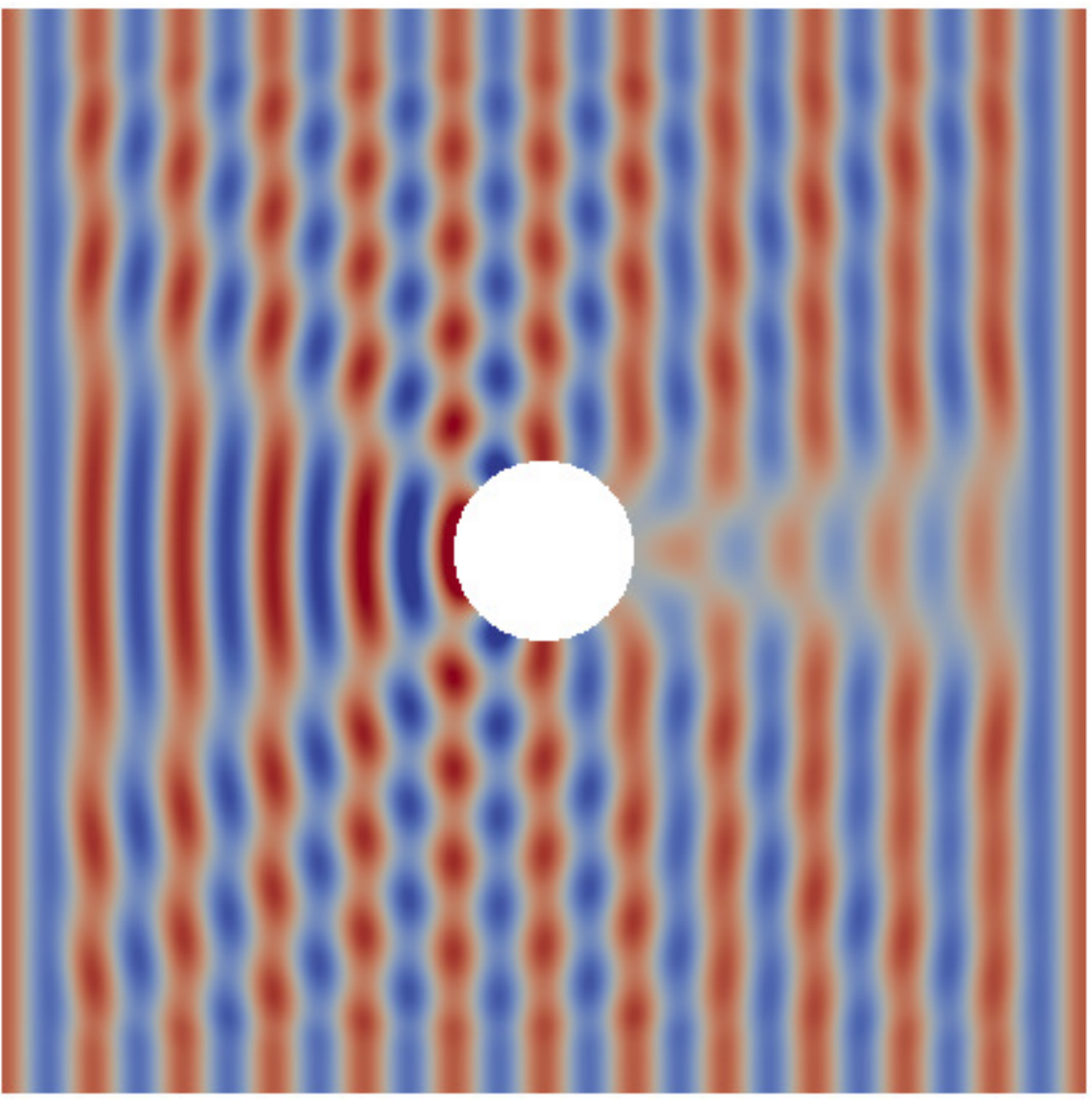}
		
		\includegraphics[scale=0.21]{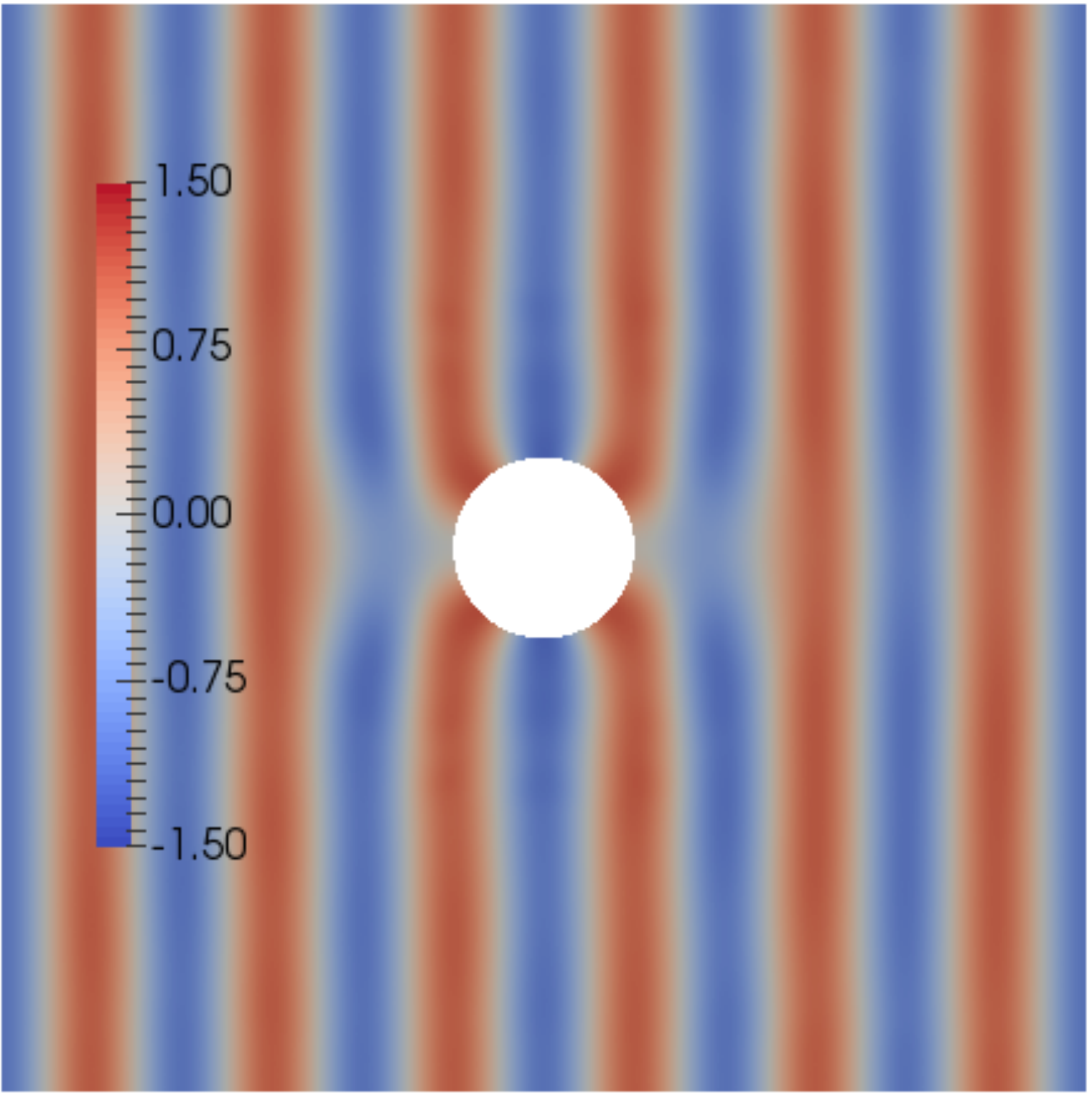}
		\includegraphics[scale=0.21]{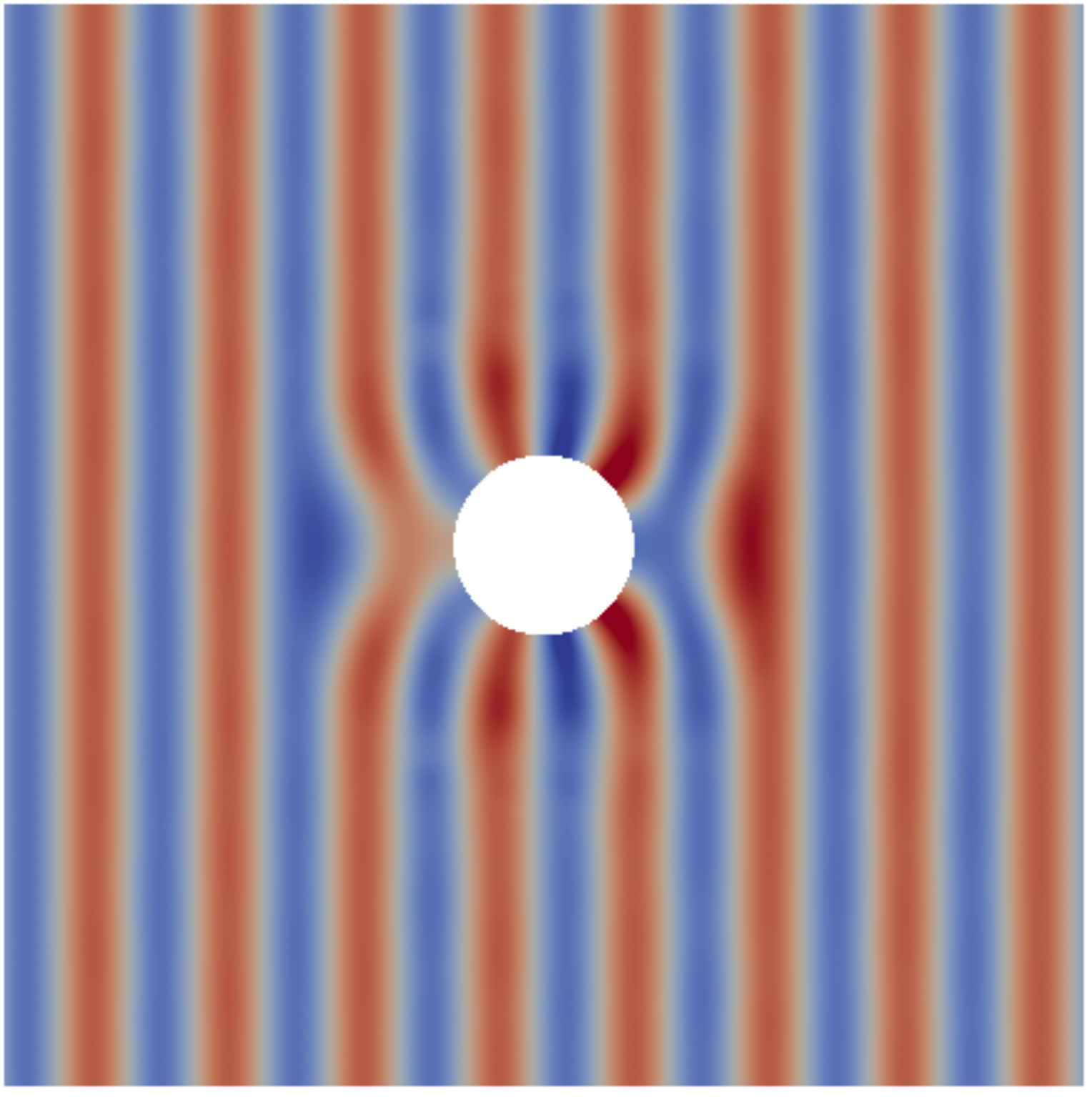}
		\includegraphics[scale=0.21]{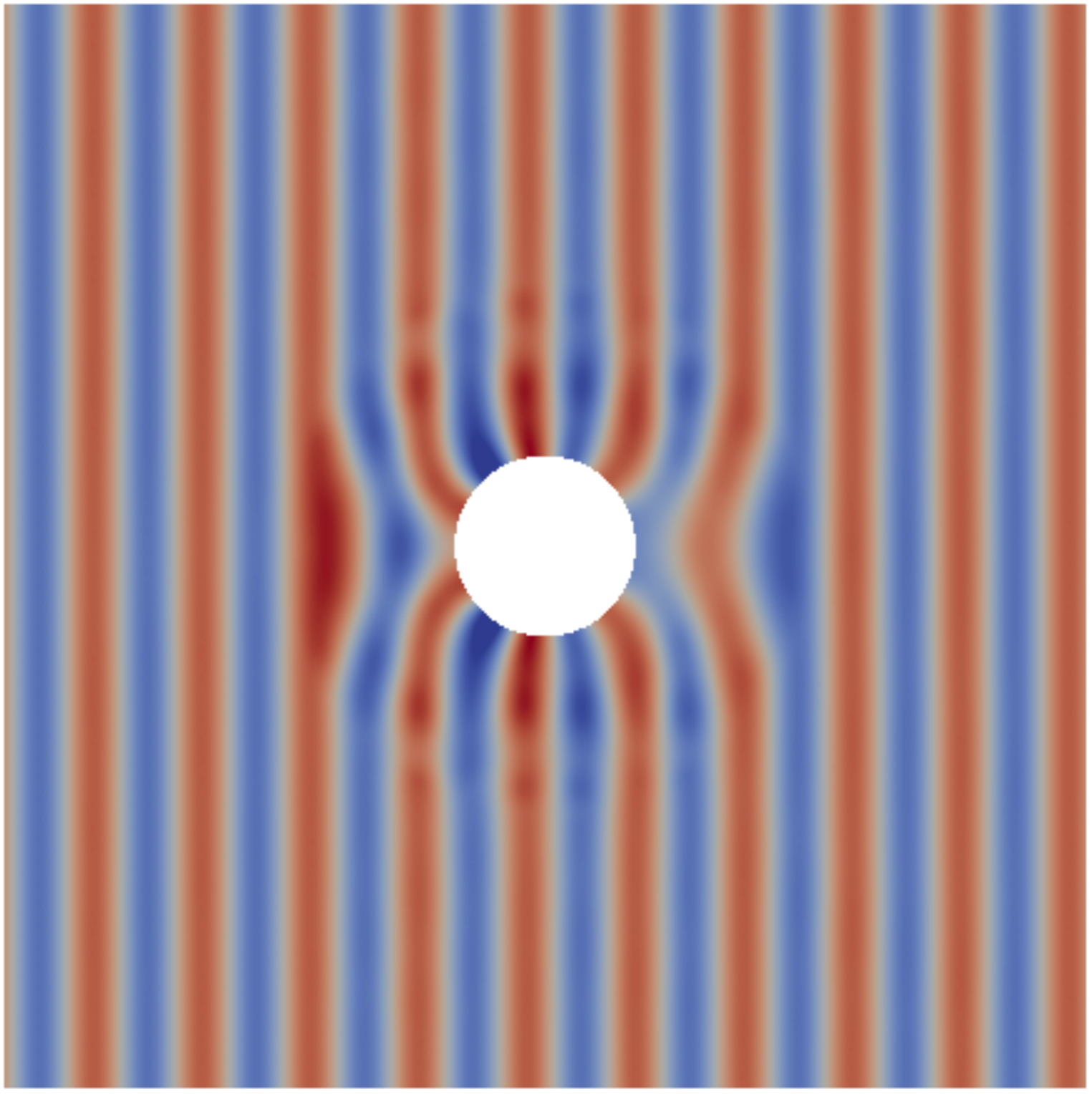}
		\includegraphics[scale=0.21]{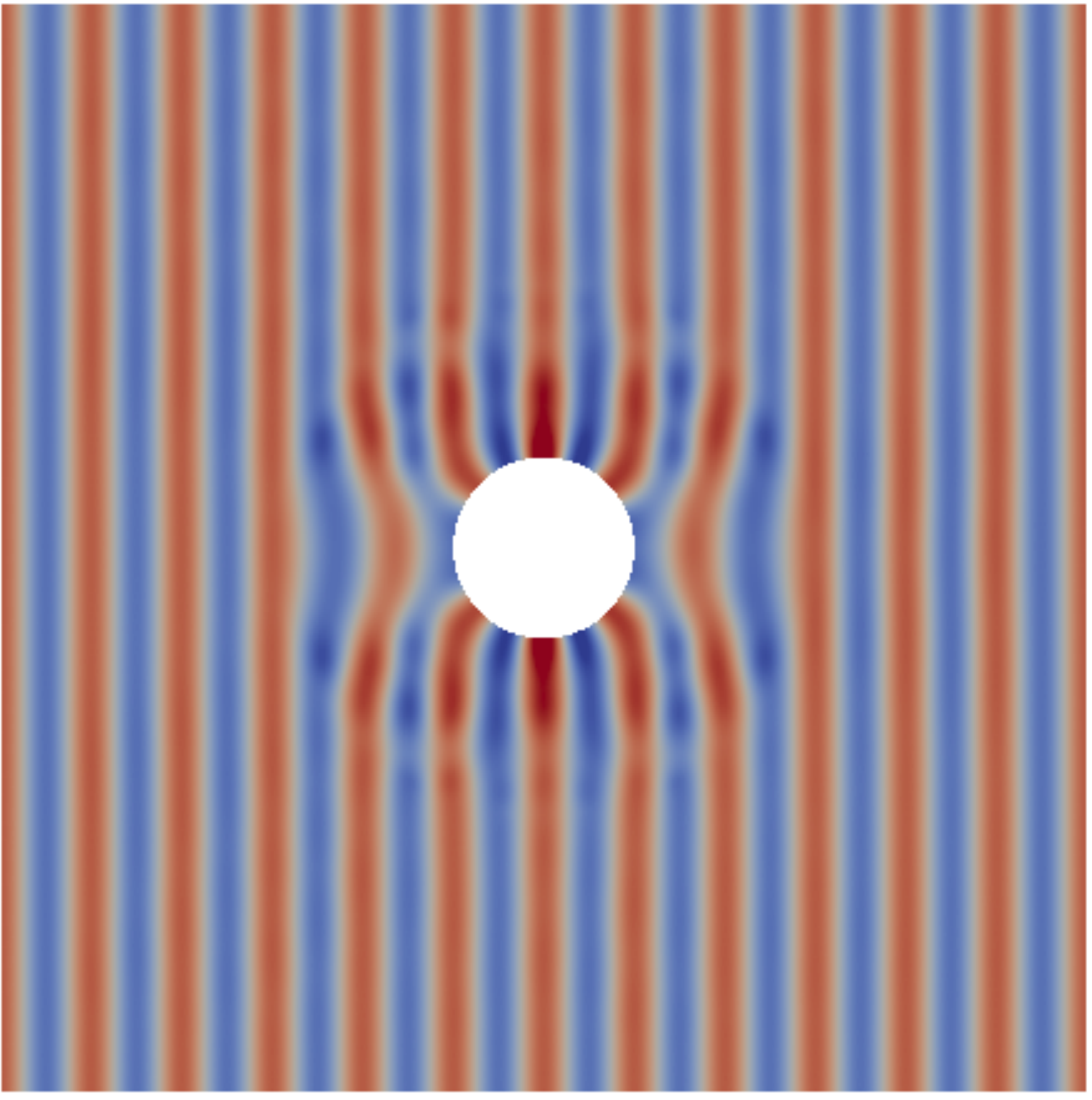}
	\end{center}
	\caption{The real part of the total wave fields without (top) and with (bottom) the cloak designed under uncertainty for the case of one direction and four frequencies.}\label{fig:4frequencies}
\end{figure}

\begin{figure}[!htb]
	\begin{center}
		\includegraphics[scale=0.28]{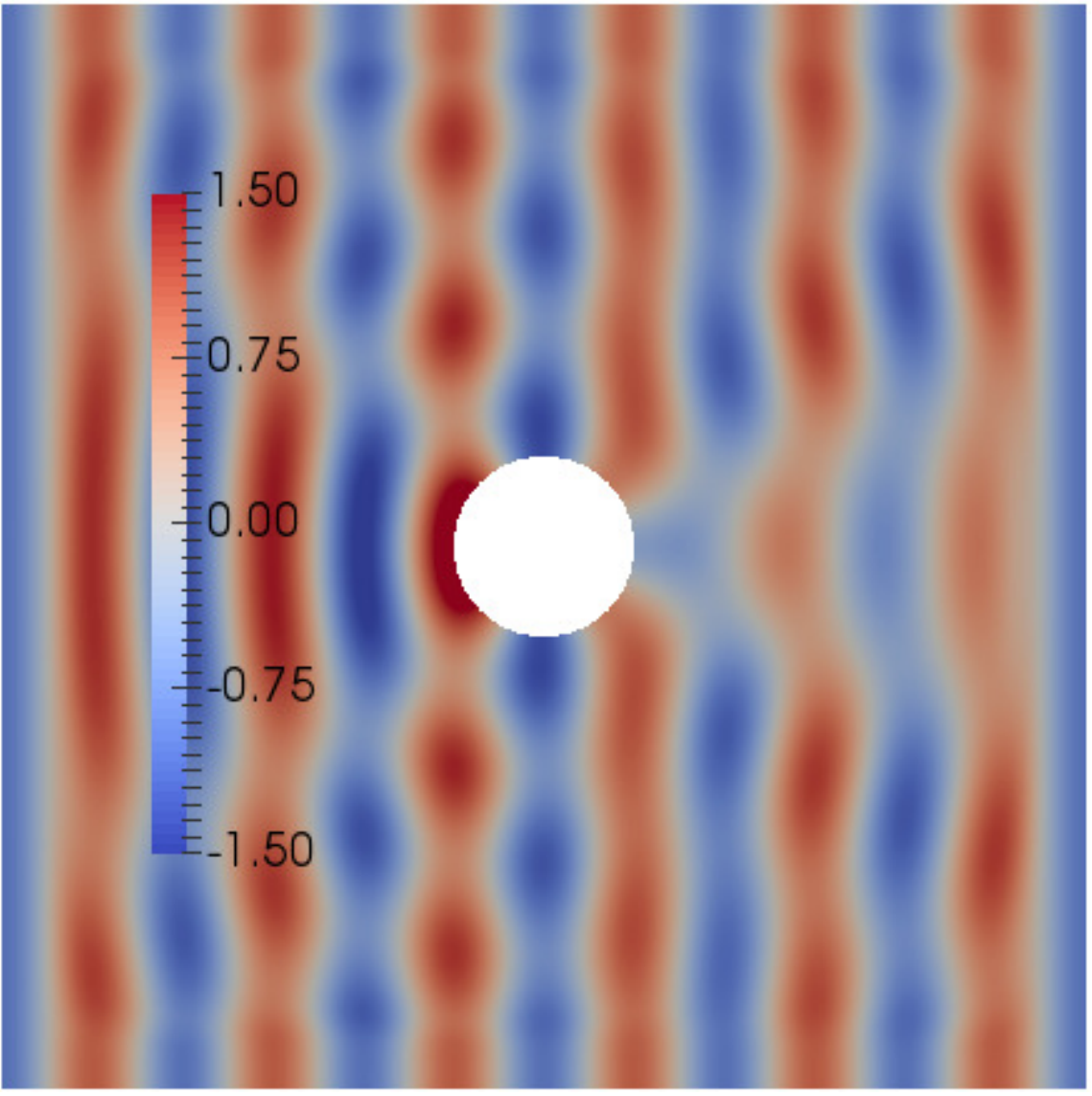}
		\includegraphics[scale=0.28]{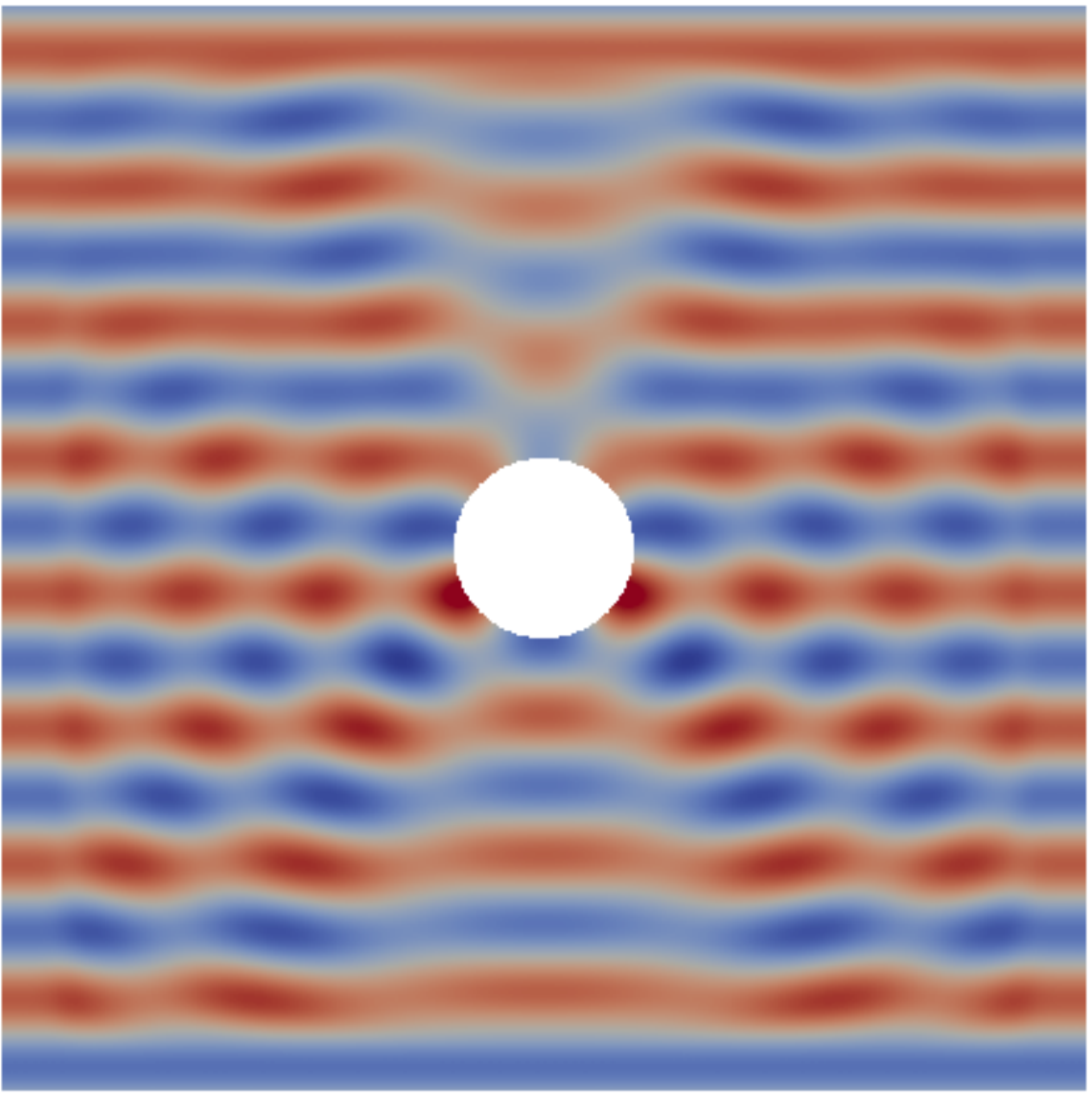}
		\includegraphics[scale=0.28]{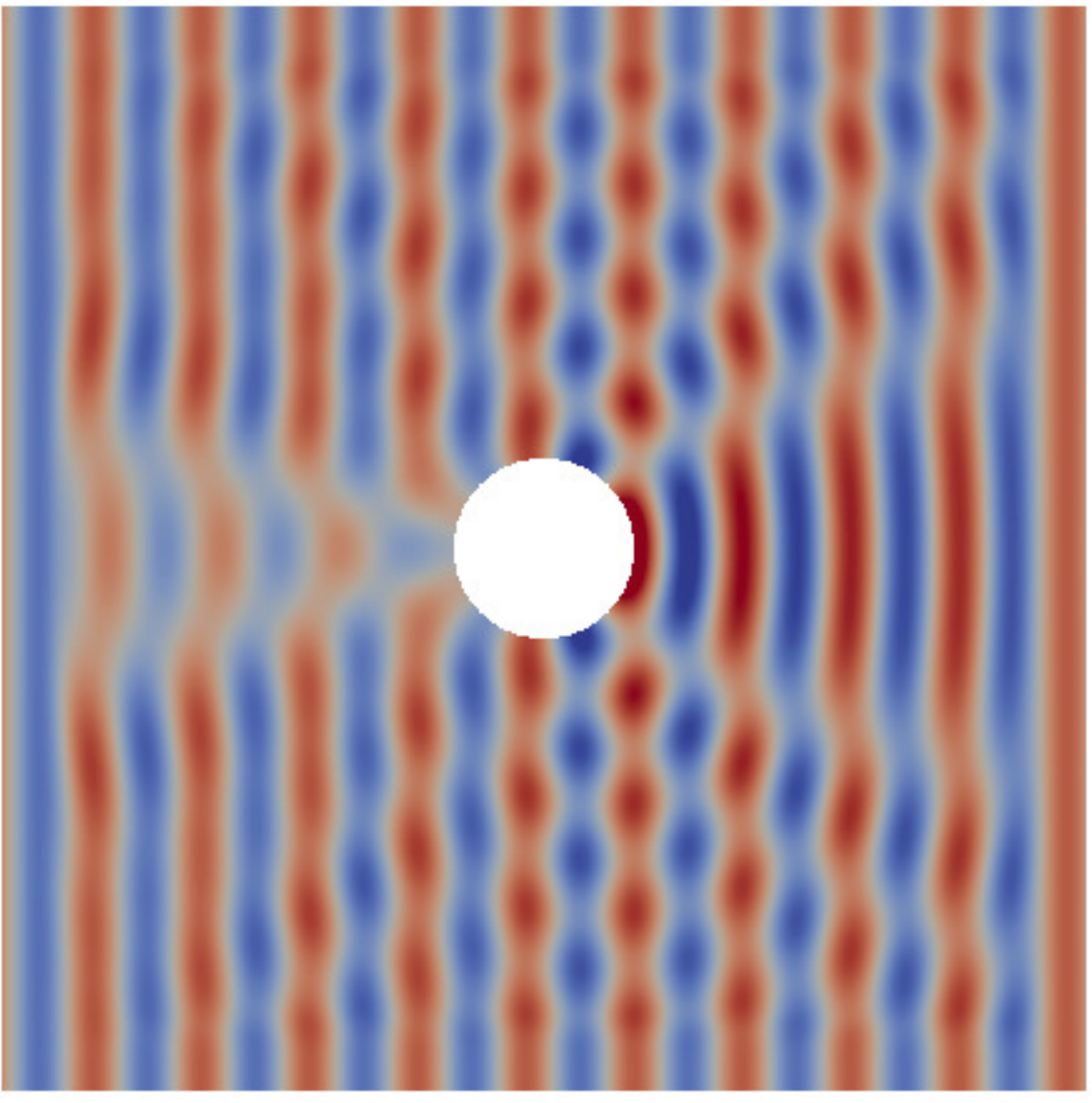}
		\includegraphics[scale=0.28]{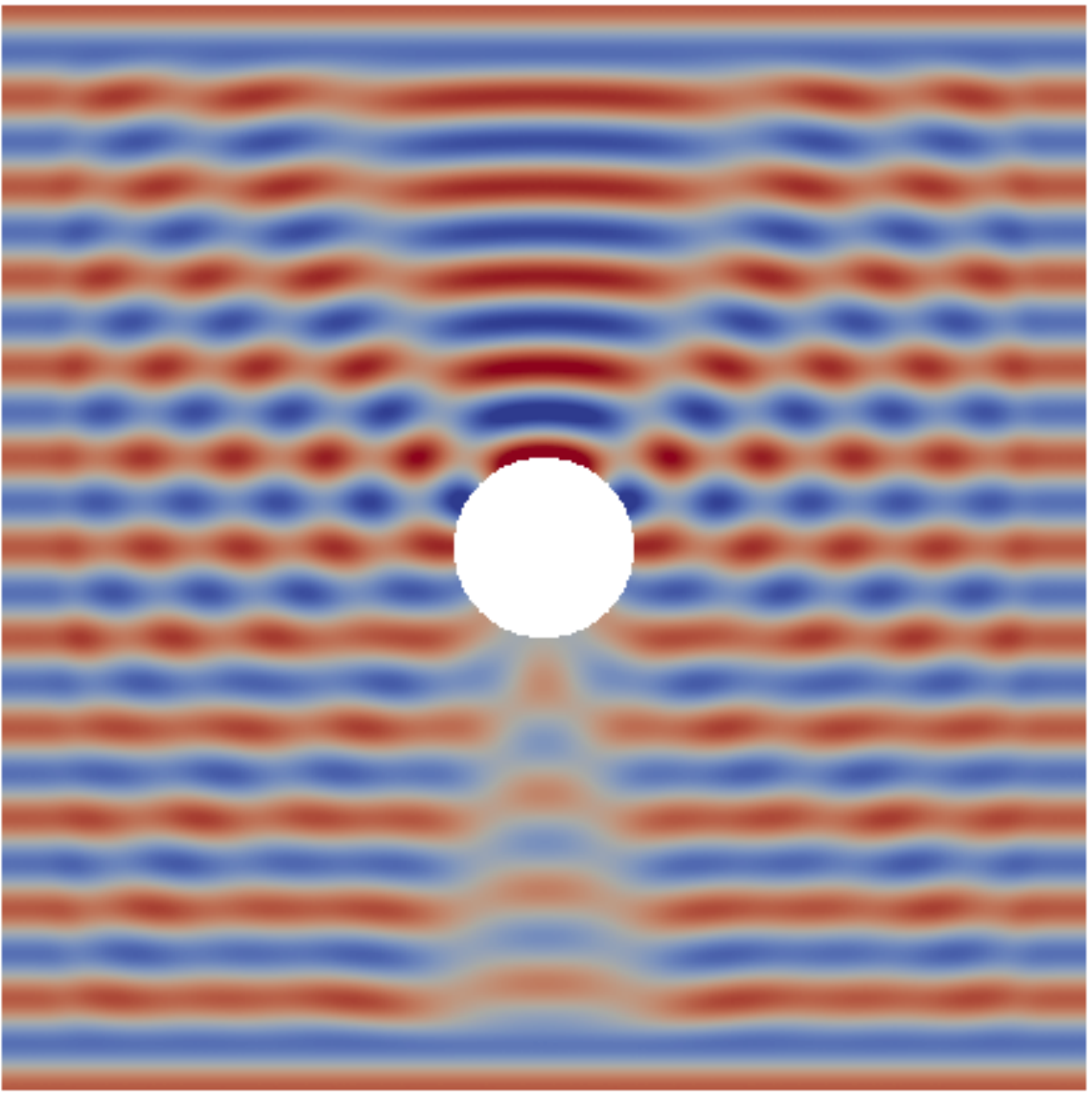}
		
		\includegraphics[scale=0.21]{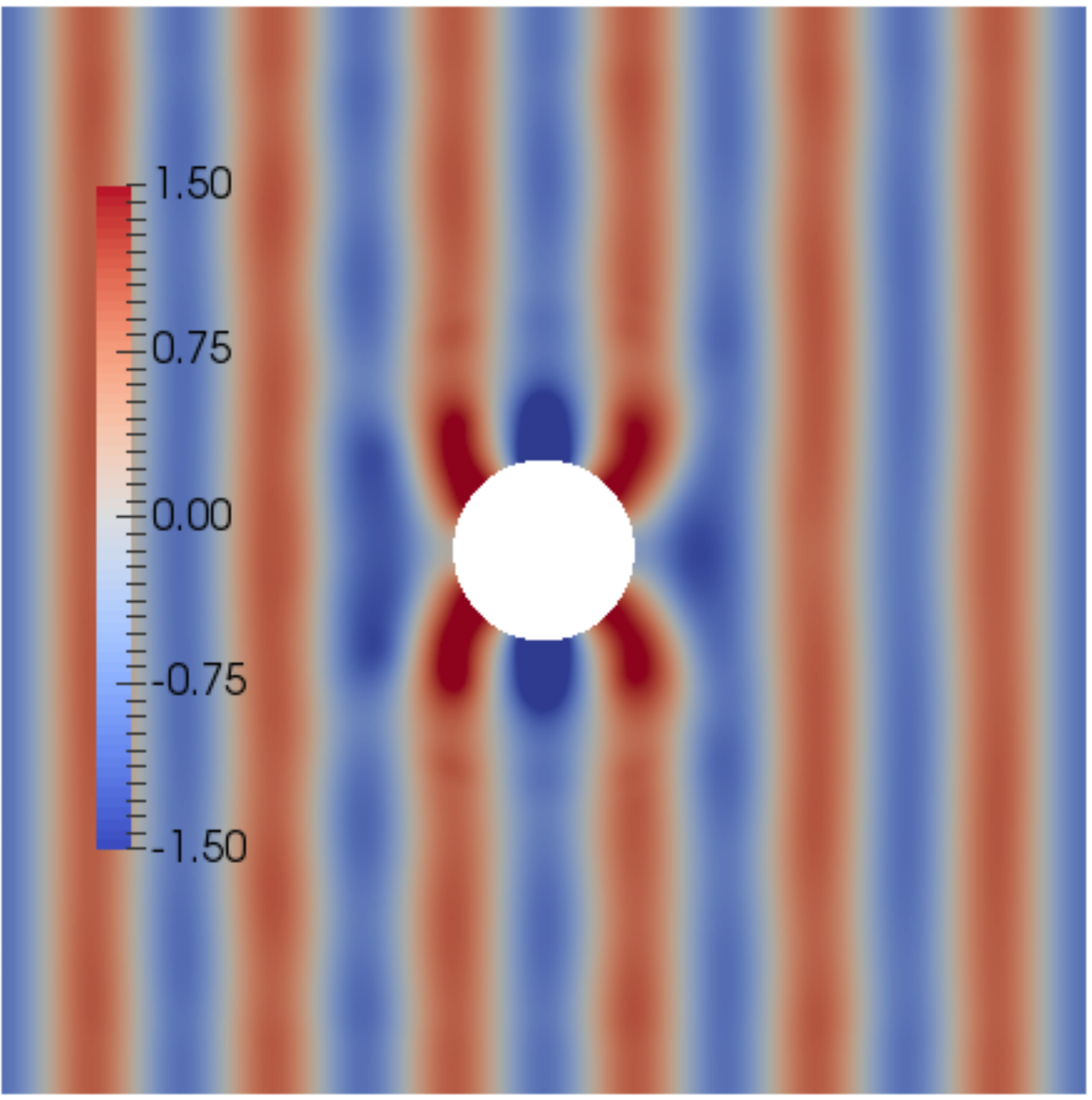}
		\includegraphics[scale=0.21]{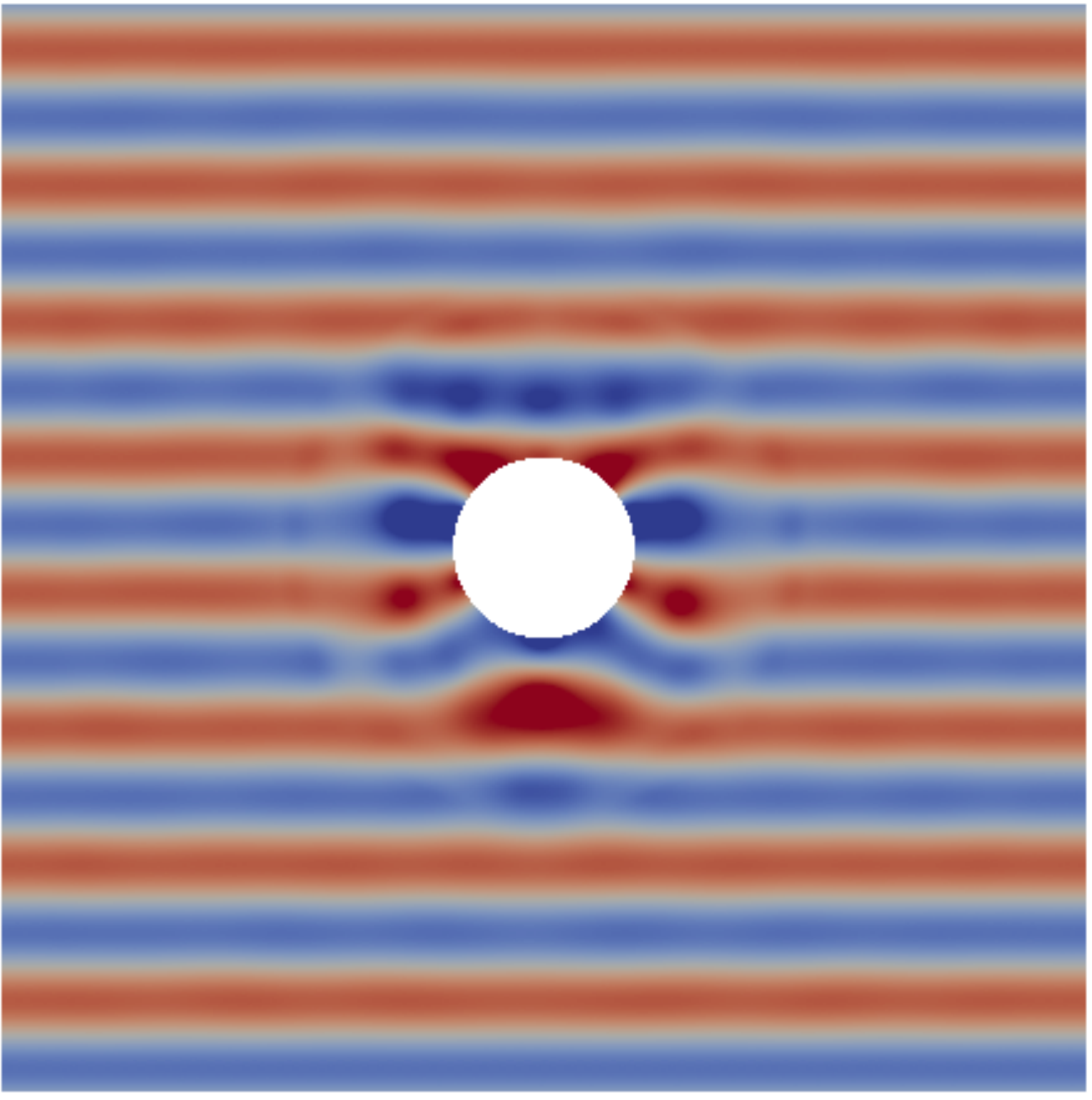}
		\includegraphics[scale=0.21]{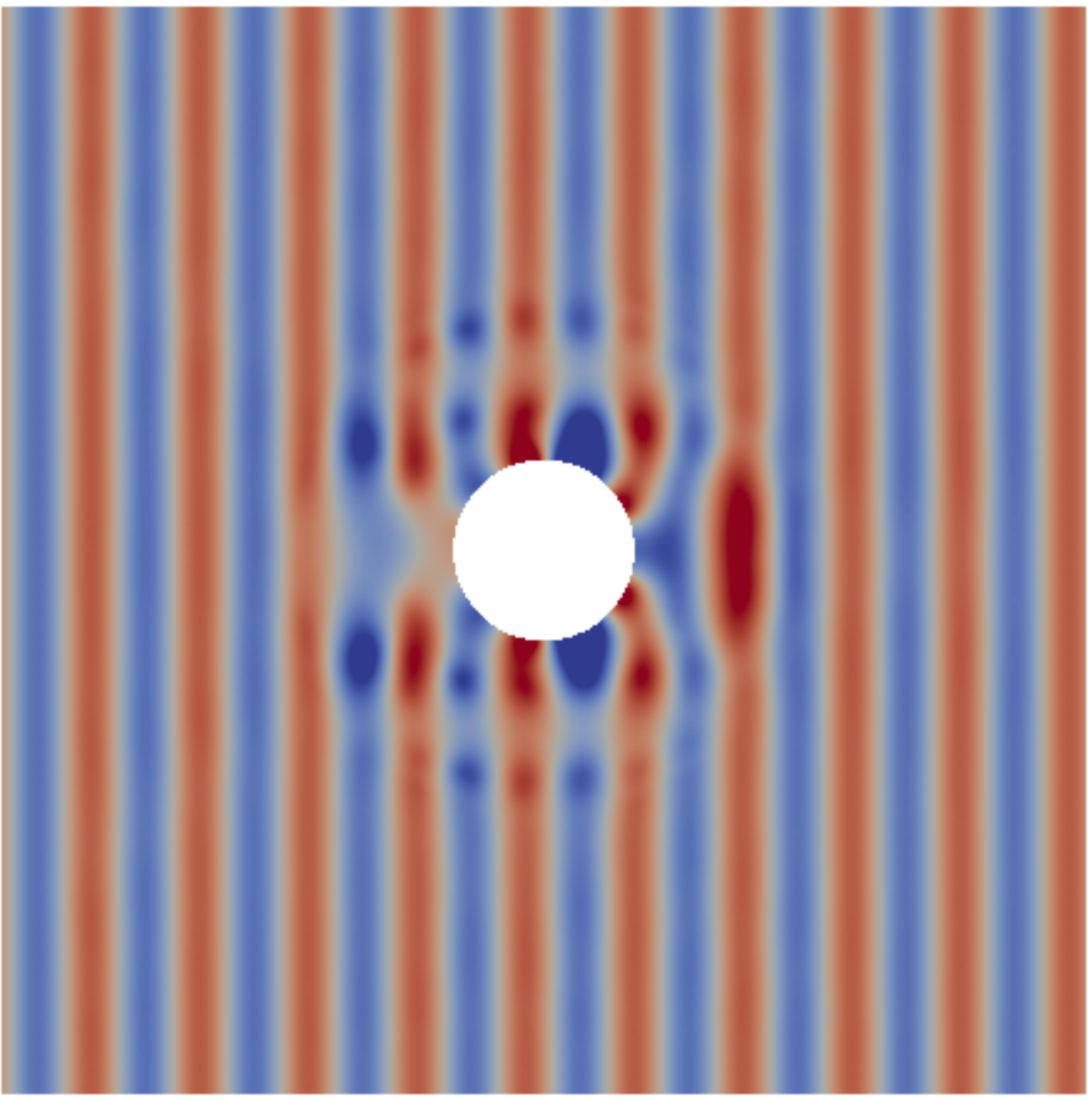}
		\includegraphics[scale=0.21]{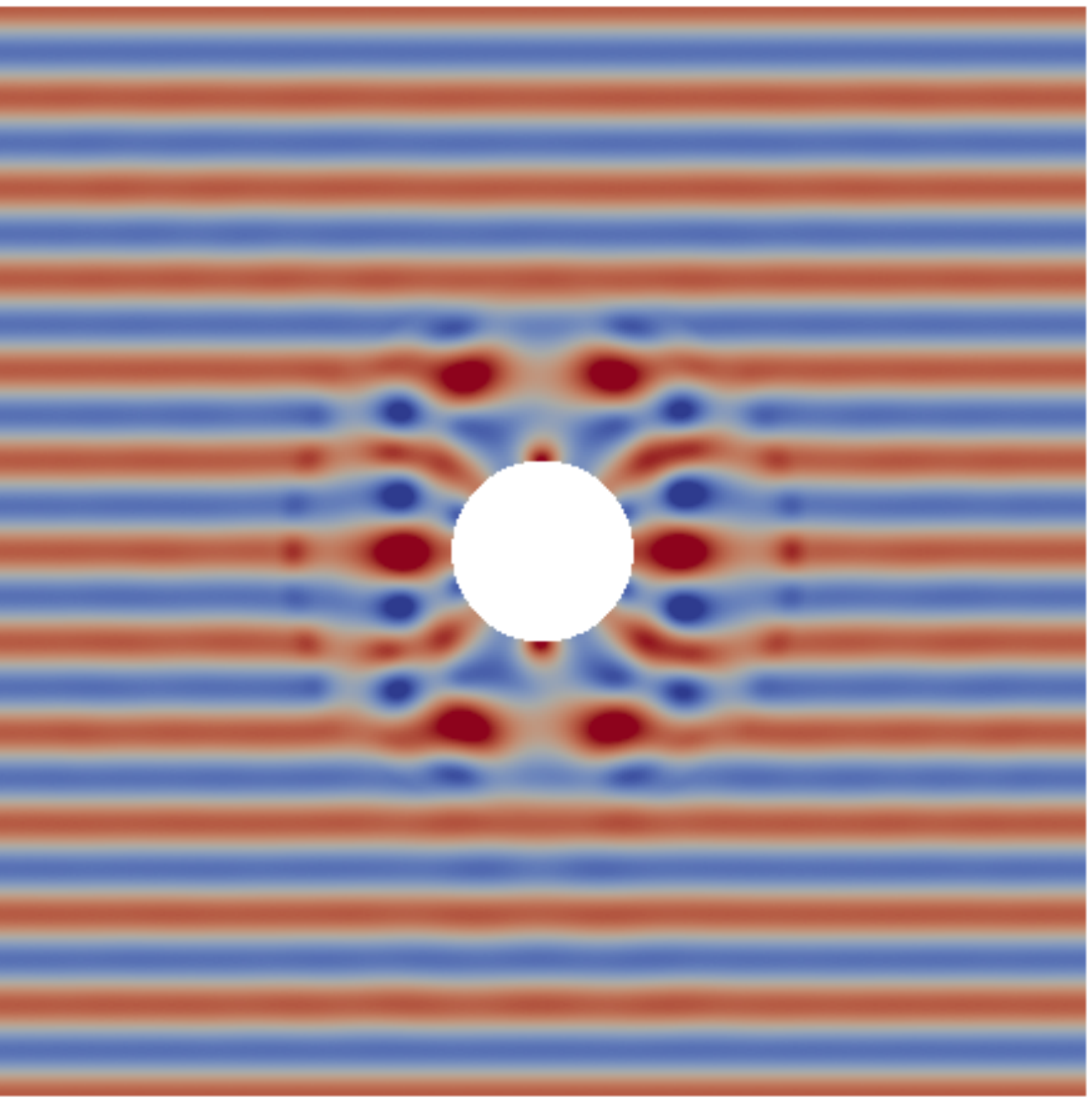}
	\end{center}
	\caption{The real part of the total wave fields without (top) and with (bottom) the cloak designed under uncertainty for the case of four directions and four frequencies.}\label{fig:4sources4frequencies}
\end{figure}

The optimal design under uncertainty using thequadratic approximation for the three different settings is shown in Fig.\ \ref{fig:4designs}, from which we can observe distinct patterns. The real parts of the total wave without and with the cloak are shown in Fig.\ \ref{fig:4sources} -- \ref{fig:4sources4frequencies}. We observe that the cloak can achieve effective cloaking for different directions with the same frequency, and can effectively reduce the scattering for different frequencies. This is expected as the characteristic length of the cloak has to accommodate all different wavelengths.

\subsection{Toward more complex geometry}

\begin{figure}
	\begin{center}
		\includegraphics[scale=0.13]{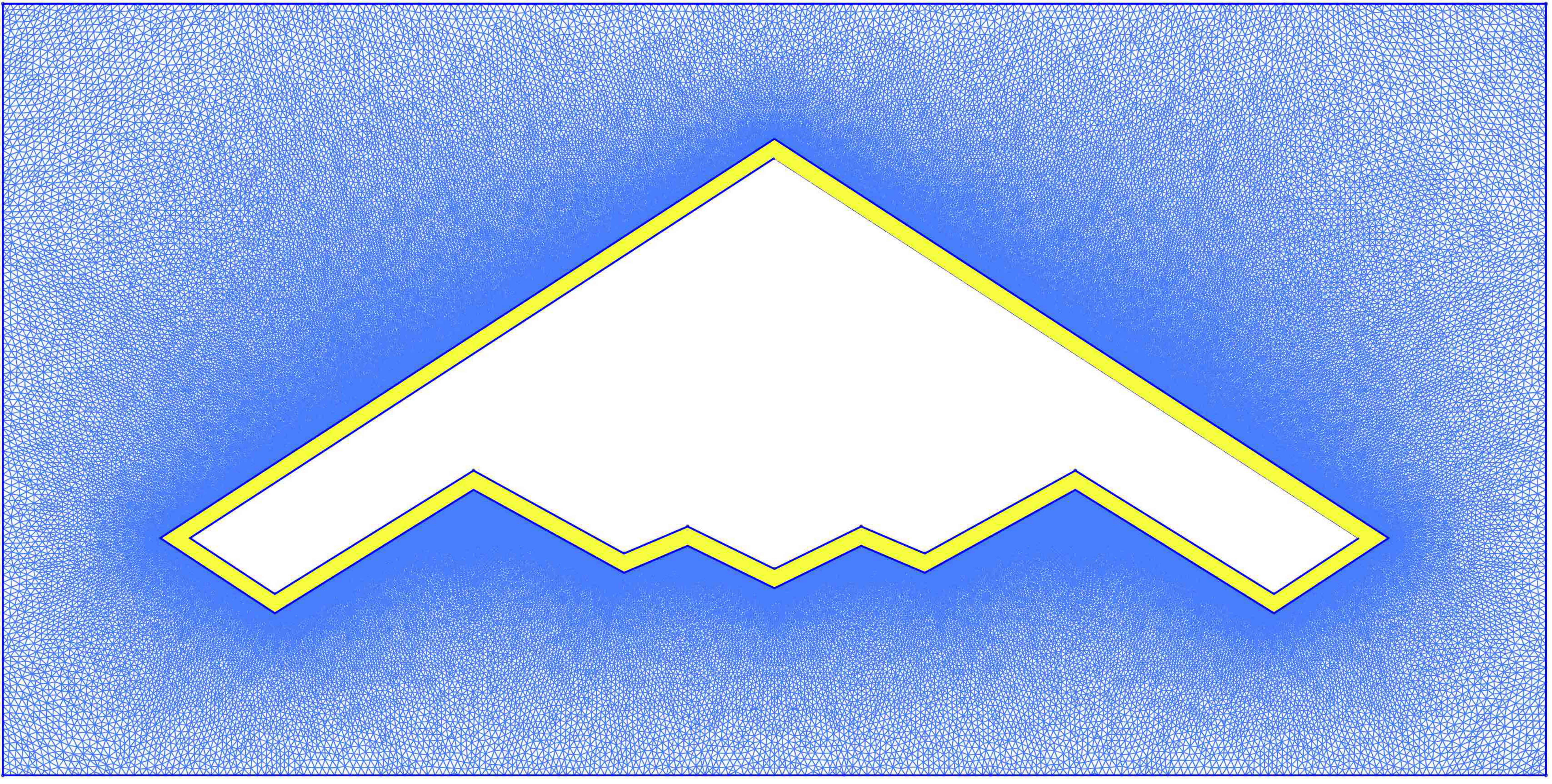}
		\includegraphics[scale=0.5]{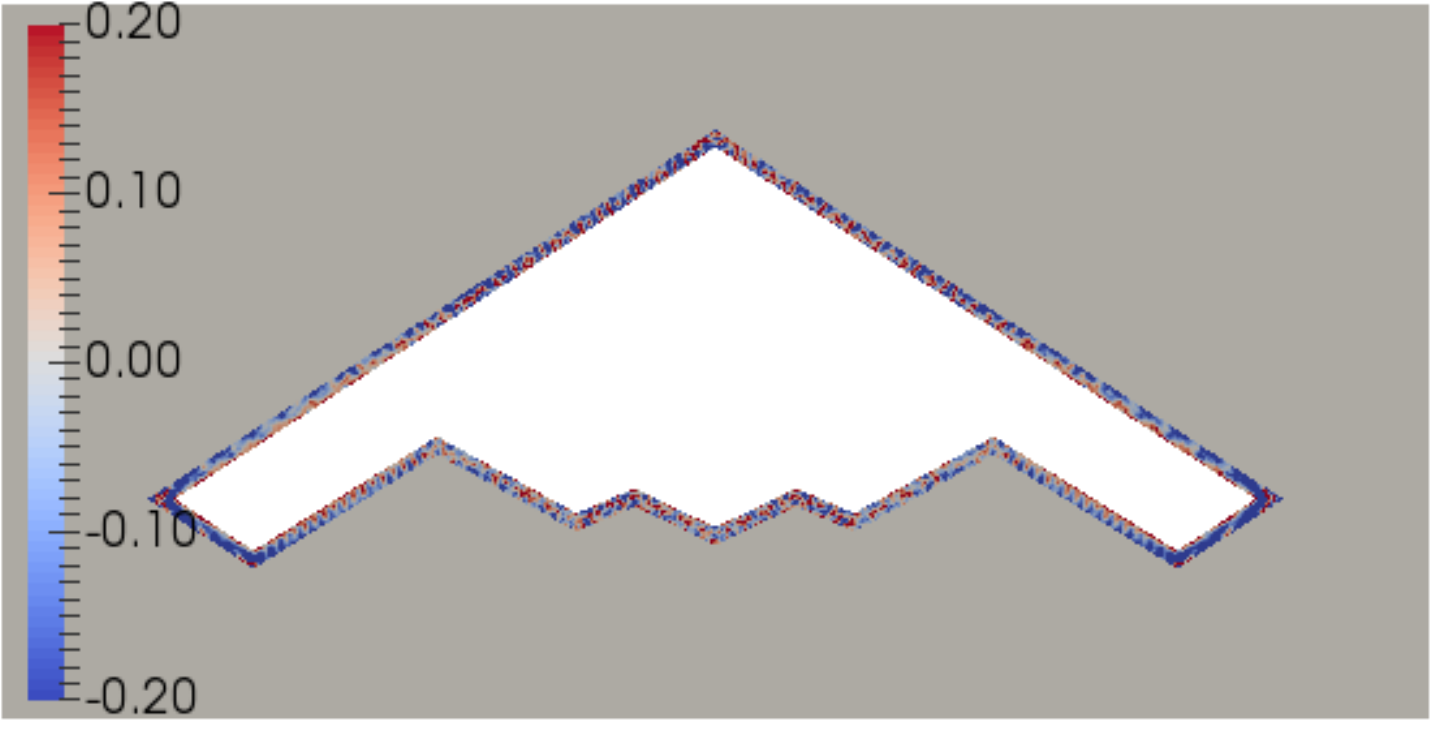}
	\end{center}
	\caption{Top: Geometry and adaptive mesh. Bottom: Optimal design field with deterministic approximation}\label{fig:b2spirit}
\end{figure}

\begin{figure}
	\begin{center}
		\includegraphics[scale=0.5]{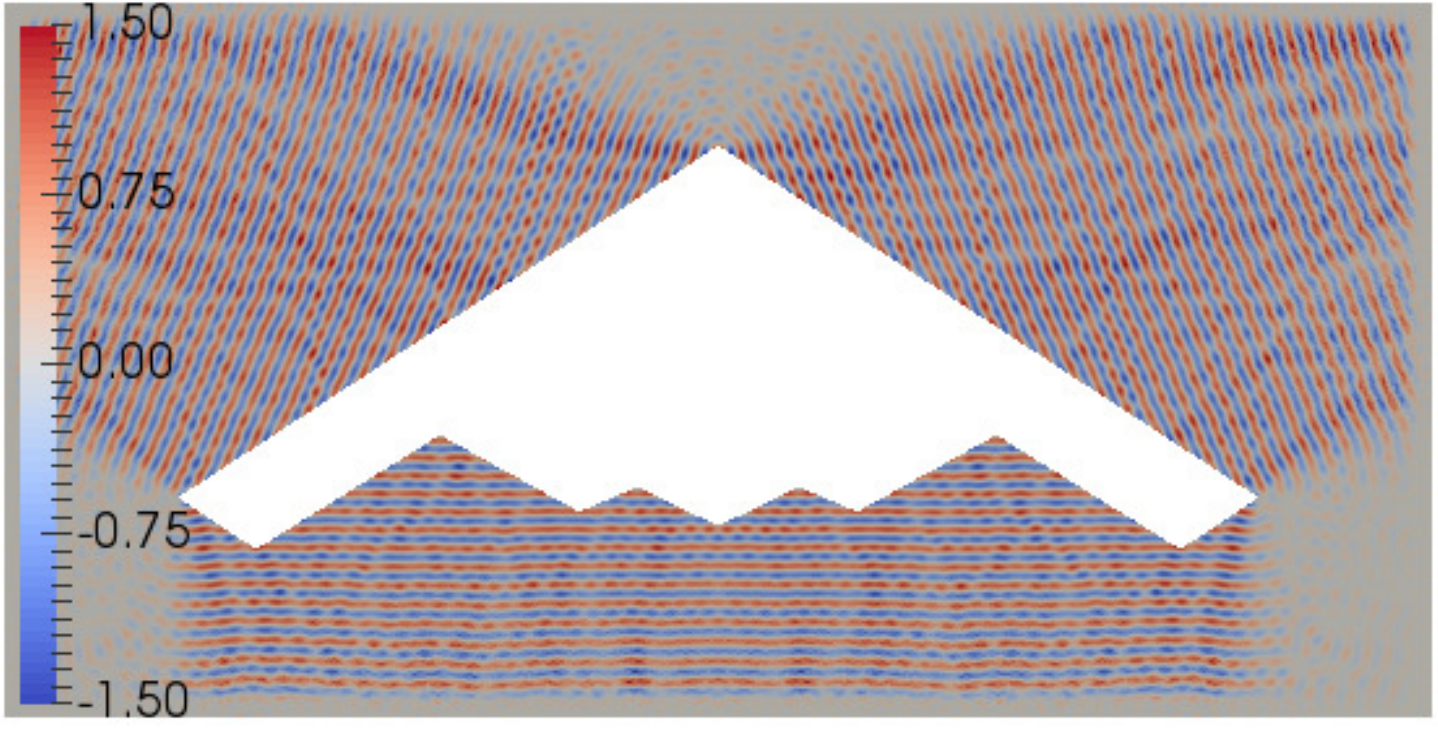}
		\includegraphics[scale=0.5]{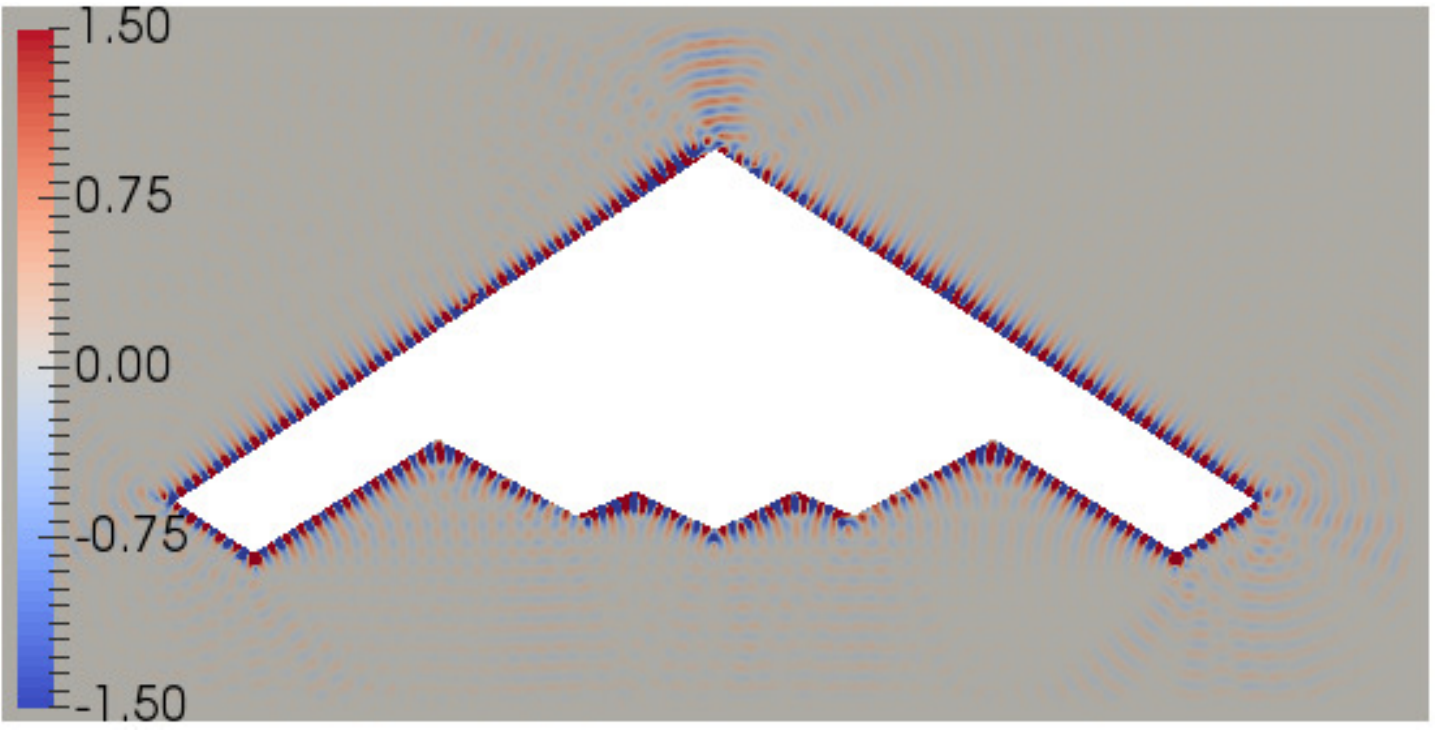}
	\end{center}
	\caption{The real part of the scattered wave fields without (top) and with (bottom)  the cloak.}\label{fig:b2scattered}
\end{figure}

Finally, we demonstrate the applicability of the proposed optimization
method for an obstacle with more complex geometry (notionally a
stealth aircraft) as shown in Fig.\ \ref{fig:b2spirit}. The design
variable field is discretized by a spatially-adapted mesh with 451,376
vertices and 898,136 elements, which results in DOF of 902,752 for the discrete state variable field with piecewise linear elements in the entire domain, 101,535 for the discrete uncertain variable field with piecewise linear elements in the thin cloaking layer (shown in
yellow), and 196,238 for the discrete optimization variable field with piecewise constant elements in the thin cloaking layer. We restrict ourselves to solution of the deterministic
optimal cloak problem, in order to demonstrate the feasibility of
computing the Hessian---which is a critical ingredient for the
approximate Newton method---for such a large problem.
Fig.\ \ref{fig:b2scattered} shows the large reduction in the scattered
wave field achieved after 200 iterations of the optimization
method. The reduction in the scattered field is striking, considering
the thinness of the cloak region and sharp corners. 



\section{Conclusions}
\label{sec:conclusion}

In this paper, we have developed a simulation-based optimal design
strategy for acoustic cloaks in the presence of material property
uncertainty. To the best of our knowledge this is the first work that
takes into account uncertainty in a systematic way for optimal design
of an acoustic cloak that is robust to material variability and
manufacturing error. Both the design variables and the uncertain
parameters are modeled by infinite-dimensional spatially-varying
fields, which become high-dimensional upon faithful discretization of
the optimal design problem. To tackle the curse of dimensionality in
the approximation of the uncertain parameter field, we employed a
scalable approximation method of the mean-variance objective based on
a Taylor expansion and a randomized SVD algorithm. To solve the
resulting high-dimensional optimization problem, we developed an
approximate Newton method in which the Hessian of the deterministic
approximation of the objective functional is used to provide an
effective approximation of the Hessian of the Taylor approximation of
the objective functional, motivated by the moderate uncertainty due to
material variability.

We demonstrated that the optimal design effectively eliminates the
scattered wave field from waves incident on simple circular
scatterers, not only for a single direction and single frequency, but
also for multiple-direction and multiple-frequency waves. We also
demonstrated that the deterministic optimization problem, on which the
approximate Hessian for the optimization under uncertainty problem is
based, can be tractably computed for an obstacle with complex
geometry. Moreover, we showed that the optimal design under
uncertainty performs better (lower variance in the scattered wave
field) in the case of random material properties than a deterministic
design does.

The proposed methodology is essentially scalable with respect to
increasing dimensions of design variables and uncertain parameters as
numerically evidenced by: the small and dimension-independent number
of forward Helmholtz solves needed to evaluate the Taylor-approximated
objective function; the weak dependence of the optimization
iterations on the problem dimension; and the dimension-independent
accuracy of the quadratic Taylor approximation.

Future research directions include (1) adding manufacturability
constraints on the design variable field stemming from additive
manufacturing processes;
(2) considering more complex three-dimensional problems with more
general objectives beyond cloaking; 
(3) developing and applying higher order Taylor
approximations (beyond quadratic) \cite{AlgerChenGhattas20} for the
objective functional for cases where large uncertainties arise; 
and (4) employing the Taylor approximations as control variates in a
variance reduction framework \cite{ChenVillaGhattas19}. 

\bibliographystyle{elsarticle-num}
\bibliography{references}

\end{document}